\RequirePackage{fix-cm}
\documentclass{svjour3}                     
\smartqed  
\usepackage[margin=1in]{geometry}
\usepackage{lineno,hyperref}
\modulolinenumbers[5]
\usepackage{cuted}
\usepackage[ansinew]{inputenc}
\usepackage[english]{babel}
\usepackage{fancyhdr}
\usepackage{anyfontsize}
\usepackage[active]{srcltx}
\usepackage{amssymb,amsmath,amsfonts,bm}
\usepackage{color}
\usepackage{array}
\usepackage{textcomp}
\usepackage[linesnumbered,ruled,vlined]{algorithm2e}
\usepackage{pgfplots}
\pgfplotsset{compat=1.17}
\renewcommand{\v}[1]{\mathbf{#1}}
\newcommand{\fun}[1]{\bm{\mathit{#1}}}

\newcommand{\p}[2][p]{{#2^{(#1)}}}
\renewcommand{\P}[2][p]{{\left[#2\right]_{#1}}}

\newcommand{\M}{\v{M}}
\newcommand{\K}{\v{K}}
\newcommand{\C}{\v{C}}
\newcommand{\G}{\fun{G}}
\renewcommand{\H}{\fun{H}}

\newcommand{\phiv}{\bm{\phi}}
\newcommand{\FU}{\bm{\mu}}
\newcommand{\FV}{\bm{\nu}}

\newcommand{\FG}{\v{\check{G}}}
\newcommand{\FH}{\v{\check{H}}}

\newcommand{\U}{\v{U}}
\newcommand{\V}{\v{V}}
\newcommand{\Ut}{\v{\dot{U}}}
\newcommand{\Vt}{\v{\dot{V}}}

\newcommand{\WU}{\v{\Psi}}
\newcommand{\WV}{\v{\Upsilon}}

\newcommand{\WUfun}{\fun{\Psi}(\v{z})}
\newcommand{\WVfun}{\fun{\Upsilon}(\v{z})}
\newcommand{\ffun}{\fun{f}(\v{z})}

\newcommand{\WUfunsim}{\fun{{\Psi}}}
\newcommand{\WVfunsim}{\fun{{\Upsilon}}}
\newcommand{\ffunsim}{\fun{{f}}}

\newcommand{\0}{\v{0}}
\newcommand{\RHS}{\v{F}}
\newcommand{\Id}{\v{I}}

\newcommand{\Is}{\mathcal{I}}
\newcommand{\Rs}{\mathcal{R}}

\DeclareMathAlphabet{\mymathbb}{U}{BOONDOX-ds}{m}{n}
\newcommand{\iu}{\mymathbb{i}}
\renewcommand{\Im}[1]{{\text{Im}[#1]}}
\renewcommand{\Re}[1]{{\text{Re}[#1]}}

\definecolor{myred}{rgb}{0.9,0,0}
\definecolor{Color_FS}{rgb}{.7,0,.7}
\definecolor{Color_NF}{rgb}{0,0,.7}
\definecolor{Color_NFTO}{rgb}{0,.5,.2}
\definecolor{Color_MD}{rgb}{.75,.35,0}
\definecolor{Color_SMD}{rgb}{1,.7,0}

\newcommand{\red}[1]{{\leavevmode\color{red!0!black}#1}}

\newcommand{\dblnk}{\begin{minipage}{12pt}\fontsize{6pt}{3pt}\selectfont{$>1$\\$<p$}\end{minipage}}

\newcommand{\Ir}{\mathrm{I}}
\newcommand{\ur}{u}
\newcommand{\vr}{v}
\newcommand{\ar}{a}
\DeclareMathOperator{\e}{e}
\usepackage{subfiles}

\journalname{Nonlinear Dynamics}
\begin{document}
\title{High order direct parametrisation of invariant manifolds for model order reduction of finite element structures: application to large amplitude vibrations and uncovering of a folding point}
\titlerunning{High order direct parametrisation of invariant manifolds for model order reduction of finite element structures}
\author{Alessandra~Vizzaccaro$^{1,3}$ \and Andrea~Opreni$^{2}$ \and Lo\"ic~Salles$^{3}$ \and Attilio~Frangi$^{2}$ \and Cyril~Touz\'e$^{4}$}

\institute{
    $^{1}$ Department of Engineering Mathematics\\
    University of Bristol\\
    Beacon House, Queens Road, Bristol BS8 1QU\\
    \\
	$^{2}$ Department of Civil and Environmental Engineering\\
	Politecnico di Milano\\
	Piazza Leonardo da Vinci, 32, 20133 Milano MI\\
	\\
	$^{3}$  Department of Mechanical Engineering\\
	Imperial College London\\
	Exhibition Rd, South Kensington, London SW7 2BX\\
	\email{a.vizzaccaro17@imperial.ac.uk}\\
	\url{https://orcid.org/0000-0002-2040-4753}\\
	\\
    $^{4}$ Institute of Mechanical Sciences and Industrial Applications (IMSIA)\\ 
    ENSTA Paris - CNRS - EDF - CEA - Institut Polytechnique de Paris\\
    828 boulevard des mar\'echaux 91762 Palaiseau cedex
}

\date{Received: date / Accepted: date}

\maketitle

\vspace{1cm}

\begin{abstract}
This paper investigates model-order reduction methods for geometrically nonlinear structures. The parametrisation method of invariant manifolds is used and adapted to the case of mechanical systems expressed in the physical basis, so that the technique is directly applicable to problems discretised by the finite element method. Two nonlinear mappings, respectively related to displacement and velocity, are introduced, and the link between the two is made explicit at arbitrary order of expansion. The same development is performed on the reduced-order dynamics which is computed at generic order following the different styles of parametrisation. More specifically, three different styles are introduced and commented: the graph style, the complex normal form style and the real normal form style. These developments allow making better connections with earlier works using these parametrisation methods. The technique is then applied to three different examples. A clamped-clamped arch with increasing curvature is first used to show an example of a system with a softening behaviour turning to hardening at larger amplitudes, which can be replicated with a single mode reduction. Secondly, the case of a cantilever beam is investigated. It is shown that the invariant manifold of the first mode shows a folding point at large amplitudes which is not connected to an internal resonance. This exemplifies the failure of the graph style due to the folding point, whereas the normal form style is able to pass over the folding. Finally, A MEMS micromirror undergoing large rotations is used to show the importance of using high-order expansions on an industrial example.
\end{abstract}


\section{Introduction}

This work is concerned with model-order reduction techniques for nonlinear vibrations of structures featuring geometric nonlinearity, with a particular emphasis on problems using the finite element (FE) procedure as space discretisation method. In this context, numerous methods have been proposed in the past in the FE community: stiffness evaluation procedure (STEP)~\cite{Mignolet08,KIM2013,mignolet13,Perez2014}, implicit condensation~\cite{Hollkamp2005,Hollkamp2008,FRANGI2019,NicolaidouIceKE}, E-STEP~\cite{KimEstep} and M-STEP~\cite{Vizza3d,givois21-CS}, modal derivatives (MD)~\cite{IDELSOHN1985,IDELSOHN1985b,Weeger2016} and quadratic manifold built from modal derivatives~\cite{Jain2017,Rutzmoser}. 

On the other hand, reduction methods for geometrically nonlinear systems have also been studied in the dynamical systems community, leading to important theoretical developments with methods which were mostly applied to partial differential equations (PDEs), and not to FE problems with large dimensions. In this direction, important contributions led to the definition of Nonlinear Normal Modes (NNMs) as invariant manifolds of the system, tangent to the linear eigenspaces~\cite{ShawPierre91,ShawPierre93}. As emphasised in numerous papers, the invariance property is key in order to derive accurate ROMs, for the simple reason that reduction to a non-invariant set leads to simulate trajectories with the reduced models that do not exist for the full system, hence immediately questioning the validity of the ROM. The idea has then been pushed forward, using either computational methods for the solution phase~\cite{PesheckJSV}, or a different methodology for the theoretical settings, i.e.\ by using the normal form approach~\cite{touze03-NNM,TOUZE:JSV:2006,TouzeCISM}.

Some steps have been recently taken 
in order to compare and assess the methods developed in the FE community against those relying on invariant manifold theory. In particular, it has been clearly demonstrated that most of the methods such as implicit condensation or quadratic manifold with MD, need a slow/fast assumption in order to deliver accurate predictions~\cite{HallerSF,VERASZTO,VizzaMDNNM,YichangICE,YichangVib}. By slow/fast assumption, it is meant that a clear frequency gap between the eigenfrequencies of the slave and of the master modes, needs to be fulfilled.

In the mathematical community, a major advancement in the understanding and formalisation of the reduction to invariant manifolds has been made thanks to the parametrisation method, first introduced by Cabr\'e, Fontich and de la Llave~\cite{Cabre1,Cabre2,Cabre3}, and then rewritten in a more computational framework, easier to understand for engineering applications, in the book by Haro {\em et al.}~\cite{Haro}. This important formalisation allows unifying different developments in the same framework. While previous works relied either on the invariant manifold computation proposed {\em e.g.} by~\cite{Carr,gucken83} (assuming a functional relationship between slave and master coordinates), or on the normal form theory~\cite{Jezequel91,touze03-NNM}, the parametrisation method allows one to show that both solutions can be derived from the {\em invariance equation}, which can be solved either with a {\em graph style} or a {\em  normal form style}.

The parametrisation method has then been first adapted to the case of vibratory systems by Haller and Ponsioen~\cite{Haller2016}. Also, whereas most of the previous studies on NNMs took advantage of existence and uniqueness of Lyapunov subcentre manifolds (LSM)~\cite{Lyapunov1907,Kelley2} to settle down the definitions in a correct mathematical framework, the situation for dissipative systems were less clear, as underlined by different investigations~\cite{NeildNF01,CIRILLO2016}. One of the main contribution of Haller and Ponsioen has thus also been to provide existence and uniqueness theorems for such invariant manifolds defined as spectral submanifolds (SSMs). In the damped case, the smoothest nonlinear continuation of a spectral subspace of the linearised system is the SSM, and is unique under general persistence and non-resonance conditions provided in~\cite{Haller2016,PONSIOEN2018}.   The link between the conservative case, with LSMs densely filled with periodic orbits, and the dissipative case, has been further investigated in~\cite{Llave2019}. Elaborating on the parametrisation method, reduction methods up to arbitrary order for two-dimensional manifolds with damping included (SSM) have then been proposed in~\cite{PONSIOEN2018}. 

One important drawback of the methods using invariant manifolds with regard to applications to large FE models was their need to express the equations of motion in the modal basis as a starting point. However, recent developments tackled this limitation and proposed direct computations in order to pass from the physical space to ROMs expressed with coordinates linked to the invariant manifolds. Elaborating on previous results on normal forms, a direct approach has been proposed in~\cite{artDNF2020} and further developed  in~\cite{AndreaROM,YichangVib}, allowing one to express the reduced dynamics with normal coordinates in an invariant-based span of the phase space. Leveraging on SSM, a direct approach has also been proposed in~\cite{JAIN2021How,MingwuLi2021_1,MingwuLi2021_2}, taking into account the damping and proposing arbitrary order approximations in an automated framework.

An overview of the nonlinear reduction methods has been proposed in~\cite{ROMGEOMNL}, allowing one to put all the developments using nonlinear techniques in perspective, and with a special emphasis on applications to FE models. In particular, SSM as defined in~\cite{Haller2016} are unique only when the order of the asymptotic development reaches the spectral quotient defined as the ratio between the maximal damping ratio of the slave modes divided by the smallest of the masters. In practice, for large FE models, this number can be very large such that all asymptotic developments are just approximations of the unique SSM. Consequently, prior developments led in~\cite{ShawPierre91,ShawPierre93,TOUZE:JSV:2006} with damping included were low-order approximations of the SSM, either with a graph style or a normal form style. 
\red{
Along the same lines, and as remarked in~\cite{ROMGEOMNL}, the computations proposed in this contribution
as well as those shown for example in~\cite{JAIN2021How,MingwuLi2021_1,MingwuLi2021_2}, are approximations of the unique SSM, which is reached at a very high order only. Importantly, all lower order approximations
of the SSM share the invariance property, up to the selected order, and can be thus used safely to provide accurate ROMs.}

This paper elaborates on the previous analysis led in~\cite{artDNF2020,AndreaROM}, with the aim of pushing the developments further  to propose an arbitrary order expansion. As a main difference, the parametrisation method~\cite{Cabre3,Haro} is used instead of starting from the normal form transformation. In short, whereas normal form expansion, as proposed in~\cite{touze03-NNM,TOUZE:JSV:2006,artDNF2020,AndreaROM}, first computes the complete nonlinear mapping and then reduces by selecting a few master normal coordinates, the parametrisation method first reduces by selecting the master coordinates, and then computes the expansions, with the added value that different solutions are possible, thus offering the possibility of using either a graph style or a normal form style. With this initial choice, the developments are thus closer to those already reported in~\cite{JAIN2021How,MingwuLi2021_1,MingwuLi2021_2}, where arbitrary order expansions have already been shown, together with the possibility of using either graph or normal form style. The main differences can be listed as follows: \red{(i) the focus here is on large FE models of mechanical systems for which the damping matrix is diagonalised by the eigenvectors of the conservative system; (ii) thanks to this assumption, displacement and velocity mappings can be treated separately allowing to show the relationship between the two at generic order and to retrieve homological equations in the sole displacement mapping; (iii) a number of implementation details on the treatment of the direct computation are reported in order to decrease the computational burden (e.g. treatment of the nonlinear tensors to reduce the memory consumption, derivation of homological equations in the sole displacement to halve the size of the linear systems to solve); (iv) two different versions of the normal form style are investigated: a complex and a real normal form style, thus pushing further the developments on realification and complexification~\cite{Haro}}.

Thanks to the computational developments, specific applications are then reported to underline the quality of the ROMs obtained. First, a clamped-clamped arch with increasing curvature is investigated in order to demonstrate that higher-order expansions are able to capture a behaviour in the backbone curve that is first softening then hardening, with a single mode reduction. Then, the fundamental mode of a cantilever beam is studied, 
putting in evidence a folding of the invariant manifold which is not due to an internal resonance, a case that had not been reported before. Due to this very particular behaviour, it is then shown that the graph style is not able to provide a correct ROM up to very large displacements. On the other hand, normal form style passes through the folding point and allows obtaining accurate results. Finally, a MEMS (Micro-Electro-Mechanical System) micromirror is used to demonstrate how the method can handle large FE structures of interest for industrial applications.

\section{Equations of motion and parametrisation method}

\subsection{Equations of motion and eigenproblem}\label{sec:EOMeigenpb}

We consider large-amplitude, geometrically nonlinear vibrations of an elastic mechanical structure discretised by the finite element method. It is assumed that the only nonlinearity comes from the strain-displacement relationship, while the constitutive law is linear elastic. In this framework, the equations of motion contain quadratic and cubic nonlinearities, and can be written in a general formulation as~\cite{holzapfel00,LazarusThomas2012,Touze:compmech:2014,Vizza3d}
\begin{equation}\label{eq:FEstart}
\M \ddot{\U} + \C \Ut + \K \U +\G(\U,\U) + \H(\U,\U,\U) = \0,
\end{equation}
where $\U$ is the $N$-dimensional time-dependent displacement vector, gathering all the degrees of freedom of the model, $\M$ and $\K$ are respectively the mass and stiffness matrix, $\C$ stands for the damping matrix. Quadratic and cubic polynomial nonlinearities are expressed through the terms $\G(\U,\U)$ and $\H(\U,\U,\U)$ which can be written as
\begin{subequations}
\begin{align}
\G(\U,\U)&=\sum^N_{r=1}\sum^N_{s=1}\G_{rs} U_r U_s,\\
\H(\U,\U,\U)&=\sum^N_{r=1}\sum^N_{s=1}\sum^N_{t=1}\H_{rst} U_r U_s U_t,
\label{eq:tensor_product_phys}
\end{align}
\end{subequations}
where $\G_{rs}$ stands for the $N$-dimensional vector of coefficients $G^p_{rs}$, for $p=1,\, ...,\, N$, and similarly $\H_{rst}$ is a vector of coefficients $H^p_{rst}$.

The eigenproblem of the corresponding conservative linear system reads
\begin{equation}
(-\omega_j^2 \M+ \K )\phiv_j = \0,
\end{equation}
where $\phiv_j$ is the eigenmode shape and $\omega_j$ the corresponding eigenfrequency. Assuming normalisation with respect to mass, the family of eigenmodes fulfils the following relationships:
\begin{equation}
\phiv_j^\text{T} \M \phiv_j =1,\qquad\phiv_j^\text{T} \K \phiv_j =\omega_j^2.
\end{equation}

Moreover, the modal displacement for a given mode $j$ is obtained by projection:
\begin{equation}\label{eq:modal_displ}
\ur_j = \phiv_j^\text{T} \M \U.
\end{equation}
Defining $\V = \Ut$ as the vector of nodal velocities, the modal velocity is also obtained by projection:
\begin{equation}\label{eq:modal_vel}
\vr_j = \phiv_j^\text{T} \M \V.
\end{equation}

In the remainder of the article, it is assumed that the damping formulation is such that the modes of the conservative system diagonalise the damping matrix $\C$ as well. The Rayleigh proportional damping law, commonly used in FE formulation, which imposes $\C$ to be a summation of mass and stiffness matrices with two independent parameters, is known to be compatible with this assumption. More generally, the reader is referred to~\cite{Caughey60,Caughey65,ADHIKARI2006} for discussions on the formulation of $\C$ such that the system possesses classical normal modes. As we are also mostly interested in lightly damped systems, it is also assumed that the damping of the master modes is small. By introducing the modal damping ratio $\xi_j$ of the $j$-th mode as
\begin{equation}
\phiv_j^\text{T} \C \phiv_j = 2\xi_j\omega_j,
\end{equation}
the assumption of light damping ($\xi_j \ll 1$) for the first modes will be generally made in the rest of the paper. The eigenproblem of the non-conservative linear system can then be written as
\begin{equation}
(\lambda_j^2 \M+ \lambda_j \C +\K )\phiv_j = \0,
\label{eq:lin_eig_secondorder}
\end{equation}
whose eigenvalues are the complex conjugate pairs
\begin{equation}
\lambda_j = -\xi_j\omega_j \pm \iu \omega_j \sqrt{1-\xi_j^2}.
\label{eq:eigenvalcc}
\end{equation}

A first-order, state-space formulation, is introduced for deriving the main part of the calculations. As a direct consequence, the size of the problem will double and become~$2N$. As shown for example in~\cite{Tisseur2001,JAIN2021How}, numerous different formulations can be used to write Eq.~\eqref{eq:FEstart}, leading to different properties in terms of the symmetry of the resulting matrices. In this contribution, $\M$ is assumed to be non-singular and the following first-order non-symmetric formulation is selected:
\begin{subequations}
\begin{align}
& \M \Vt + \C \V + \K \U +\G(\U,\U) + \H(\U,\U,\U) = \0,\label{eq:eom-a}
\\
& \M\Ut = \M\V.\label{eq:eom-b}
\end{align}\label{eq:eom}
\end{subequations}
In Eq.~\eqref{eq:eom-b}, the mass matrix has been added for symmetry reasons in the upcoming formula. Other choices, leading to symmetric formulations, could have been used. This choice is justified by the following arguments. First, as it will be shown in the next developments, the state-space formulation is used essentially for readability and to recover important symmetry properties. However, a special emphasis will be put throughout the calculations in order to solve $N$-dimensional problems rather than $2N$, by exploiting the relationship between displacement and velocity arising from the fact that the initial problem is second-order. Second, further extensions of the methods will finally lead to a non-symmetric formulation when including forces that will break this property. Thus, for the sake of generality, it has been found more convenient to directly work in such a setting, which will involve defining two projection basis with right and left eigenvectors.



In vibration theory, eigenvalues are complex conjugate and come by pairs following Eq.~\eqref{eq:eigenvalcc}, and two of them are needed to form a vibration mode. In state-space form, one can sort them either one next to the other, or put the first $N$ ({\em e.g.} with positive sign on the imaginary part) and complete the sorting by the last $N$ complex conjugates. This second choice is here retained such that the $j$-th vibration mode of the second-order system,  is now split into two  complex conjugate modes, corresponding to the $j$-th  and  $(j+N)$-th lines. Consequently the eigenspectrum is sorted according to the following order, $\forall \, j \, \in \, [1,N]$:
\begin{subequations}
\begin{align} 
\lambda_j &= -\xi_j \omega_j + \iu \omega_j \sqrt{1-\xi_j^2},\label{eq:def_complex_eigvals-a}
\\
\lambda_{j+N} &= \bar{\lambda}_j = -\xi_j \omega_j - \iu \omega_j \sqrt{1-\xi_j^2}.\label{eq:def_complex_eigvals-b}
\end{align}\label{eq:def_complex_eigvals}
\end{subequations}
This choice is appealing since the second half of the complex problem is simply given by the complex conjugate of the first. This has consequences in all the upcoming expressions as, in most of the derivations, the index $j$ can span only the first half, $j\in [1,N]$, the second half being implicitly verified using the conjugation operation without extra work. In order to come back to the real $j$-th vibration mode, one needs to pick the pair $(j,j+N)$ in the complex eigenproblem.  The corresponding right complex eigenvectors $\v{Y}_j$ of the first-order problem, following the same classification,  read, for $j \in [1,N]$:
\begin{subequations}\begin{align} 
\v{Y}_j = &\begin{bmatrix}
\phiv_j\lambda_j\\ \phiv_j
\end{bmatrix},
\\
\v{Y}_{j+N} = \bar{\v{Y}}_j = &\begin{bmatrix}
\phiv_j\bar{\lambda}_j\\ \phiv_j
\end{bmatrix}.
\end{align}\label{eq:def_complex_eig}\end{subequations}
Again, the second half for index ranging from $N+1$ to $2N$ is simply given by the complex conjugate.  Note that the right eigenvectors are expressed directly in terms of the real modes $\phiv_j$ of the second-order system. The right eigenvectors $\v{Y}_j$ are solution of the following eigenproblem:
\begin{equation}
\left(
\lambda_s
\begin{bmatrix}
\M & \0
\\
\0 & \M
\end{bmatrix}
+
\begin{bmatrix}
\C & \K
\\
-\M & \0
\end{bmatrix}
\right)
\v{Y}_s
=
\0.
\label{eq:lin_eig_right}
\end{equation}
This eigenvalue problem is valid for all $s$, nevertheless it is sufficient to span $s \in [1,N]$ since the second half is the complex conjugate, thanks to the relationships $\lambda_{s+N} = \bar{\lambda}_s$ and $\v{Y}_{s+N} = \bar{\v{Y}}_s$. Since the retained first-order system is not symmetric, one also needs to define the complex left eigenvectors $\v{X}_j$ as
\begin{subequations}\begin{align} 
\v{X}_j = \frac{1}{\lambda_j - \bar{\lambda}_j}
&\begin{bmatrix}
\phiv_j \\
- \phiv_j\bar{\lambda}_j
\end{bmatrix},
\\
\v{X}_{j+N} = \bar{\v{X}}_j = \frac{1}{\bar{\lambda}_j - \lambda_j}
&\begin{bmatrix}
\phiv_j \\
-\phiv_j\lambda_j 
\end{bmatrix},
\end{align}\label{eq:def_complex_eig_left}\end{subequations}
where $j \in [1,N]$ spans the real modes of the second-order system.

The left eigenmodes are solutions of the linear problem:
\begin{equation}
\v{X}_s^\text{T}
\left(
\lambda_s
\begin{bmatrix}
\M & \0
\\
\0 & \M
\end{bmatrix}
+
\begin{bmatrix}
\C & \K
\\
-\M & \0
\end{bmatrix}
\right)
=\0,
\label{eq:lin_eig_left}
\end{equation}
where, as before, it is sufficient to write Eq.\eqref{eq:lin_eig_left} for $s \in [1,N]$. 

As usual in vibration theory, the left and right eigenvectors share important orthogonality properties. For real vibration modes, both orthogonality with respect to mass and stiffness matrices are fulfilled. For the first-order non-symmetric system considered herein, the equivalent of the mass orthonormalisation reads
\begin{equation}
\v{X}_r^\text{T}
\begin{bmatrix}
\M & \0
\\
\0 & \M
\end{bmatrix}
\v{Y}_s = \delta_{sr},
\label{eq:massnorm_complex}
\end{equation}
with  $s,r \in [1,2N]$ and $\delta_{sr}$ the Kroneker delta. The equivalent of the orthogonality condition with respect to stiffness reads
\begin{equation}
\v{X}_r^\text{T}
\begin{bmatrix}
\C & \K
\\
-\M & \0
\end{bmatrix}
\v{Y}_s = -\lambda_r\delta_{sr},
\label{eq:stiffnorm_complex}
\end{equation}
with  $s,r \in [1,2N]$. Due to the first-order formulation and the complex conjugate eigenfrequencies \eqref{eq:def_complex_eigvals}, a {\em complexification} of the problem is used to conduct most of the calculations. Coming back to real coordinates using a {\em realification} will be addressed in Section~\ref{sec:compreal} to close the developments.

%
%
%
%
%
%
%


\subsection{Parametrisation method and invariance equations}

In this section, the definitions needed for the nonlinear mapping and the reduced dynamics are introduced. The method relies on the parametrisation method of invariant manifolds, first introduced by Cabr{\'e}, Fontich and de la Llave in~\cite{Cabre1,Cabre2,Cabre3}. The book by Haro {\em et al.}~\cite{Haro} details the method with the aim of developing effective computations for physical problems. It has already been applied in vibration theory for systems in modal space in~\cite{Haller2016,PONSIOEN2018}, \red{where the denomination SSM has been firstly introduced. The existence and uniqueness of these sought invariant manifolds under appropriate smoothness and non-resonance conditions have been demonstrated in~\cite{Cabre1,Haller2016}.} More recent progress focuses on working directly in the physical space, with in view application to structures modelled with the FEM. This has been realised using either a normal form approach~\cite{artDNF2020,AndreaROM,YichangVib}, or the parametrisation method~\cite{JAIN2021How,MingwuLi2021_1,MingwuLi2021_2}. Here we elaborate on the parametrisation method having numerous technical differences in the course of the development as compared to~\cite{JAIN2021How,MingwuLi2021_1,MingwuLi2021_2}.

\red{The general idea is to reduce the dynamics to the invariant manifold tangent to the eigenvectors selected as master modes, which is to say, to reduce the dynamics of the whole system to that on the approximated master SSM}. Since the invariant manifold is a curved subset in phase space, a nonlinear mapping is defined. Let us assume that $n$ master coordinates are selected, with $n\ll N$. These master coordinates are linked to their corresponding vibration modes and the searched invariant manifold is the nonlinear continuation of the subspace spanned by the $n$ second-order vibration modes. In phase space, the invariant manifold is $2n$-dimensional due to the fact that two coordinates (basically displacements and velocity) are needed. In order to describe the reduced dynamics on this manifold, we introduce $2n$ \textit{normal} coordinates $\v{z}$, following the denomination introduced in~\cite{touze03-NNM,TOUZE:JSV:2006}. The $2N$ original coordinates $\U$ and $\V$ are then expressed as a function of the new normal coordinates $\v{z}$ as
\begin{subequations}
\begin{align}
&\U = \WUfun,
\\
&
\V = \WVfun,
\end{align}\label{eq:mappings_compact}
\end{subequations}
where the two nonlinear mapping functions $\WUfunsim$ and $\WVfunsim$ are the unknowns to be computed. Note that in contrast to~\cite{JAIN2021How,MingwuLi2021_1,MingwuLi2021_2}, two mappings are introduced, a feature that will be key to reduce numerous computations, by making explicit the link between them, and recovering whenever possible a $N$-dimensional problem.

The reduced dynamics governs the evolution onto the corresponding invariant manifold. At this stage it is also unknown and is assumed to write
\begin{equation}
\v{\dot{z}} = \ffun.
\label{eq:reddyn_compact}
\end{equation}
The aim of the method is to compute $\WUfunsim$, $\WVfunsim$ and $\ffunsim$. The reduced-order dynamics is then given by $\ffunsim$, while the nonlinear mappings $\WUfunsim$ and $\WVfunsim$ allow one to pass from the physical space to the invariant manifold and give the functional relationship between the {\em normal} and original (physical) coordinates.

In order to solve for the unknowns, the key is to derive the so-called {\em invariance equation}~\cite{Cabre3,Haro} which states that the computed manifold is indeed invariant. The general formulation of the invariance equation given in~\cite{Cabre3,Haro} is here adapted to the mechanical context. Note that the invariance equation is also used in~\cite{Haller2016} for mechanical problems. In this contribution, the distinctive feature relies in the fact that both lines of the mappings, related to displacement and velocities, are followed during the calculations through $\WUfunsim$ and $\WVfunsim$. This allows one in particular to keep track of the mechanical characteristic features (mass matrix, linear and nonlinear stiffness) throughout the calculations, express more closely the relationships existing between $\WUfunsim$ and $\WVfunsim$, and make a clear connections to earlier works. Finally, it will also allow us to provide numerous expressions with $N$-dimensional matrices instead of $2N$. The procedure to derive the invariance equation consists in differentiating Eq.~\eqref{eq:mappings_compact} with respect to time, and then replace all time dependencies thanks to Eq.~\eqref{eq:reddyn_compact} to eventually arrive at a time-independent equation. Deriving Eq.~\eqref{eq:mappings_compact} with respect to time and using Eq.~\eqref{eq:reddyn_compact} leads to
\begin{subequations}\begin{align}
&\Ut 
= \nabla_\v{z}\WUfun \, \v{\dot{z}} 
= \nabla_\v{z}\WUfun \,  \ffun
= \sum_{s=1}^{2n} \dfrac{\partial \WUfun}{\partial z_s} f_s(\v{z}),
\\
&\Vt 
= \nabla_\v{z}\WV(\v{z}) \, \v{\dot{z}} 
= \nabla_\v{z}\WVfun \,  \ffun
= \sum_{s=1}^{2n} \dfrac{\partial \WVfun}{\partial z_s} f_s(\v{z}).
\end{align}\label{eq:dotted_compact}\end{subequations}

Substituting in the first-order equations of motion \eqref{eq:eom}, one arrives at the invariance equation which reads, for the mechanical problem with geometric nonlinearities
\begin{subequations}
\begin{align}
& \M \;\nabla_\v{z} \WVfun\; \ffun + \C \WVfun + \K \WUfun
 +\G(\WUfun,\WUfun) + \H(\WUfun,\WUfun,\WUfun) = \0,
\\
& \M\nabla_\v{z}\WUfun\; \ffun = \M\WVfun.
\end{align}\label{eq:invariance_compact}
\end{subequations}
These nonlinear equations can be solved locally by using asymptotic expansions in the unknown (the normal coordinate $\v{z}$), as proposed in~\cite{Cabre3,Haro}. The remainder of the paper details how this procedure can be written for an arbitrary order such that high-order converged solutions can be computed. One of the main difficulty resides in tracking all the terms having the same order, since they can originate from different sources. This process is handled step by step in the next sections.

\subsection{Asymptotic expansions and homological equations}

Let us assume that $n$ master modes have been selected for the analysis. In the first-order form, this corresponds to $2n$ complex conjugate modes $\v{Y}_s$ such that the index $s$ will span from 1 to $2n$. The choice of the master modes is left to the user and is guided by physical reasoning and the dynamics one wants to simulate with the ROM, see {\em e.g.}~\cite{ROMGEOMNL} for a discussion.

Both unknown nonlinear mappings and reduced dynamics can be expressed as polynomial expansions of the $2n$ normal coordinates. Let us denote as $o$ the maximum order reached by the expansion, which is arbitrary at the moment and will be seen as a convergence parameter for the solution. The unknown functions are thus expanded as
\begin{equation}\label{eq:expand1}
\WUfun =  \sum_{p=1}^o \P{\WUfun}, \quad \WVfun = \sum_{p=1}^o \P{\WVfun}, \quad \ffun =  \sum_{p=1}^o \P{\ffun},
\end{equation}
where the shortcut notation $\P{.}$ is used to indicate a polynomial term of order $p$. Detailed expressions of the polynomial expansions will be given when needed in the remainder of the paper, but here this simple notation is  used to underline the main points of the method. The constant terms for $p=0$ are not taken into account in all the expansions since it is assumed that the fixed point of the system \eqref{eq:eom}, which represents the structure at rest, is at the origin of the phase space. 

The order-$p$ holomological equations correspond to selecting all the terms of order $p$ from the invariance equation \eqref{eq:invariance_compact}. Using the notation $\P{.}$ they can be simply written as
\begin{subequations}
\begin{align}
& \M\P{\nabla_\v{z} \WVfun\; \ffun} + \C\P{\WVfun} + \K\P{\WUfun} + \P{\G(\WUfun,\WUfun)} + \P{\H(\WUfun,\WUfun,\WUfun)} = \0,
\\
& \M\P{\nabla_\v{z}\WUfun\; \ffun} = \M\P{\WVfun}.
\end{align}\label{eq:homological_compact}
\end{subequations}


\subsection{First-order solution: tangency to linear eigenspaces}

The process of solving the order-$p$ homological equations is sequential in nature since orders lower than $p$ will create new order-$p$ terms, due to the presence of the nonlinearity. In this section, the first-order is solved to initiate the process, showing that the linear solution is retrieved. In terms of geometry in phase space, this means that the searched manifold is tangent at origin to the linear space spanned by the master modes.

Rewriting Eq.~\eqref{eq:homological_compact} for $p=1$ yields
\begin{subequations}\begin{align}
& \M\P[1]{\nabla_\v{z} \WVfun\; \ffun} + \C\P[1]{\WVfun} + \K\P[1]{\WUfun} = \0,
\\
& \M\P[1]{\nabla_\v{z}\WUfun\; \ffun} = \M\P[1]{\WVfun}.
\end{align}\label{eq:homological1_compact}
\end{subequations}
Since only linear terms are here retained by application of the operation $\P[1]{\cdot}$, the nonlinear quadratic and cubic terms $\G$ and $\H$ are simply discarded at this order. The linear terms of the three unknowns can be rewritten as matrix-vector products as
\begin{subequations}\begin{align}
& \P[1]{\WVfun} = \p[1]{\WV} \v{z}, \\
& \P[1]{\WUfun} = \p[1]{\WU} \v{z},\\
& \P[1]{\ffun} = \p[1]{\v{f}} \v{z}, 
\end{align}\end{subequations}
where  $\p[1]\WV$ and $\p[1]\WU$ are matrices of size $N\times 2n$ and $\p[1]{\v{f}}$ is a square matrix $2n\times 2n$. Using the fact that the gradient of linear functions are simple, one can then rewrite \eqref{eq:homological1_compact} as
%
\begin{subequations}\begin{align}
& \M\p[1]{\WV}\; \p[1]{\v{f}}\v{z} + \C\p[1]{\WV}\v{z} + \K\p[1]{\WU}\v{z}= \0,
\\
& \M\p[1]{\WU}\; \p[1]{\v{f}}\v{z} = \M\p[1]{\WV}\v{z}.
\end{align}\end{subequations}

Collecting into $2N\times 2N$ matrices, and using the fact that the previous equations 
must be fulfilled for any $\v{z}$, leads to
\begin{equation}
\begin{bmatrix}
\M & \0
\\
\0 & \M
\end{bmatrix}
\begin{bmatrix}
\p[1]{\WV}
\\
\p[1]{\WU}
\end{bmatrix}
\; \p[1]{\v{f}}
+
\begin{bmatrix}
\C & \K
\\
-\M & \0
\end{bmatrix}
\begin{bmatrix}
\p[1]{\WV}
\\
\p[1]{\WU}
\end{bmatrix}
=\0.
\label{eq:homological1_eig}
\end{equation}
One recognises the linear eigenproblem as stated in Eq.~\eqref{eq:lin_eig_right}. In order to write the solutions in compact form, one can introduce the $2n\times2n$ matrix of master mode eigenvalues $\v{\Lambda} = \text{diag}[\lambda_1,\;\ldots,\;\lambda_{2n}]$. Using the fact that eigenvalues are complex conjugate such that $\lambda_{i+n} = \bar{\lambda}_i$, it yields
\begin{equation}\label{eq:master_eigv}
\v{\Lambda} = 
\text{diag}[\lambda_1,\;\ldots,\;\lambda_n,\;\bar{\lambda}_1,\;\ldots \;\bar{\lambda}_n].
\end{equation}
In the same lines, one can introduce the $2N\times2n$ matrix $\v{\Phi}$ of real master eigenfunctions $\phiv_i$, for $i\in [1,n]$, as
\begin{equation}\label{eq:master_eigs}
\v{\Phi} = 
\begin{bmatrix}
\phiv_1 &\phiv_2 &\ldots &\phiv_n  &\phiv_1 &\phiv_2 &\ldots& \phiv_n 
\end{bmatrix}.
\end{equation}
In general when dealing with second-order real problems, this matrix, which allows one to go directly from physical to modal space, has only $n$ columns. Here due to the use of first-order formulation, this matrix needs to have $2n$ columns by repeating the master eigenvectors.


Then the solution to Eq.~\eqref{eq:homological1_eig}  is given by the linear eigenvectors and eigenvalues:
\begin{subequations}\label{eq:solorder1lin}
\begin{align}
& \begin{bmatrix}
\p[1]{\WV}
\\
\p[1]{\WU}
\end{bmatrix}
=
\begin{bmatrix}
\v{\Phi}\v{\Lambda}\\\v{\Phi}
\end{bmatrix}=\begin{bmatrix}
\v{Y}_1 &\ldots& \v{Y}_n &\bar{\v{Y}}_1 &\ldots& \bar{\v{Y}}_n
\end{bmatrix},\label{eq:solorder1lin-map}
\\
& \p[1]{\v{f}} = \v{\Lambda}.\label{eq:solorder1lin-rdyn}
\end{align}
\end{subequations}


This result shows that the linear part of the mapping is simply given by the eigenfunction of the selected master modes.   The higher-order terms will bring corrections to the mappings, by taking into account the non-resonant nonlinear couplings between the modes. The linear part of the reduced dynamics is left unchanged since the eigenvalues are retrieved, meaning that at the linear level, the dynamics is governed by the modal uncoupled linear oscillator equations. Nonlinear terms will then introduce the needed corrections.

\section{Arbitrary order expansion}

In this Section, the detailed expressions of the order-$p$ homological equations are derived for an arbitrary order $p$. To that purpose, the asymptotic expansions need to be emphasised. 

\subsection{Nonlinear mappings and reduced dynamics}

Now that the first-order solutions are known, the asymptotic expansions of the unknown nonlinear mappings $\WUfunsim$ and $\WVfunsim$ can be rewritten up to the maximum order of the expansion $o$ as
\begin{subequations}\begin{align}
&\WUfun = \;\;\v{\Phi} \v{z}
+ \sum_{p=2}^o \P{\WUfun},
\\
&\WVfun = \v{\Phi}\v{\Lambda} \v{z}
+ \sum_{p=2}^o \P{\WVfun}.
\end{align}\label{eq:mappings_byorder}\end{subequations}
The generic order $p$ term for each of the two mappings  is a polynomial of order $p$ in the normal coordinate $\v{z}$. By making appear the different monomials of order $p$, 
one can formally write, :
\begin{subequations}\label{eq:polyexpandmon}
\begin{align}
&\P{\WUfun} = 
\sum_{i_1=1}^{2n} \sum_{i_2=1}^{2n} \ldots \sum_{i_p=1}^{2n} 
\p{\WU_{i_1 i_2 \ldots i_p}} \; z_{i_1} z_{i_2} \ldots z_{i_p},
\\
&\P{\WVfun} = 
\sum_{i_1=1}^{2n} \sum_{i_2=1}^{2n} \ldots \sum_{i_p=1}^{2n} 
\p{\WV_{i_1 i_2 \ldots i_p}} \; z_{i_1} z_{i_2} \ldots z_{i_p}.
\end{align}
\end{subequations}
In these expressions, $z_{i_1} z_{i_2} \ldots z_{i_p}$ represents a generic order-$p$ monomial having as coefficient a vector  $\p{\WU_{i_1 i_2 \ldots i_p}}$ (and similarly for $\p{\WV}$). Each index $i_k$ spans all the master modes from $1$ to $2n$ so that the summations span all the possible combinations of order-$p$ monomials.

To introduce a more compact notation, let us  define $\Is$ as the generic set of indices of order $p$:
\begin{equation}
\Is= \lbrace i_1 i_2 \ldots i_p \rbrace,
\end{equation}
which gathers all indices involved in a given monomial. 
The monomial associated to $\Is$, i.e.\ the order-$p$ product of normal coordinates, 
will be denoted as  $\p{\pi_\Is}$, with
\begin{equation}\label{eq:defpiI}
\p{\pi_\Is} = z_{i_1} z_{i_2} \ldots z_{i_p}.
\end{equation}
It is important to notice that the way the set $\Is$ is constructed, does not involve grouping of repeated indices nor specification of their multiplicity. For instance, let us take the order-$p$ monomial $z_2 z_5^2 z_6^{p-3}$; for this monomial the set would be $\Is=\lbrace 2 5 5 6 \ldots \rbrace$. This shows that the cardinal number of $\Is$ is always $p$, so index with multiplicity higher than one are simply repeated multiple times inside $\Is$.

Substituting in \eqref{eq:polyexpandmon} allows writing the generic order-$p$ of the nonlinear mappings in the form
\begin{subequations}\begin{align}
&\P{\WUfun} = 
\sum_\Is \p{\WU_\Is} \; \p{\pi_\Is},
\\
&\P{\WVfun} = 
\sum_\Is \p{\WV_\Is} \; \p{\pi_\Is},
\end{align}\label{eq:mappings_I}\end{subequations}
where the summation spans all possible $\Is$ of order $p$.


The same expansions are needed for the reduced dynamics which intervenes in the order-$p$ homological equations. Following the same notations, one can expand  $\ffun$ as a polynomial function of the normal coordinates, the linear term being known thanks to Eq.~\eqref{eq:solorder1lin-rdyn}, as
\begin{equation}
\ffun = \v{\Lambda}\v{z} + \sum_{p=2}^o \P{\ffun}.
\end{equation}
Similarly, the generic order-$p$ writes
\begin{equation}
\P{\ffun} = 
\sum_{i_1=1}^{2n} \sum_{i_2=1}^{2n} \ldots \sum_{i_p=1}^{2n} 
\p{\v{f}_{i_1 i_2 \ldots i_p}} \; z_{i_1} z_{i_2} \ldots z_{i_p}
= \sum_\Is \p{\v{f}_\Is} \; \p{\pi_\Is}.
\label{eq:red_dyn_I}
\end{equation}
Using Eq.~\eqref{eq:reddyn_compact}, one can rewrite explicitly the reduced dynamics for each $s$ normal coordinate, where $s\in[1,2n]$ spans the master modes, as the following order-$o$ approximation
\begin{equation}
\dot{z}_s = \lambda_s z_s + \sum_{p=2}^o  \sum_\Is \p{f_{s\,\Is}}
\; \p{\pi_\Is} +\mathcal{O}(\p[o+1]{|\v{z}|})	
\label{eq:reddyn_orderp}
\end{equation}
At this stage, all the unknowns have been expressed with asymptotic expansions. The solutions at arbitrary order are given by replacing all the developments into the order-$p$ homological equation.

\subsection{Order-\textit{p} homological equations}

This sections aims at providing explicit expressions for 
the order-$p$ homological equations~\eqref{eq:homological_compact} by using the expansions derived in the previous section and considering that,
starting from the second-order,
the nonlinear polynomial terms $\G$ and $\H$  will generate contributions.
In order to collect all terms of order $p$ in \eqref{eq:homological_compact}, it is important to make the distinction between the terms that are directly of order $p$, from those that are created by products of terms with a lower order. 
It is also important to understand the {\em sequential} nature of the procedure. When arriving at order $p$, all the lower order mappings and reduced dynamics functions are already known and are denoted by $\P[<p]{\WUfun}$, $\P[<p]{\WVfun}$, $\P[<p]{\ffun}$, where we used the shortcut notation $\P[<p]{.}$ to describe all terms of order strictly lower than $p$.
Consequently, the unknowns are $\P{\WUfun}$, $\P{\WVfun}$, $\P{\ffun}$.


Let us examine how each term in Eqs.~\eqref{eq:homological_compact} can be made explicit. The two terms coming from the nonlinear polynomial restoring force $\G$ and $\H$ obviously depends on previously calculated mappings at order $<p$, so that one can write
\begin{subequations}\begin{align} 
&\P{\G(\WUfun,\WUfun)}
=
\P{\G(\P[<p]{\WUfun},\P[<p]{\WUfun})},
\\
&\P{\H(\WUfun,\WUfun,\WUfun)}
=
\P{\H(\P[<p]{\WUfun},\P[<p]{\WUfun},\P[<p]{\WUfun})}.
\end{align}\end{subequations}


Using the same notation as in the previous section and introducing $\p{\pi_\Is}$ given by Eq.~\eqref{eq:defpiI} as the generic order-$p$ monomial, one can now simply expand the nonlinear terms as polynomials of order $p$ with given tensor of coefficients $\p{\FG_\Is}$ and $\p{\FH_\Is}$ as
\begin{subequations}\label{eq:expandGHNEW}
\begin{align} 
&\P{\G(\WUfun,\WUfun)} = \sum_\Is \p{\FG_\Is}\p{\pi_\Is},
\\
&\P{\H(\WUfun,\WUfun,\WUfun)} = \sum_\Is \p{\FH_\Is}\p{\pi_\Is}.
\end{align}
\end{subequations}
With algebraic manipulations of polynomial representations, and using the set of indices $\Is = \lbrace i_1\, i_2 \,\ldots i_p \rbrace$ already introduced, explicit expressions of $\p{\FG_\Is}$ and $\p{\FH_\Is}$ can be derived. One needs just to notice that since $\G$ groups the quadratic terms, a term of order $p$ is necessarily the product of two terms of orders $k$ and $p-k$, with $k$ ranging from 1 to $p-1$. The same can be written for $\H$, being a cubic term and involving products of three lower order terms, such that
\begin{subequations}
\begin{align}
&
\p{\FG_\Is} = 
\sum_{k=1}^{p-1}
\G(\WU^{(k)}_{i_1 \ldots i_k},
\WU^{(p-k)}_{i_{k+1} \ldots i_p}),
\\
&
\p{\FH_\Is} = 
\sum_{k=1}^{p-2}\sum_{l=1}^{p-k-1}
\H(\WU^{(k)}_{i_1 \ldots i_k},
\WU^{(l)}_{i_{k+1} \ldots i_{k+l}},
\WU^{(p-k-l)}_{i_{k+l+1} \ldots i_p}) .
\end{align}
\end{subequations}


The most cumbersome terms to handle from Eq.~\eqref{eq:homological_compact} are those  composed by the gradient of the mapping functions contracted with the reduced dynamics:
\begin{subequations}\label{eq:nablaZfz}
\begin{align} 
&\P{\nabla_\v{z} \WUfun\; \ffun} = \P{\sum_{s=1}^{2n}\dfrac{\partial \WUfun}{\partial z_s} f_s(\v{z})},
\\
&\P{\nabla_\v{z} \WVfun\; \ffun} = \P{\sum_{s=1}^{2n}\dfrac{\partial \WVfun}{\partial z_s} f_s(\v{z})}.
\end{align}
\end{subequations}

In order to keep track of the different contributions and collect terms of the same order by separating known from unknown quantities, a simple formulation consists in dividing  both the mapping functions and the reduced dynamics  in  Eqs.~\eqref{eq:nablaZfz} into three terms: a linear term, an order $p$ term, and the intermediate ones, of order lower than $p$ but larger than $1$ that we will denote using the shortcut notation  $\P[\dblnk]{\cdot}$. Using these notations, one arrives at
\begin{subequations}\begin{align} 
&\P{\nabla_\v{z} \WUfun\; \ffun} = 
\P{\sum_{s=1}^{2n}
\left(\phiv_s + \dfrac{\partial \P[\dblnk]{\WUfun}}{\partial z_s}+\dfrac{\partial \P{\WUfun}}{\partial z_s} \right)
\left(\lambda_s z_s + \P[\dblnk]{f_s(\v{z})}+ \P{f_s(\v{z})}\right)},
\\
&\P{\nabla_\v{z} \WVfun\; \ffun} = 
\P{\sum_{s=1}^{2n}
\left(\phiv_s\lambda_s + \dfrac{\partial \P[\dblnk]{\WVfun}}{\partial z_s}+\dfrac{\partial \P{\WVfun}}{\partial z_s} \right)
\left(\lambda_s z_s + \P[\dblnk]{f_s(\v{z})}+ \P{f_s(\v{z})}\right)}. 
\end{align}\end{subequations}
This separation is meaningful since the operator $\P{\cdot}$  solely selects the terms of order $p$; consequently the terms from the linear mapping, $\phiv_s$, create an order $p$  only when multiplied with 
the order $p$ reduced dynamics $\P{f_s(\v{z})}$. The same applies for the linear term of the reduced dynamics, $\lambda_s z_s$, with the last term of the first parenthesis of order $p-1$.  One can then write
\begin{subequations}\begin{align}
&\P{\nabla_\v{z} \WUfun\; \ffun} = 
\sum_{s=1}^{2n}\left(  
\dfrac{\partial \P{\WUfun}}{\partial z_s} \lambda_s z_s
+ \phiv_s \P{f_s(\v{z})}
+ \P{\dfrac{\partial \P[\dblnk]{\WUfun}}{\partial z_s} \P[\dblnk]{f_s(\v{z})}}
\right),
\\
&\P{\nabla_\v{z} \WVfun\; \ffun} = 
\sum_{s=1}^{2n} \left(  
\dfrac{\partial \P{\WVfun}}{\partial z_s} \lambda_s z_s
+ \phiv_s \lambda_s \P{f_s(\v{z})}
+ \P{\dfrac{\partial \P[\dblnk]{\WVfun}}{\partial z_s}\P[\dblnk]{f_s(\v{z})}}
\right).
\end{align}\label{eq:diff_terms}\end{subequations}
In these two expressions, the only unknowns are the order $p$ mapping terms $\P{\WUfun}$ and $\P{\WVfun}$, as well as the reduced dynamics $\P{f_s(\v{z})}$. All other quantities are known, coming from lower order mappings and reduced dynamics. 
They have have been collected into the last term of the summation. Let us now focus on this last known term and expand its expression. Since it is a product, the order $p$ term is the product of two terms such that the sum of their orders equals $p$. Hence one can write
\begin{subequations}
\begin{align}
&
\sum_{s=1}^{2n}
\P{\dfrac{\partial \P[\dblnk]{\WUfun}}{\partial z_s} \P[\dblnk]{f_s(\v{z})}}
=\sum_{s=1}^{2n}\sum_{k=2}^{p-1}
\dfrac{\partial \P[p-k+1]{\WUfun}}{\partial z_s} \P[k]{f_s(\v{z})},
\\
&
\sum_{s=1}^{2n}
\P{\dfrac{\partial \P[\dblnk]{\WVfun}}{\partial z_s}\P[\dblnk]{f_s(\v{z})}}
=\sum_{s=1}^{2n}\sum_{k=2}^{p-1}
\dfrac{\partial \P[p-k+1]{\WVfun}}{\partial z_s}\P[k]{f_s(\v{z})}.
\end{align}\end{subequations}
These terms are also polynomials of order $p$, so they can be written in compact form by introducing the tensors of coefficients $\p{\FU_\Is}$ and $\p{\FV_\Is}$:
\begin{subequations}\begin{align} 
& \sum_{s=1}^{2n}\P{\dfrac{\partial \P[\dblnk]{\WUfun}}{\partial z_s} \P[\dblnk]{f_s(\v{z})}}
= \sum_\Is \p{\FU_\Is}\p{\pi_\Is},
\\
& \sum_{s=1}^{2n}\P{\dfrac{\partial \P[\dblnk]{\WVfun}}{\partial z_s} \P[\dblnk]{f_s(\v{z})}}
= \sum_\Is \p{\FV_\Is}\p{\pi_\Is}.
\end{align}\end{subequations}
Comparing the last two expressions and recalling the definition of the generic set $\Is = \lbrace i_1\, i_2 \,\ldots i_p \rbrace$, one can arrive at the following expressions for the newly introduced coefficients $\p{\FU_\Is}$ and $\p{\FV_\Is}$:
\begin{subequations}
\begin{align}
&
\p{\FU_\Is} =  
\sum_{s=1}^{2n}\sum_{k=2}^{p-1}\sum_{l=0}^{p-k}
\WU^{(p-k+1)}_{i_1 \ldots i_l s\, i_{l+k+1} \ldots i_p}
f^{(k)}_{s\, i_{l+1} \ldots i_{l+k}},
\\
&
\p{\FV_\Is} =  
\sum_{s=1}^{2n}\sum_{k=2}^{p-1}\sum_{l=0}^{p-k}
\WV^{(p-k+1)}_{i_1 \ldots i_l s\, i_{l+k+1} \ldots i_p}
f^{(k)}_{s\, i_{l+1} \ldots i_{l+k}}.
\end{align}
\end{subequations}
The unknown nonlinear mappings have been expanded into their polynomial form through compact expressions given in Eqs.~\eqref{eq:mappings_I}, which can be used to rewrite the first term in the right-hand side (RHS) of Eqs.~\eqref{eq:diff_terms} as
\begin{subequations}\label{eq:dzphiexpand}
\begin{align} 
& \sum_{s=1}^{2n}\dfrac{\partial \P{\WUfun}}{\partial z_s} \lambda_s z_s = 
\sum_{s=1}^{2n}\dfrac{\partial }{\partial z_s}\left(\sum_\Is \p{\WU_\Is} \p{\pi_\Is}\right)\lambda_s z_s,
\\
& \sum_{s=1}^{2n}\dfrac{\partial \P{\WVfun}}{\partial z_s} \lambda_s z_s = 
\sum_{s=1}^{2n}\dfrac{\partial}{\partial z_s}\left(\sum_\Is \p{\WV_\Is} \p{\pi_\Is}\right)\lambda_s z_s.
\end{align}
\end{subequations}
These terms can be simplified by noticing that the derivative with respect to $z_s$ is different from zero only if $z_s$ is contained inside $\p{\pi_\Is} = z_{i_1} z_{i_2} \ldots z_{i_p}$. Hence \eqref{eq:dzphiexpand} can be rewritten explicitly as
\begin{subequations}
\begin{align} 
& \sum_{s=1}^{2n}\dfrac{\partial}{\partial z_s}\left(\sum_\Is \p{\WU_\Is} \p{\pi_\Is}\right)\lambda_s z_s
=
\sum_\Is \p{\WU_\Is} (\lambda_{i_1} + \lambda_{i_2} + \ldots + \lambda_{i_p}) \p{\pi_\Is},
\\
& \sum_{s=1}^{2n}\dfrac{\partial}{\partial z_s}\left(\sum_\Is \p{\WV_\Is} \p{\pi_\Is}\right)\lambda_s z_s
=
\sum_\Is \p{\WV_\Is} (\lambda_{i_1} + \lambda_{i_2} + \ldots + \lambda_{i_p}) \p{\pi_\Is}.
\end{align}
\end{subequations}
Indeed, if we consider for instance $s=i_k$ in the summation on the left, for each set $\Is$ a non-vanishing contribution is obtained only if $i_k\in\Is$, because of the derivative of $\p{\pi_\Is}$ with respect to $z_{i_k}$ which is then  equal to $\p{\pi_\Is}/z_{i_k}$. Since this term is multiplied by $\lambda_{i_k} z_{i_k}$, the whole process makes reappear $\p{\pi_\Is}$, multiplied by the sum of all the $\lambda_{i_k}$ with $i_k\in\Is$ in the RHS term. This process is not influenced by a possible repetition of index in~$\Is$. Indeed, if one has for instance $i_k = i_{k+1}$, the derivative of $\p{\pi_\Is}$ with respect to $z_{i_k}$ can be still seen as the summation of multiple derivatives of the same variable\footnote{This is a simple application of the fact that  $\partial (x^2)/\partial x = 2 x$ can be seen as $\partial (x_1 x_2)/ \partial x_1 + \partial (x_1 x_2)/ \partial x_2 = x_2+x_1$ also when $x_1=x_2$.}.


The appearance of the summation of eigenvalues is of crucial importance for the rest of the development. Let us denote this term by $\sigma_\Is$:
\begin{equation}\label{eq:defsigmaIsum}
\sigma_\Is = \lambda_{i_1} + \lambda_{i_2} + \ldots + \lambda_{i_p},
\end{equation}
with $\Is = \lbrace i_1 \,i_2 \,\ldots i_p \rbrace$. This term is responsible for the nonlinear resonance and the different solutions (or {\em styles}) one can select to solve for the homological equations. This will be further explained in the next sections.
Setting the different contributions together, the terms composed by the gradient of the mapping functions contracted with the reduced dynamics, introduced in Eqs.~\eqref{eq:nablaZfz}, can finally be rewritten as
\begin{subequations}
\begin{align} 
&\P{\nabla_\v{z} \WUfun\; \ffun} =  
\sum_\Is \left(
\p{\WU_\Is} \sigma_\Is
+ \sum_{s=1}^{2n}(\phiv_s \,\p{f_{s\,\Is}})
+ \p{\FU_\Is}
\right)\p{\pi_\Is},
\\
&\P{\nabla_\v{z} \WVfun\; \ffun} = 
\sum_\Is \left(
\p{\WV_\Is} \sigma_\Is
+ \sum_{s=1}^{2n}(\phiv_s \lambda_s \,\p{f_{s\,\Is}})
+ \p{\FV_\Is}
\right)\p{\pi_\Is},
\end{align}\label{eq:diff_terms_I}
\end{subequations}
where the terms $\P{f_s(\v{z})}$ have also been expressed in terms of the monomials $\p{\pi_\Is}$ using Eq.~\eqref{eq:red_dyn_I}.

We are now in position of giving a detailed expression of the order-$p$ homological equations \eqref{eq:homological_compact}. Thanks to the previous developments, all terms have been rewritten, and using Eqs.~\eqref{eq:expandGHNEW} and \eqref{eq:diff_terms_I} allows one to expand \eqref{eq:homological_compact} as
\begin{subequations}
\begin{align}
& \sum_\Is\left(
\M \p{\WV_\Is}\sigma_\Is 
+ \sum_{s=1}^{2n}\left(\M\phiv_s \lambda_s \,\p{f_{s\,\Is}} \right)
+ \M \p{\FV_\Is}
+ \C \p{\WV_\Is}
+ \K \p{\WU_\Is} + \p{\FG_\Is} + \p{\FH_\Is}
\right)\p{\pi_\Is} = \0,
\\
& 
\sum_\Is\left(
\M \p{\WU_\Is}\sigma_\Is 
+ \sum_{s=1}^{2n}\left(\M\phiv_s  \p{f_{s\,\Is}}\right)
+ \M\p{\FU_\Is}-\M\p{\WV_\Is}\right)\p{\pi_\Is}
=
\0.
\label{eq:homo_eqs}
\end{align}
\end{subequations}
These equations have to be verified for any monomial term $\p{\pi_\Is}$, i.e. for any value of $\v{z}$, meaning that each term of the summation over all possible sets $\Is$ needs to be equal to zero. This allows one to rewrite a general equation for the unknowns of the problem at this stage. By keeping the unknowns on the left-hand side and moving the known terms on the right-hand side, one has:
\begin{equation}
\left(\sigma_\Is 
\begin{bmatrix}
\M & \0
\\
\0 & \M
\end{bmatrix}
+
\begin{bmatrix}
\C & \K
\\
-\M & \0
\end{bmatrix}
\right)
\begin{bmatrix}
\p{\WV_\Is}\\ \p{\WU_\Is}
\end{bmatrix}
+ \sum_{s=1}^{2n}
\p{f_{s\,\Is}} 
\begin{bmatrix}
\M & \0
\\
\0 & \M
\end{bmatrix}
\v{Y}_s
=
\begin{bmatrix}
- \M \p{\FV_\Is} - \p{\FG_\Is} - \p{\FH_\Is}\\
- \M \p{\FU_\Is}
\end{bmatrix} \, ,
\label{eq:homo}
\end{equation}
where the eigenvectors $\v{Y}_s$ are made appear  for compactness. In short, Eq.~\eqref{eq:homo} details the order-$p$ homological equation for a given monomial $\p{\pi_\Is}$ defined by a given set of indices $\Is$. This vectorial equation is $2N$-dimensional 
and must be solved for any set of indices $\Is= \lbrace i_1 \,i_2 \,\ldots i_p \rbrace$, with $i_k \in [1,2n]$.
 
The system \eqref{eq:homo} is underdetermined because both the mappings $\p{\WV_\Is},\p{\WU_\Is}$ and the reduced dynamics $\p{f_{s\,\Is}}$ are unknown; therefore there are different solution strategies which give rise to different {\em styles} of parametrisation~\cite{Haro}. This will be emphasised in the next development. Before proceeding, it is important to remark that one can separate the $2n$ lines corresponding to the master coordinates from the $2(N-m)$ lines corresponding to the slave ones. This gives rise to the so-called tangent and normal part of the order-$p$ homological equations~\cite{Haro}, which are obtained by projecting respectively onto the subspace spanned by the master and the slave modes, respectively. Separating these two contributions is crucial in order to provide general solutions for the unknown mappings and reduced dynamics. Indeed, the projection on the slave modes (normal part) does not generate terms from the reduced dynamics, because the $\v{Y}_s$ contains  only master modes such that their projection on the slave modes vanishes. Consequently, the normal part of the homological equation is not underdetermined, as opposed to the tangent part. For this reason, only the tangent part is derived hereafter in order to focus the discussion on the choice of the styles.


The idea being of projecting Eq.~\eqref{eq:homo}  onto the set of master modes, one can define $\p{{\theta}_{r\Is}}$ as the projection of the two mappings $\p{\WU_\Is}$ and $\p{\WV_\Is}$ on the left eigenvector  $\v{X}_r$  as 
\begin{equation}\label{eq:defthetaproj}
\p{{\theta}_{r\Is}}
\doteq
\v{X}_r^\text{T}
\begin{bmatrix}
\M & \0
\\
\0 & \M
\end{bmatrix}
\begin{bmatrix}
\p{\WV_\Is}\\ \p{\WU_\Is}
\end{bmatrix}.
\end{equation}
Here, the index $r$ spans only the master modes, $r \in [1,2n]$, such that the size of the tangent part is very small as compared to the complete problem.
In short, $\p{{\theta}_{r\Is}}$  coincides with the $r$-th component mappings in the modal basis. Indeed, recalling that we are using complex modes, it is logical that this component involves both parts of the mapping $\p{\WU_\Is}$ and $\p{\WV_\Is}$.

In light of the orthogonality properties \eqref{eq:massnorm_complex} and \eqref{eq:stiffnorm_complex}, one can easily demonstrate that:
\begin{equation}
\v{X}_r^\text{T}
\begin{bmatrix}
\C & \K
\\
-\M & \0
\end{bmatrix}
\begin{bmatrix}
\p{\WV_\Is}\\ \p{\WU_\Is}
\end{bmatrix}
=
-\lambda_r
\p{{\theta}_{r\Is}}.
\end{equation}
Let us finally denote as $\p{{g}_{r\Is}}$ the projection of the right-hand side of Eq.~\eqref{eq:homo} on the left eigenvector~$\v{X}_r$:
\begin{equation}
\p{{g}_{r\Is}} \doteq \v{X}_r^\text{T} 
\begin{bmatrix}
- \M \p{\FV_\Is} - \p{\FG_\Is} - \p{\FH_\Is}\\
- \M \p{\FU_\Is}
\end{bmatrix}.
\label{eq:gtilde}
\end{equation}
With these quantities, one can arrive at a compact expression of the tangent  homological equation, obtained by  projecting Eq.~\eqref{eq:homo} on each of the master mode $\v{X}_r$. A generic row of the tangent homological equation then reads, $\forall r \in [1,2n]$,
\begin{equation}
(\sigma_\Is - \lambda_r ) \p{{\theta}_{r\Is}} + \p{f_{r\Is}} = \p{{g}_{r\Is}}.
\label{eq:homo_tang}
\end{equation}

In \eqref{eq:homo_tang}, $\p{{g}_{r\Is}}$ is known, the two unknowns being $\p{{\theta}_{r\Is}}$, which is related to the coefficients of the nonlinear mappings, and $\p{f_{r\Is}}$, which is the term describing the reduced-order dynamics. The underdeterminacy is now obvious and one can see that different solutions for either $\p{{\theta}_{r\Is}}$ or $\p{f_{r\Is}}$ are possible. These solutions will be discussed in Section~\ref{sec:styles}. However, since $\p{{\theta}_{r\Is}}$ is multiplied by $\sigma_\Is - \lambda_r$, a first discussion is needed to treat the case when this factor vanishes or tends to very small values. This defines the well-known {\em resonance conditions} that are a cornerstone of dynamical systems~\cite{gucken83}, and appear explicitly in the normal form theory~\cite{IoossAdel,Murdock}. The next section is devoted to the analysis of  nonlinear resonances and, most importantly for vibratory systems, the appearance of {\em trivial} resonances will be underlined and separated from {\em internal} resonances, following the terminology used for example in~\cite{touze03-NNM,TOUZE:JSV:2006,TouzeCISM}.

\subsection{Resonances}\label{sec:res}

The emergence of nonlinear resonances as in Eq.~\eqref{eq:homo_tang} is a well-known fact in dynamical system theory. It has been underlined since the pioneering works by Poincar{\'e} and Dulac on normal form theory~\cite{Poincare,Dulac1912} and is discussed in all classical mathematical textbooks~\cite{Kuznetsov,Wiggins}. In this contribution, since we are interested in vibratory systems, real and imaginary parts do not play the same role. In nonlinear vibration theory, one is generally interested in lightly damped systems since the presence of large damping mostly inhibits the appearance of peculiar nonlinear phenomena and enforces the predominance of linear behaviour with strong temporal decays. Even though the damping ratios are usually increasing with the frequency, it is thus very common to assume lightly damped master modes such that $\forall\,s \in [1,n], \; \xi_s \ll 1$. In this case the predominant part of the spectrum is driven by the pairs of complex conjugate terms 
such that one can assume $\lambda_s \simeq \pm \iu \omega_s$. It is also important to understand that, even though the resonance are not exactly fulfilled with the presence of the real parts due to damping, the closeness of the resonance is enough to make appear the problem of small divisors. The writing of resonance relationships must thus be written for exact and close fulfilment 
of the condition in order to have a uniform treatment. 

This simplification of the eigenspectrum comes with two consequences. First, the appearance of {\em trivial} resonances at each odd order $p$. Second, the definition of {\em internal resonances} as commensurability relationships between the eigenfrequencies only, a common feature in nonlinear vibration theory~\cite{Nayfeh79,Nayfeh00}. Let us first discuss the trivial resonances\footnote{The wording {\em trivial resonance} has been introduced in~\cite{touze03-NNM,TOUZE:JSV:2006}, following~\cite{MannevilleENG}.}.

Trivial resonances are present for vibratory systems at each odd order. The first ones appear at the third-order due to the fulfilment of the trivial relationship $+\iu \omega_p =  +\iu \omega_j - \iu \omega_j + \iu \omega_p$, $\forall (j,p)$. Since we are working here at arbitrary order, one needs to take care of all trivial resonances appearing at all odd orders. This remark also explains why odd and even orders play two very different roles in nonlinear vibration and reduction methods. In order to select all the trivial resonances for any order $p$, one has to group together pairs of complex conjugates eigenvalues $(\lambda_s,\bar{\lambda}_s) \simeq (\iu \omega_s,-\iu \omega_s)$ to cancel them two-by-two. Let us define the upper asterisk $\,^*$ as the operator that selects the conjugate of a generic index, that is to say:
\begin{equation}\label{eq:conjugateindex}
i_k^* = 
\begin{cases}
i_k+n \qquad \text{if } i_k \leq n,\\
i_k-n \qquad \text{if } i_k > n,\\
\end{cases}
\end{equation}
where this choice is related to our initial ordering of eigenvalues, see Eqs.~\eqref{eq:def_complex_eigvals} and the discussion in Section~\ref{sec:EOMeigenpb}. 
Then with this notation one can easily select all the sets of indices $\Is$ that are related to a trivial resonance of order $p$, for $p$ odd only. A set $\Is$ is trivially resonant with the index $i_r$  if and only if it writes as
\begin{equation}\label{eq:setItrivial}
\Is = \lbrace i_1\; i_2\; \ldots i_{\frac{p-1}{2}}\;\; i_1^*\; i_2^*\; \ldots i_{\frac{p-1}{2}}^*\;\; i_r \rbrace.
\end{equation}
which simply states that the complex conjugates just need to cancel two-by-two in the summation defining $\sigma_\Is$, and only one remaining index is needed.
In an arbitrary order framework, one then needs to track all these sets, for all odd order $p$, and for all indices $i_r$ such that mode $r$ is a master. 

The second case is that of {\em internal resonances}, which are defined in vibration theory as a commensurability relationship between eigenfrequencies~\cite{Nayfeh79,Nayfeh00}. The low-order internal resonances are the most well known and have given rise to a vast literature investigating their solutions, see {\em e.g.}~\cite{Miles84,nayfehcarboChin99,manevitch2003free,Givois11,Gobatres12}  to cite only a few entries. Second-order internal resonances are related to quadratic nonlinearities, and give rise to the 1:2 case where $\omega_j \simeq 2\omega_p$, as well as combination of resonance such as $\omega_p \simeq \omega_k \pm \omega_l$. Third-order internal resonances encompass the 1:1 case ($\omega_j \simeq \omega_k$), the 1:3 case ($\omega_l \simeq 3\omega_j$) as well as combinations ($\omega_l = 2\omega_j \pm \omega_s$, $\omega_p = \omega_k \pm \omega_j \pm \omega_r$). For deriving a method  
which tracks all the possible high-order resonances at a generic order, one needs to identify all the sets $\Is$ associated to such resonances. This is more involved than for trivial resonances where the explicit writing of {\em all} sets $\Is$ giving rise to a trivial resonance gives directly Eq.~\eqref{eq:setItrivial}. Indeed, internal resonance will encompass all the cases where there is more than a single index remaining in the summation. One can then rewrite the summation of eigenvalues $\sigma_\Is$ defined in Eq.~\eqref{eq:defsigmaIsum} by dividing the sets of indices $\Is$ into two parts: a first subset  containing all the indices whose conjugate pair is not included in the set and a second subset containing the pairs of two conjugate indices. We can then rewrite the generic $\sigma_\Is$ as the sum of these two subsets:
\begin{equation}\label{eq:two_subsets}
\sigma_\Is  =
\sum_{i_k\in\Is:\,i_k^* \notin \Is} \lambda_{i_k} 
+
\sum_{i_l\in\Is:\,i_l^* \in \Is} \lambda_{i_l} 
\approx
\sum_{i_k\in\Is:\,i_k^* \notin \Is} \lambda_{i_k} \;,
\end{equation}
where the last simplification stems from our assumption of small damping. Indeed, when $\lambda_s \simeq \pm \iu \omega_s$, if both an eigenvalue and its conjugate are in the summation, then their sum is close to zero. An internal resonance then occurs if this summation is close to an eigenfrequency of the system:
\begin{equation}
\sum_{i_k\in\Is:\,i_k^* \notin \Is} \lambda_{i_k} 
\approx 
\lambda_r
\end{equation}
In such case, the associated $\sigma_\Is$  is resonant with the $r$-th eigenvalue. Finally two different cases can be distinguished. If $r$ is included in the list of master modes, then, following the denomination introduced in~\cite{Haro,Haller2016}, an {\em inner resonance} occurs, since the internal resonance relationship concerns only the set of master modes. On the other hand, if $r$ belongs to the set of slave modes, then an {\em outer resonance} is at hand. In this second case, an ill-conditioning of the system to be solved will appear in the subsequent calculations, since a strong coupling exists between one slave mode and some of the master modes. This means that the set of master modes is not complete as they are strongly coupled with a slave mode. Consequently the only possible strategy is to include mode $r$ in the set of master modes such that the outer resonance becomes an inner one. In the case of inner resonance, then no ill-conditioning appears and the solution strategy is discussed in the next section where the different styles of parametrisation are detailed.
\red{Theoretically speaking, the maximum order of $\sigma_\Is$ to check for possible resonances should exceed the spectral quotient as underlined in~\cite{Haller2016}; however, being the spectral quotient the ratio between the real part of the most damped mode and that of the least damped one, in the case of large FE models it can reach very high values, thus making such check unfeasible. Checking for resonances of order higher than the spectral quotient is needed from a mathematical viewpoint to ensure uniqueness but in real-world engineering applications, very high order internal resonances are difficult to appear. In the practice, it is not uncommon to a priori assume the absence of outer resonances of order higher than the parametrisation order. This is the assumption retained in this paper as well as in previous similar works dealing with FE models (\cite{artDNF2020,AndreaROM,JAIN2021How,MingwuLi2021_1,MingwuLi2021_2}).}

\section{Solutions at arbitrary order: styles of parametrisation}

This Section is devoted to detailing the different possible solutions of homological equations and thus the different {\em styles} of parametrisation. The wording adopted here follows the book by Haro {\em et al.}~\cite{Haro} and also previous discussions, 
see {\em e.g.}~\cite{Murdock} or~\cite{KahnZarmi} where the wording {\em free functions}
is used. As shown in~\cite{Haro}, there are two main styles of solution, namely the graph style and the normal form style. Here we will adapt this discussion to the case of vibratory systems. Taking the peculiarity of nonlinear vibrations into account will also lead to distinguish between complex and real normal form styles, and make a better link with previous developments on NNMs and SSMs.

\subsection{Styles of parametrisation}\label{sec:styles}

The two different styles of parametrisation can be simply understood by looking at the tangent order-$p$ homological equation expressed in~\eqref{eq:homo_tang}. Indeed a first possibility to unfold the underdeterminacy consists in setting $\p{{\theta}_{r\Is}}$ to zero and let as only unknown $\p{f_{r\Is}}$. The reason for this choice is to select the tangent part of the nonlinear mapping as simple as possible (by just cancelling it) and leave all the complexity to the nonlinear reduced dynamics $\p{f_{r\Is}}$. This option is that of the {\em graph style}. 
On the other hand, one could choose to cancel as many nonlinear terms as possible in the reduced dynamics in order to derive its simplest expression. This is common to the idea of normal transformation and it leads to the {\em normal form style}. The main idea is to cancel $\p{f_{r\Is}}$ in~\eqref{eq:homo_tang}, but this can be done only in case of no resonance. In fact, recalling Eq.~\eqref{eq:homo_tang}, when the factor $(\sigma_\Is - \lambda_r )$ is very small, i.e. when a resonance occurs, one must set $\p{{\theta}_{r\Is}}$ to zero or the problem will be ill-conditioned.
 
As discussed in the previous section, resonances are of great importance and play a specific role in nonlinear vibration theory when small damping is assumed. As a matter of fact, resonances cannot be avoided due to the occurrence of numerous trivial resonances. For the next developments, let us introduce ${\mathcal R}$ as the set of indices such that the choice of setting $\p{{\theta}_{r\Is}}$ to zero is done; for each set $\Is$, a set ${\mathcal R}$ must be defined, and the generic solutions of Eq.~\eqref{eq:homo_tang} can be rewritten as
\begin{equation}
\begin{cases}
\p{f_{r\Is}} = 0 ,\qquad (\sigma_\Is - \lambda_r ) \p{{\theta}_{r\Is}} = \p{{g}_{r\Is}}  ,
\qquad\qquad &\text{if }r\notin \Rs,
\\
\p{{\theta}_{r\Is}} = 0 ,\qquad  \p{f_{r\Is}} = \p{{g}_{r\Is}},
\qquad\qquad &\text{if }r\in \Rs.
\end{cases}
\label{eq:homo_tang_sol}
\end{equation}
As previously discussed, the second line can be selected for two different reasons. The former, 
which corresponds to the graph style choice, consists in stating that one does not want to track all the complicated resonance relationships throughout the calculations. The latter 
is that of an existing resonance relationship. Indeed in this specific case, $\sigma_\Is - \lambda_r $ vanishes and no other option is viable.

Let us now discuss in greater detail the different solutions, starting with the graph style. In this case, the second line of Eq.~\eqref{eq:homo_tang_sol} is selected for every master coordinate, $\forall r\in[1,2n]$ and for each set $\Is$. As a consequence all the associated monomials $\p{f_{r\Is}}$ can be solved for and are kept in the reduced dynamics, and all the 
master modal projections $\p{{\theta}_{r\Is}}$ are set to zero. In terms of the set $\Rs$, the graph style corresponds to selecting, for each $\Is$, $\Rs_\text{Graph}$ as
\begin{equation}\label{eq:mathcalRgraph}
\Rs_\text{Graph} = \lbrace 1\,2\,\ldots 2n \rbrace.
\end{equation}
A first remark is that the nonlinear terms of the change of coordinates are selected as simple as possible, but the price to pay is the maximal complexity of the reduced dynamics. As we will see later, this choice has for consequence that a functional relationship can be deduced between slave and master coordinates such that the master coordinates are equal to the modal ones, hence defining a graph relationship, which gives the name to this style of solution.

Turning now to the normal form style, we distinguish between two solutions, namely the complex and the real normal form. In vibration theory, complex and real normal forms have already been used in different contexts, see {\em e.g.}~\cite{Jezequel91,touze03-NNM,TOUZE:JSV:2006,NeildNF00,LamarqueUP,NeildNF01,Wagg2019}. It appears meaningful to discuss these two different strategies in the present context of the parametrisation of invariant manifold in order to synthetically present the origin of their difference. 

The general choice of the normal form style consists in selecting the first line in Eq.~\eqref{eq:homo_tang_sol} as long as no resonance occurs. In case of resonance, then the second line is selected. The main advantage resides in the fact that the simplest form of the reduced dynamics is found since all non-resonant terms are cancelled. Also, a complete nonlinear mapping is retrieved at the end of the process, including the master modes, as opposed to the graph style.  The only difference between complex and real normal form resides in the indices retained in the set ${\mathcal R}$. Since Eq.~\eqref{eq:homo_tang_sol} must be solved for each monomial and thus for each set $\Is$, it is also important to understand here that for the normal form style, there is a different ${\mathcal R}$ for each $\Is$. In case $\Is$ is not associated to a resonance relationship then ${\mathcal R}$ is empty and the first line of \eqref{eq:homo_tang_sol} is selected $\forall r\in[1,2n]$. Let us discuss now the choice when there is a resonance, be it a trivial or an internal one, for both complex and real normal form.

The complex normal form (CNF) is the choice made in mathematical textbooks~\cite{Haro,HaragusIooss,Murdock} and is also retained for example in~\cite{Jezequel91,Haller2016,JAIN2021How}. For a given set $\Is$, it consists in finding all the $r$ eigenvalues such that $\sigma_\Is  \approx \lambda_r$ is fulfilled. Assuming a set $\Is$ such that $r$ has a resonance relationship, then the associated ${\mathcal R}$ denoted as $\Rs_\text{CNF}$ (for complex normal form) contains only this $r$:
\begin{equation}
\Rs_\text{CNF} = \lbrace r \rbrace \qquad \text{with } r: \sigma_\Is  \approx \lambda_r.
\end{equation}
On the other hand, the real normal form (RNF) imposes one more condition. 
In order to come back more easily to real coordinates, and since the two complex eigenvalues are related to the same real normal mode with eigenfrequency $\omega_r$, one would like to treat similarly the two conjugates. In terms of the resonance condition, one can select both complex conjugate terms by imposing the fulfilment of its square:  
\begin{equation}
\sigma_\Is^2 \approx \lambda_r^2.
\end{equation}
Indeed, one can notice that for each value of $\sigma_\Is$, both $\lambda_r$ and its conjugate $\bar{\lambda_r}$ satisfy the condition. For a given $\Is$, such that $r$ is resonant, the set of resonant indices $\Rs$ is composed, for the real normal form style, of both $r$ and $r^*$
\begin{equation}\label{eq:mathcalR_RNF}
\Rs_\text{RNF} = \lbrace r \, r^* \rbrace \qquad \text{with } r,r^*: \sigma_\Is^2 \approx \lambda_r^2 \approx \lambda_{r^*}^2.
\end{equation}
In other words, the CNF style separates the two conjugates eigenvalues $+\iu \omega_r$ and $-\iu \omega_r$, and treats them distinctly to track the resonances. This leads to the simplest 
and most symmetric normal form. On the other hand, RNF groups the two eigenvalues $\pm \iu \omega_r$ throughout the process. One justification comes from the fact that with RNF, it is easier to go back to second-order, oscillatory-like equations. Indeed, in the case of no internal resonance, calculations show that, with RNF up to third-order, the derivative of the normal displacement is equal to the normal velocity, without extra nonlinear terms as those appearing for the CNF. The same property holds for the graph style at any order of the expansion. The link with Cartesian coordinates is then simpler and more direct to write.


At this point it is important to emphasise that the real normal form introduced here is different from the one derived in~\cite{touze03-NNM,TOUZE:JSV:2006,artDNF2020}, and appears to be closer to the one introduced in~\cite{NeildNF00,NeildNF01}. Indeed, the real normal form used in our previous derivations~\cite{touze03-NNM,TOUZE:JSV:2006,artDNF2020} has more non-empty sets~$\Is$, meaning that more terms are considered as resonant and kept in the reduced dynamics. This development was meaningful in~\cite{touze03-NNM} since it allowed to keep oscillator-like equations with real modes and coordinates throughout the calculations. However for the purpose of arbitrary order expansion, it has been chosen to discard this choice which induces an important loss of symmetry. \ref{app:oldrealNF} gives more details on how this third normal form style could be recovered from the present analysis and makes a link with other works to better unify the different approaches available in the literature.

\subsection{Solutions of the homological equation}

In this Section, we explain how to write down the general solutions directly from the physical space, for each of the three styles discussed before. As a matter of fact, all the developments starting 
from Eq.~\eqref{eq:defthetaproj} used a projection to the modal space in order to shed light on the different possible solutions to the order-$p$ homological equation, as also derived for example in~\cite{Haro}. This is meaningful since it allows to foresee the different styles, understand the appearance of the resonances, and figure out that the matrices multiplying the terms of the nonlinear mappings in Eq.~\eqref{eq:homo} might in fact be ill-conditioned because of the presence of resonance. Unfortunately an efficient computational strategy cannot rely on the projected equations since the numerical burden would increase dramatically, losing all the advantages of a direct approach. Consequently, one must now take advantage of the solutions identified and rewrite them  
directly from \eqref{eq:homo}, without resorting to the projection.

The main problem originates from the resonances, since such a relationship renders the problem~\eqref{eq:homo} ill-conditioned and the matrix multiplying the nonlinear mapping terms not invertible. One needs thus to rewrite the problem in such a way that it can be solved. To that purpose, let us use the solutions given in the previous sections in modal space. For a resonant index  $r\in\Rs$, according to the second line of Eq.~\eqref{eq:homo_tang_sol},  $\p{{\theta}_{r\Is}}$ must be set to zero, which implies, by using \eqref{eq:defthetaproj}:
\begin{equation}
\p{{\theta}_{r\Is}}
=
\begin{bmatrix}\v{X}_r\end{bmatrix}^\text{T}
\begin{bmatrix}
\M & \0
\\
\0 & \M
\end{bmatrix}
\begin{bmatrix}
\p{\WV_\Is}\\ \p{\WU_\Is}
\end{bmatrix}
= 
0.
\end{equation}
By expanding the products and simplifying, one easily arrives at
\begin{equation}\label{eq:border}
\phiv_r^\text{T}\M\p{\WV_\Is} = \bar{\lambda}_r\phiv_r^\text{T}\M \p{\WU_\Is}.
\end{equation}
This condition can now be used in order to rewrite Eq.~\eqref{eq:homo} in a solvable way, by concatenating all the unknowns in a single vector and by adding Eq.~\eqref{eq:border} to the system. Assuming for the moment a generic case where the set $\Rs$ of resonant indices contains $m$ terms: $\Rs = \lbrace r_1\,r_2\,\ldots,r_m\rbrace$; one can rewrite Eq.~\eqref{eq:homo} as
\begin{equation}\label{eq:homoaugmented}
\begin{bmatrix}
\sigma_\Is\M +\C & \K & \M\phiv_{r_1}\lambda_{r_1} & \M\phiv_{r_2}\lambda_{r_2} & \ldots & \M\phiv_{r_m}\lambda_{r_m}
\\
-\M & \sigma_\Is \M & \M\phiv_{r_1} & \M\phiv_{r_2} & \ldots & \M\phiv_{r_m}
\\
\rule{0pt}{11pt}\phiv_{r_1}^\text{T}\M & -\phiv_{r_1}^\text{T}\M\bar{\lambda}_{r_1} & 0 & 0 & \ldots & 0\\
\rule{0pt}{10pt}\phiv_{r_2}^\text{T}\M & -\phiv_{r_2}^\text{T}\M\bar{\lambda}_{r_2} & 0 & 0 & \ldots & 0\\
\vdots					& \vdots			& \vdots & \vdots &  \vdots & \vdots\\
\phiv_{r_m}^\text{T}\M & -\phiv_{r_m}^\text{T}\M\bar{\lambda}_{r_m} & 0 & 0 & \ldots & 0\\
\end{bmatrix}
\begin{bmatrix}
\p{\WV_\Is}\\ \p{\WU_\Is}\\ \p{f_{r_1\Is}}\\\p{f_{r_2\Is}}\\\vdots\\\p{f_{r_m\Is}}
\end{bmatrix}
=
\begin{bmatrix}
- \M \p{\FV_\Is} - \p{\FG_\Is} - \p{\FH_\Is}\\
- \M \p{\FU_\Is}\\
0\\0\\\vdots\\0
\end{bmatrix}.
\end{equation}
This system leverages on the knowledge of the solution types and rewrites the order-$p$ homological equation directly from the physical space in an augmented way, by adding the last $m$ lines corresponding to resonances, such that the whole problem is now well defined. It is interesting to notice that the singular matrix is bordered by the right and left eigenvectors of its kernel, thus the invertibility of the whole system is automatically fulfilled. 

The size of system \eqref{eq:homoaugmented} is $(2N+m)\times(2N+m)$, with $m$ the number of indices in $\mathcal R$, which can be up to $n$, the number of master modes in the worst case (graph style). Thanks to the splitting between displacement and velocity mappings selected at the beginning, Eq.~\eqref{eq:mappings_compact}, it is however possible to make this system symmetric and of size $(N+m)\times (N+m)$, which has very important computational consequences as it halves the size of the problem being $m < 2n \ll N$. To this purpose, one can notice that in the second row (corresponding to lines $N+1$ to $2N$), the unknown vector is everywhere multiplied by the mass matrix. Consequently, rewriting the equation corresponding to that row alone, allows one to make explicit the link between displacement and velocity mappings $\p{\WU_\Is}$ and $\p{\WV_\Is}$. The link between these two has already been underlined in numerous studies, see {\em e.g.}~\cite{PesheckBoivin,Pesheck00,artDNF2020}, but has never been presented at arbitrary order. It reads:
\begin{equation}
\p{\WV_\Is} = 
\sigma_\Is \p{\WU_\Is} + 
\sum_{r_j \in \Rs}\left(\phiv_{r_j}\p{f_{r_j \Is}}\right) + 
\p{\FU_\Is}.
\label{eq:mapping_vel}
\end{equation}
This important property is intimately related to the fact that the original system is second-order in time. Even if it has been rewritten as a first-order problem, important features of the initial formulation are preserved, such as the relationship between displacement and velocity. 
Also, the choice retained to make the problem first-order has important consequence at this stage of the procedure. Indeed, other choices might lead to a different form of the second line in \eqref{eq:homoaugmented} and Eq.~\eqref{eq:mapping_vel} might not be easily accessible, thus also preventing from halving the size of \eqref{eq:homoaugmented}.

Using Eq.~\eqref{eq:mapping_vel} allows rewriting the first row of \eqref{eq:homoaugmented} (lines 1 to $N$) as
\begin{equation}
(\sigma_\Is^2 \M +\sigma_\Is\C +\K)\p{\WU_\Is}+
\sum_{r_j \in \Rs}\left((\sigma_\Is+\lambda_{r_j})\M +\C\right)\phiv_{r_j}\p{f_{r_j \Is}}
=
\p{\RHS_\Is},
\label{eq:first_row_pre}
\end{equation}
where the vector $\p{\RHS_\Is}$ has been introduced to collect all known terms into a single right-hand side; it is defined as:
\begin{equation}
\p{\RHS_\Is} = 
- \p{\FG_\Is} 
- \p{\FH_\Is}
- \M \p{\FV_\Is}
- (\sigma_\Is \M +\C) \p{\FU_\Is}.
\label{eq:RHS}
\end{equation}
Eq.~\eqref{eq:first_row_pre} can be further simplified thanks to the following equation:
\begin{equation}\label{eq:symm_of_syst}
\left((\sigma_\Is+\lambda_{r_j})\M +\C\right)\phiv_{r_j} = 
(\sigma_\Is-\bar{\lambda}_{r_j})\M\phiv_{r_j}.
\end{equation}
This property stems naturally from the eigenproblem and is shown in  \ref{app:eigs}. Finally the system derived from the order-$p$ homological equation, the solutions of which will give the unknown mappings and reduced dynamics at order $p$, can be rewritten as
\begin{equation}
(\sigma_\Is^2 \M +\sigma_\Is\C +\K)\p{\WU_\Is}+
\sum_{r_j \in \Rs}(\sigma_\Is - \bar{\lambda}_{r_j})\M\phiv_{r_j}\p{f_{r_j \Is}}
=
\p{\RHS_\Is}.
\label{eq:first_row}
\end{equation}
In order to make the system complete, a final simplification on the last $m$ rows of \eqref{eq:homoaugmented} can be performed. Let us write one of these rows with generic index $r_k$:
\begin{equation}
\phiv_{r_k}^\text{T}\M \left(\sigma_\Is \p{\WU_\Is} +
\sum_{r_j \in \Rs}\left(\phiv_{r_j}\p{f_{r_j \Is}}\right) +
\p{\FU_\Is}\right)
-\phiv_{r_k}^\text{T}\M\bar{\lambda}_{r_k} \p{\WU_\Is} =0.
\end{equation}
Thanks to the orthonormality property of the eigenmodes,  $\phiv_{r_k}^\text{T}\M\phiv_{r_j}=1$ if $r_k$ is equal to $r_j$ or to $r_j^*$, and $\phiv_{r_k}^\text{T}\M\phiv_{r_j}=0$ otherwise. Using this property, the row associated to $r_k$ finally  reads:
\begin{equation}\label{eq:augmentedrowslast}
\begin{cases}
(\sigma_\Is -\bar{\lambda}_{r_k})\phiv_{r_k}^\text{T}\M\p{\WU_\Is} +
\p{f_{r_k \Is}} =
-\phiv_{r_k}^\text{T}\M\p{\FU_\Is},	&\qquad\text{if }r_k^*\notin \Rs,
\\
(\sigma_\Is -\bar{\lambda}_{r_k})\phiv_{r_k}^\text{T}\M\p{\WU_\Is} +
\p{f_{r_k \Is}} + \p{f_{r_k^* \Is}}=
-\phiv_{r_k}^\text{T}\M\p{\FU_\Is},	&\qquad\text{if }r_k^*\in \Rs.
\end{cases}
\end{equation}
To conclude, the order-$p$ homological equation can be rewritten as a $(N+m)\times (N+m)$ symmetric problem given by both Eqs.~\eqref{eq:first_row} and \eqref{eq:augmentedrowslast}, where $m$ is the cardinal number of $\mathcal R$. As compared to Eq.~\eqref{eq:homoaugmented}, the size of the system has been divided by two which has important consequences for computational speed-up. In order to derive the explicit solutions, one must now discuss, as a final step, 
how each style of parametrisation will affect the problem to solve. In order to make this last discussion easier, we focus on the simple case of a system free of internal resonances. 
Consequently, only trivial resonances are present in the system and need to be tracked by the method.

\subsection{Systems free of internal resonances}

To better understand this simplified case, let us recall Eq.~\eqref{eq:two_subsets}, where, under the small damping assumption, the sum of two complex conjugate eigenvalues was neglected; if both indices of a complex conjugate pair are included in the generic set $\Is$ then $\sigma_\Is$ simplifies to the sum of the sole eigenvalues whose conjugate is not included in $\Is$:
\begin{equation}
\sigma_\Is 
\approx
\sum_{i_k\in\Is:\,i_k^* \notin \Is} \lambda_{i_k} \;.
\end{equation}

In the case of internally resonant systems, not only $\sigma_\Is$ can be close to an eigenvalue for a larger number of sets $\Is$ but also it can be close to multiple eigenvalues. Let us suppose for example that a two-modes system is investigated and their frequencies satisfy the relationship $\omega_1 \approx \omega_2$. Then the set $\Is = \lbrace 1 1^* 1 \rbrace$ is not only trivially resonant with $ r = 1$ but also with $r = 2$, due to the presence of a 1:1 internal resonance between the two modes. The resonant sets for the normal form styles will be in this case $\Rs_\text{CNF} = \lbrace 1,2\rbrace$ and $\Rs_\text{RNF} = \lbrace 1,1^*,2,2^*\rbrace$; only the graph style will not be affected by the presence of internal resonances. Moreover, the number of sets $\Is$ that will be resonant is obviously larger when internal resonances occur. In the scenario of 1:1 internal resonance between two modes, not only the set $\Is = \lbrace 1 1^* 1\rbrace$ is resonant with the indices $ r = 1$ and $r = 2$, but also the sets $\Is = \lbrace 1 1^* 2\rbrace$, $\Is = \lbrace 2 2^* 1\rbrace$, $\Is = \lbrace 2 2^* 2\rbrace$,
$\Is = \lbrace 1 2^* 1\rbrace$, $\Is = \lbrace 1 2^* 2\rbrace$, $\Is = \lbrace 2 1^* 1\rbrace$, and $\Is = \lbrace 2 1^* 2\rbrace$.

Conversely, if there are no internal resonances, the distinction between resonant and non-resonant sets $\Is$, together with the individuation of the eigenvalues they resonate with, is much simpler, thus allowing to give explicit expressions. In the case of two-modes reduction at third order, only the set $\Is = \lbrace 1 1^* 1\rbrace$ would resonate with $r=1$ and the set $\Is = \lbrace 2 2^* 2\rbrace$ with $r=2$. In fact, any set $\Is$ is either non-resonant or resonant with one eigenvalue because the value $\sum_{i_k\in\Is:\,i_k^* \notin \Is} \lambda_{i_k} $ is either zero or equal to a single eigenvalue $\lambda_{i_k}$. 
It follows that, the normal form styles will have either an empty set $\Rs$ or a set composed of a single index in the case of CNF, and two indices in the case of RNF. Furthermore, only the odd order sets can be resonant, whereas any even order set will not be.

Before moving to the expressions of the homological solutions for the different styles in the case of systems with no internal resonances, it is worth highlighting that, from a computational point of view, the treatment of the resonances is the same for both internally resonant and non-internally resonant systems, the only difference being the individuation of the set $\Rs$ for each $\Is$. However, from a presentation point of view, it is much easier to restrict ourselves to the case of systems free of internal resonances, as the expressions provided in the following are much more readable than those of the general case.

Now the solutions to Eqs.~\eqref{eq:first_row} and \eqref{eq:augmentedrowslast} are detailed for each style of parametrisation. Indeed, according to the discussion led in Section~\ref{sec:styles}, the set $\mathcal R$ is differently filled out for each style. Let us begin with the complex normal form style which leads to the smallest set $\mathcal R$.

In the normal form style, one needs to distinguish the resonant case from the non-resonant one. In the non-resonant case, the set  $\Is$ is not resonant with any $r$. Consequently, $\p{f_{r \Is}} = 0$ and the system reduces to 
\begin{equation}\label{eq:cnfnonreso}
\left( \sigma_\Is^2\M + \sigma_\Is \C + \K \right) \p{\WU_\Is} = \p{\RHS_\Is}.
\end{equation}
Solving \eqref{eq:cnfnonreso} allows finding the nonlinear mapping term $\p{\WU_\Is}$, which, together with $\p{f_{r \Is}} = 0$ and Eq.~\eqref{eq:mapping_vel}, gives the full solutions for the monomial associated to $\Is$ in the complex normal form style.

Otherwise, if the set $\Is$ is resonant with $r$, then the corresponding monomial of the reduced dynamics cannot be cancelled. The system is composed of two terms, stemming respectively from \eqref{eq:first_row} and \eqref{eq:augmentedrowslast}, and reads:
\begin{equation}
\begin{bmatrix}
\sigma_\Is^2\M + \sigma_\Is \C + \K & (\sigma_\Is-\bar{\lambda}_r)\M\phiv_r
\\
(\sigma_\Is - \bar{\lambda}_r) \phiv_r^\text{T}\M & 1
\end{bmatrix}
\begin{bmatrix}
\p{\WU_\Is}\\ \p{f_{r\Is}}
\end{bmatrix}
=
\begin{bmatrix}
\p{\RHS_\Is}\\
- \phiv_r^\text{T}\M \p{\FU_\Is} 
\end{bmatrix}\;.
\label{eq:homo_reduced_cnf}
\end{equation}

As mentioned before, the matrix $\sigma_\Is^2\M + \sigma_\Is \C + \K$ is singular because in this case $\sigma_\Is\approx\lambda_r$; however, the whole system is invertible thanks to the bordering of the singular matrix with the eigenvector of its kernel $\phiv_r$. In the case of internally resonant systems, some sets $\Is$ could be resonant with more than one eigenvalue, and in that case the bordering would consists in more columns and rows instead of just one.

The coefficient of the reduced dynamics can be made explicit thanks to the relationship $\p{f_{r\Is}} = \p{{g}_{r\Is}}$, see Eq.~\eqref{eq:homo_tang_sol}. Using Eqs.~\eqref{eq:gtilde} together with the definition of the complex left eigenvectors $\v{X}_j$ given in Eq.~\eqref{eq:def_complex_eig_left} leads to:
\begin{equation}
\p{f_{r\Is}} =
\dfrac{
\phiv_r^\text{T} (
- \M \p{\FV_\Is} - \p{\FG_\Is} - \p{\FH_\Is}
+ \lambda_{r^*} \M \p{\FU_\Is})
}{\lambda_r - \lambda_{r^*}},
\end{equation}
which is an explicit expression of the coefficient of the resonant monomial at arbitrary order in the complex normal form style.


We now turn to the real normal form style where the set $\mathcal R$, when not empty, is composed of two conjugated indices, see Eq.~\eqref{eq:mathcalR_RNF}. Since it is also a normal form style, one still needs to distinguish between resonant and non-resonant cases. If the set $\Is$ is not a resonant one, then the same solution as for the CNF style applies and the system reduces to the same equation:
\begin{equation}
\left( \sigma_\Is^2\M + \sigma_\Is \C + \K \right) \p{\WU_\Is} = \p{\RHS_\Is}.
\end{equation}

Let us consider the case of a set $\Is$ which is resonant with $r$. Since $\Rs_\text{RNF} = \lbrace r \, r^* \rbrace$, then two lines needs to be considered in Eq.~\eqref{eq:homoaugmented}, which can be now rewritten as:
\begin{equation}
\begin{bmatrix}
\sigma_\Is^2\M + \sigma_\Is \C + \K & 
(\sigma_\Is-\bar{\lambda}_r)\M\phiv_r &
(\sigma_\Is-\lambda_r)\M\phiv_r
\\
(\sigma_\Is - \bar{\lambda}_r) \phiv_r^\text{T}\M & 1 & 1
\\
(\sigma_\Is - \lambda_r) \phiv_r^\text{T}\M & 1 & 1
\end{bmatrix}
\begin{bmatrix}
\p{\WU_\Is}\\ \p{f_{r\Is}}\\ \p{f_{r^*\Is}}
\end{bmatrix}
=
\begin{bmatrix}
\p{\RHS_\Is}\\
- \phiv_r^\text{T}\M \p{\FU_\Is} \\
- \phiv_r^\text{T}\M \p{\FU_\Is} 
\end{bmatrix}.
\label{eq:homo_reduced_rnf}
\end{equation}
The main difference with the complex normal form style can be highlighted by properly interpreting the last two rows of  Eq.~\eqref{eq:homo_reduced_rnf} in the case of a resonant $\Is$. Indeed, thanks to the last row and the presence of the conjugate, not taken into account in the complex normal form, the last two rows are equivalent to imposing the orthogonality of the mapping with the $r$-th real mode: $\phiv_r^\text{T}\M \p{\WU_\Is}=0$. To demonstrate this property, one needs to take the difference and the sum of the last two lines of \eqref{eq:homo_reduced_rnf}, yielding:
\begin{subequations}\begin{align}
&(2\sigma_\Is - \bar{\lambda}_r - \lambda_r) \phiv_r^\text{T}\M \p{\WU_\Is} + 
 2\p{f_{r\Is}}+2\p{f_{r^*\Is}} = - 2\phiv_r^\text{T}\M \p{\FU_\Is},
\\
&(\bar{\lambda}_r - \lambda_r) \phiv_r^\text{T}\M \p{\WU_\Is}= 0,
\end{align}\end{subequations}
which simplifies to:
\begin{subequations}\begin{align}
&\p{f_{r\Is}}+\p{f_{r^*\Is}} = - \phiv_r^\text{T}\M \p{\FU_\Is}, \label{eq:ff_mu}
\\
&\phiv_r^\text{T}\M \p{\WU_\Is}= 0,
\end{align}\end{subequations}
hence showing the property for the displacement mapping $\p{\WU_\Is}$. The same orthogonality condition holds also for the velocity mapping; using Eq.~\eqref{eq:mapping_vel}, one obtains:
\begin{equation}
\phiv_r^\text{T}\M \p{\WV_\Is}= 0.
\end{equation}
To conclude with the real normal form style, we now provide the explicit expressions for the coefficients of the resonant monomials. In the real normal form style, since two indices are contained in $\mathcal R$, two coefficients are kept for each resonant monomial. Using   $\p{f_{r\Is}} = \p{{g}_{r\Is}}$ and $\p{f_{r^*\Is}} = \p{{g}_{r^*\Is}}$ together with  Eq.~\eqref{eq:gtilde} leads to:
\begin{subequations}\begin{align}
\p{f_{r\Is}} =
\dfrac{
\phiv_r^\text{T} (
- \M \p{\FV_\Is} - \p{\FG_\Is} - \p{\FH_\Is}
+ \lambda_{r^*} \M \p{\FU_\Is})
}{\lambda_r - \lambda_{r^*}},
\\
\p{f_{r^*\Is}} =
\dfrac{
\phiv_r^\text{T} (
- \M \p{\FV_\Is} - \p{\FG_\Is} - \p{\FH_\Is}
+ \lambda_r \M \p{\FU_\Is})
}{\lambda_{r^*} - \lambda_r}.
\end{align}\end{subequations}
Importantly, the reduced dynamics term $\p{f_{r\Is}}$ has the same expression as in the complex style but its value is different since $\p{\FV_\Is}$ and $\p{\FU_\Is}$ are different for complex and real normal form styles, as they depends on lower order mappings and reduced dynamics coefficients.

Let us now conclude the Section by giving the solutions in the case of the graph style, where the second line in Eqs.~\eqref{eq:homo_tang_sol} is always taken. For every set $\Is$ related to each monomial, $\mathcal R$ is fully populated with all the indices of the master modes, as stated in Eq.~\eqref{eq:mathcalRgraph}. So, for every set $\Is$, the system to be solved, composed of Eq.~\eqref{eq:homoaugmented} with its subsequent simplifications, now reads:
\begin{equation}
\begin{bmatrix}
\sigma_\Is^2\M + \sigma_\Is \C + \K & 
\M\v{\Phi}\; (\Id_{[2n]} \sigma_\Is - \bar{\v{\Lambda}})
\\
\rule{0pt}{20pt}
\left(\M\v{\Phi} \;(\Id_{[2n]} \sigma_\Is - \bar{\v{\Lambda}})\right)^\text{T} & 
\begin{bmatrix}
\Id_{[n]} & \Id_{[n]} \\ \Id_{[n]} & \Id_{[n]}
\end{bmatrix}
\end{bmatrix}
\begin{bmatrix}
\p{\WU_\Is}\\ \rule{0pt}{15pt}\p{\v{f}_{\Is}}
\end{bmatrix}
=
\begin{bmatrix}
\p{\RHS_\Is}\\
\rule{0pt}{15pt}- \v{\Phi}^\text{T}\M \p{\FU_\Is}
\end{bmatrix}\;,
\label{eq:homo_reduced_cnf_GS}
\end{equation}
where the matrix of master eigenvectors $\v{\Phi}$ defined in Eq.~\eqref{eq:master_eigs} and the matrix of master eigenvalues $\v{\Lambda}$ defined in Eq.~\eqref{eq:master_eigv} have been used for compactness; also the $2n\times 2n$ and $n\times n$ identity matrices, 
$\Id_{[2n]}$ and $\Id_{[n]}$ respectively, have been introduced.

A further development is needed to understand an important property of the graph style, which gives the name to the solution. Let us first remark that  both mapping vectors, for displacement and velocity $\p{\WU_\Is}$ and $\p{\WV_\Is}$, are orthogonal to each of the real master modes $\phiv_{r}$, for $r\in[1,n]$. The proof of this property strictly follows the same lines as in the real normal form case, since it is the grouping of the two conjugate indices that allows fulfilling the property. Unlike the real normal form style, in this case the property holds  $\forall r\in[1,n]$. One can then write:
\begin{subequations}\begin{align}\label{eq:graph_orthog}
&\v{\Phi}^\text{T}\M \p{\WU_\Is} = 0,
\\
&\v{\Phi}^\text{T}\M \p{\WV_\Is} = 0.
\end{align}\end{subequations}
This means that the nonlinear terms of the mappings are orthogonal to the linear subspace spanned by the master eigenvectors. Due to  the cancellation of the nonlinear terms of the mappings that stems from the choice made in Eq.~\eqref{eq:homo_tang_sol}, $\p{{\theta}_{r\Is}} = 0$,  the nonlinear mappings defined in Eq.~\eqref{eq:mappings_byorder} have no nonlinear terms in the modal space for the master coordinates. This is easily demonstrated by   projecting the whole change of coordinates $\U = \WUfun$ and $\V = \WVfun$ on a given real master mode $r$, thus retrieving the modal displacement  ${\ur}_r $ and modal velocity ${\vr}_r$ of the $r$-th mode defined in Eq.\eqref{eq:modal_displ} and Eq.~\eqref{eq:modal_vel}, which simply writes:
\begin{subequations}\label{eq:projgraphstylemod}
\begin{align}
&{\ur}_r = \phiv_{r}^\text{T}\M \U = \phiv_{r}^\text{T}\M (\v{\Phi} \v{z})= z_r + \bar{z}_r,
\\
&{\vr}_r = \phiv_{r}^\text{T}\M \V = \phiv_{r}^\text{T}\M (\v{\Phi}\v{\Lambda} \v{z})= \lambda_r z_r + \bar{\lambda}_r\bar{z}_r .
\end{align}
\end{subequations}
Due to the orthogonality property of Eq.~\eqref{eq:graph_orthog}, all the nonlinear mappings are orthogonal to the master modes, and only the linear terms are not. Interestingly, the left-hand side of these equations, which represents the modal coordinates, is linearly related to the right-hand side, which represents the normal coordinates. In other words, Eqs.~\eqref{eq:projgraphstylemod} shows that with the graph style parametrisation, $z_r + \bar{z}_r$ represents the modal displacement of mode $r$, and $\lambda_r z_r + \bar{\lambda}_r\bar{z}_r$ its velocity.
Consequently the nonlinear mapping actually operates solely on the slave modes, whereas the modal coordinates are left unchanged. Note that this result does not hold with the normal form styles, where the projections to modal space using left-hand sides of \eqref{eq:projgraphstylemod} would show a nonlinear relationship.
Finally, the last lines of the nonlinear mapping \eqref{eq:mappings_compact}, corresponding to the slave modal coordinates, can be interpreted as a functional relationship, or a graph, between slave and master modal coordinates. This explains the name selected for the graph style, and allows one to understand that this method is exactly equivalent to the choice proposed by Shaw and Pierre in their first derivations of nonlinear normal modes as invariant manifolds, see {\em e.g.}~\cite{ShawPierre91,ShawPierre93,ROMGEOMNL}. Since a graph relationship is assumed, 
the graph style is theoretically limited up to a point where the invariant manifold folds. This is in contrast with the normal form style, which proposes a completely nonlinear relationship and is thus a priori able to pass over the folding. This property will be illustrated in Section~\ref{sec:cantilever}.


\section{Complex to real coordinates}\label{sec:compreal}

In the previous Section important improvements in terms of computational efficiency have been achieved recalling that the initial system is of second-order in time. Consequently the velocity mapping can be expressed directly as a function of the displacement mapping, Eq.~\eqref{eq:mapping_vel}, which has been used to halve the size of system \eqref{eq:homoaugmented}. In the same spirit, the fact that the initial problem is real can also be used in order to gain computational time and memory, and  efficacy in the output processing. Indeed, since the initial problem is real, the final problems should also be real. The introduction of the complex number is an important tool which is helpful in order to better highlight the symmetries of the underlying problem, but at the end of the process, one should be able to come back to real quantities and see the complexification as a side help for conducting the inner calculations. This process is called complexification/realification and is commented in a general framework for example in~\cite{Haro}.


Let us first begin by noticing that the initial displacement and velocity vectors $\U$ and $\V$ are real, they thus must fulfil the relationships $\bar{\U} = \U$, and $\bar{\V} = \V$. Using the asymptotic expansion defining the mappings as given in Eq.~\eqref{eq:mappings_I}, one can  write:
\begin{subequations}\label{eq:UbarUVbarV}
\begin{align}
&\U= 
\sum_{p=1}^o\sum_\Is \p{\WU_\Is} \; \p{\pi_\Is}=
\bar{\U}=
\sum_{p=1}^o\sum_\Is \p{\bar{\WU}_\Is} \; \p{\bar{\pi}_\Is}=
\sum_{p=1}^o\sum_{\Is^*} \p{\bar{\WU}_{\Is^*}} \; \p{\pi_\Is},
\\
&\V= 
\sum_{p=1}^o\sum_\Is \p{\WV_\Is} \; \p{\pi_\Is}=
\bar{\V}=
\sum_{p=1}^o\sum_\Is \p{\bar{\WV}_\Is} \; \p{\bar{\pi}_\Is}=
\sum_{p=1}^o\sum_{\Is^*} \p{\bar{\WV}_{\Is^*}} \; \p{\pi_\Is}.
\end{align}
\end{subequations}
In these equations, $\Is^*$ has been introduced as the conjugate of the set of indices $\Is$, which is simply obtained by conjugating each of the indices of  $\Is$ following the rule defined in Eq.~\eqref{eq:conjugateindex}. The last identity has been written by using the fact that since the summation over the monomials spans all possible sets  $\Is$, it can be equivalently interpreted as a sum over $\Is^*$.  Using the definition of the generic monomial $\p{\pi_\Is}$ in  Eq.~\eqref{eq:defpiI}, one immediately has $ \p{\bar{\pi}_{\Is^*}}=\p{\pi_\Is}$. Term-by-term identification of the first and last summations in \eqref{eq:UbarUVbarV} shows that the following relationship must hold for the coefficients of the mappings:
\begin{subequations}
\begin{align}
    &\p{\bar{\WU}_{\Is^*}} = \p{{\WU}_{\Is}},
    \\
    &\p{\bar{\WV}_{\Is^*}} = \p{{\WV}_{\Is}}.
\end{align}
\label{eq:symm_complex}
\end{subequations}

A similar property can also be deduced for each of the coefficients of the monomial terms in the reduced dynamics. Indeed, the normal coordinates being complex conjugate, they have of course to verify  $\bar{z}_{s^*} = z_s$, but also $\dot{\bar{z}}_{s^*} = \dot{z}_s$. Using the asymptotic expansion used for the reduced dynamics, one can thus write the following equalities:
\begin{equation}
    \dot{\bar{z}}_{s^*} = \sum_{p=1}^o\sum_{\Is^*} \p{\bar{f}_{s^*\Is^*}} \; \p{\bar{\pi}_{\Is^*}} 
    = \dot{z}_s = \sum_{p=1}^o\sum_{\Is} \p{f_{s\Is}} \; \p{{\pi}_\Is} .
\end{equation}
Term-by-term identification then leads to the following property:
\begin{equation}\label{eq:symm_complex_f}
    \p{\bar{f}_{s^*\Is^*}} = \p{f_{s\Is}}.
\end{equation}

It could be shown that both Eqs.~\eqref{eq:symm_complex} and~\eqref{eq:symm_complex_f} are verified if $\Rs(\Is) = (\Rs(\Is^*))^*$, which means that the set of indices considered resonant with the set $\Is$ is equal to the conjugate of the set of indices considered resonant with the set $\Is^*$. Since a symmetric treatment of resonances with respect to their conjugate is logical, this is a very general case and it holds for all the styles mentioned here.

Thanks to the above properties, one is now in the position to write all the needed quantities, nonlinear mappings and reduced dynamics, in real coordinates. This will have some implications on the computational aspects in terms of burden and memory requirements, but also on the post-processing of the results. Indeed, providing real results for the reduced-order dynamics is much more comfortable since ROMs are generally aimed at being used for either direct time integration or more generally for interfacing with a numerical continuation method for analysing the bifurcation scenario and predict the vibratory solutions of the structure. In this context, real quantities are needed as input to continuation codes.

%

Since complex quantities are still included inside $\p{{\pi}_\Is}$, new real variables have to be introduced: $(\ar_j,\ar_{j+n})$, which correspond the Cartesian representation of the normal coordinates $\v{z}$. 

We define these purely real coordinates as twice the real and imaginary part of the complex conjugate pair $(z_j,z_{j+n})$:
\begin{subequations}\label{eq:cartes_def}\begin{align}
&\ar_j = z_j + z_{j+n} = 2\;\Re{z_j},
\\
&\ar_{j+n} = \dfrac{z_j - z_{j+n}}{\iu} = 2\;\Im{z_j},
\end{align}\end{subequations}
with $j\in[1,n]$.

Since the expression of $\Re{\p{\pi_\Is}}$ and $\Im{\p{\pi_\Is}}$ as a function of $\bm{\ar}$ are cumbersome, they are not reported here. They are however easy to compute automatically. 
Both the reduced model and the change of coordinates can be written in terms of $\bm{\ar}$, which is what we output and solve for. Let us define the generic $p$ order product of the Cartesian normal coordinates as:
\begin{equation}
\p{\tilde{\pi}_{{\Is}}} = \ar_{i_1} \ar_{1_2} \ldots \ar_{i_p},\qquad 
{\Is} = \lbrace i_1 i_2 \ldots i_p \rbrace,
\end{equation}
then it is easy to see that this monomial $\p{\tilde{\pi}_{{\Is}}}$ in $\ar_i$ can be expressed as a linear combination of some monomials $\p{\pi_\Is}$ in $z_i$ and vice versa.
For the nonlinear change of coordinates, realification is obtained by
\begin{subequations}\label{eq:real_nlmap}
\begin{align}
&\U= 
\sum_{p=1}^o\sum_\Is \p{\WU_\Is} \; \p{\pi_\Is}=
\sum_{p=1}^o\sum_{{\Is}} \p{\tilde{\WU}_{{\Is}}} \; \p{\tilde{\pi}_{{\Is}}},
\\
&\V= 
\sum_{p=1}^o\sum_\Is \p{\WV_\Is} \; \p{\pi_\Is}=
\sum_{p=1}^o\sum_{{\Is}} \p{\tilde{\WV}_{{\Is}}} \; \p{\tilde{\pi}_{{\Is}}},
\end{align}
\end{subequations}
where the nonlinear mapping tensors in Cartesian coordinates, $\p{\tilde{\WU}_{{\Is}}}$ and $\p{\tilde{\WV}_{{\Is}}}$, are now purely real quantities.

Similarly, for the reduced dynamics, realification leads to
\begin{equation}\label{eq:real_reddyn}
\dot{\ar}_s = \sum_{p=1}^o\sum_{{\Is}} \p{\tilde{f}_{s{\Is}}} \; \p{\tilde{\pi}_{{\Is}}},\qquad s\in[1,2n],
\end{equation}
where the reduced dynamics coefficients in Cartesian coordinates, $\p{\tilde{f}_{s{\Is}}} $, are purely real and are obtained by imposing the following equalities:
\begin{subequations}\begin{align}
&\Re{\dot{z}_j} = \sum_{p=1}^o\sum_{\Is} \Re{\p{f_{s\Is}} \p{{\pi}_\Is}}
=
\dfrac{1}{2} \dot{\ar}_j = 
\sum_{p=1}^o\sum_{{\Is}} \dfrac{1}{2} \p{\tilde{f}_{j{\Is}}} \; \p{\tilde{\pi}_{{\Is}}},\qquad j\in[1,n],
\\
&\Im{\dot{z}_s} = \sum_{p=1}^o\sum_{\Is} \Im{\p{f_{s\Is}} \p{{\pi}_\Is}}
=
\dfrac{1}{2} \dot{\ar}_{j+n} = 
\sum_{p=1}^o\sum_{{\Is}} \dfrac{1}{2} \p{\tilde{f}_{j+n{\Is}}} \; \p{\tilde{\pi}_{{\Is}}},\qquad j\in[1,n].
\end{align}\end{subequations}

It is worth mentioning that another possible choice for realification consists in polar coordinates; one could in fact express the complex normal coordinates as:
\begin{subequations}\label{eq:cnf_reddyn_polar}\begin{align}
&z_j = \dfrac{1}{2}\rho_j \e^{+\iu \alpha_j},
\\
&z_{j+n} = \dfrac{1}{2}\rho_j \e^{-\iu \alpha_j},\qquad\forall j\in[1,n].
\end{align}\end{subequations}
where the scaling factor $1/2$ has been chosen coherently with the choice made for the Cartesian coordinates that have been defined to be twice the real and imaginary part of the complex coordinates. With this choice of the scaling factor, the relationship between polar and Cartesian coordinates would read:
\begin{subequations}\label{eq:cnf_reddyn_carttopolar}\begin{align}
&\ar_j = \rho_j \cos(\alpha_j)\\
&\ar_{j+n} = \rho_j \sin(\alpha_j),\qquad\forall j\in[1,n]
\end{align}\end{subequations}

In the case of single mode reduction, the polar form can be particularly attractive if the complex normal form style is used because the nonlinear frequency would be directly given by the reduced dynamics expressed in polar coordinates, as it will be shown in the next section; however, this advantage does not extend to the general case of internally resonant systems. For the sake of generality, we choose here to express the reduced dynamics in Cartesian coordinates; in fact, in the general case of internally resonant systems, the reduced dynamics given by Eq.~\eqref{eq:real_reddyn} is better suited for continuation algorithms than its equivalent polar form. Moreover, Cartesian coordinates bear more resemblance with the choice made in previous works~\cite{touze03-NNM,TOUZE:JSV:2006,artDNF2020,AndreaROM}.

\section{Summary of the main results and illustrative examples}

The parametrisation method, in one of the three different styles presented,
allows one to compute both the nonlinear mappings and the reduced-order dynamics at a generic order of expansion $o$. The main difficulty lies in the tracking of all monomials of order $p$ and of the resonance conditions that need to be treated with care. The coefficients of the unknowns at order $p$ depend only on previous calculations at lower order, such that the computation needs to keep track of these incoming terms, which takes the main part of the analysis. 

Previous works already considered similar developments. For example, in ~\cite{PesheckBoivin} the invariant manifold approach was used to propose order-3 developments that are here recovered with the graph style subcase. 
Third-order expansions~\cite{touze03-NNM,TOUZE:JSV:2006} or even higher-order~\cite{LEUNG1998,LEUNG1998b,LamarqueUP} were also obtained applying the normal form approach~\cite{Jezequel91}. The complex versions \cite{Jezequel91,LEUNG1998} are recovered by the CNF, while the RNF allows retrieving the real approach developed in~\cite{NeildNF00,NeildNF01}. The real formulation proposed in~\cite{touze03-NNM,TOUZE:JSV:2006} represents a different parametrisation which is not further investigated in this article, but \ref{app:oldrealNF} collects some computations allowing to retrieve this style. The approach proposed in this article generalises these prior developments, using the same framework and offering arbitrary order expansions. Whereas previous contributions focused on giving analytical expressions for all the coefficients as function of the input (see {\em e.g.}~\cite{PesheckBoivin,Jezequel91,touze03-NNM}), 
this objective is left aside here for the sake of efficiency. Instead, the arbitrary order expansion provides an automatic reasoning in order to compute numerically all the coefficients while never searching for their analytical expressions. This is a purely numerical approach, different from previous developments, relying on a more symbolic representation.

The technique as presented in this paper uses the parametrisation method for vibratory systems and is thus very close to the developments shown in~\cite{JAIN2021How,MingwuLi2021_1,MingwuLi2021_2}. However, as shown during the previous sections, numerous distinctive features render the approach presented here interesting, by highlighting how computational efficiency can be improved by using the fact that the initial system is second-order. Also, three different styles are derived, allowing a better understanding of previous works.

To give more insight into this last point, let us show how the three different styles are treating differently the reduced dynamics for a simple case: a generic system with a single master mode with development only up to order three.
In the case of the complex normal form style, the reduced dynamics reads:
\begin{subequations}\label{eq:cnf_reddyn_singlemode}\begin{align}
&\dot{z}_1 = \lambda_1 z_1 + 
\p[3]{f_{1\lbrace 1 1 1^* \rbrace}} z_{1}^2\bar{z}_1 ,
\\
&\dot{\bar{z}}_1 = \bar{\lambda}_1 \bar{z}_1+ 
\p[3]{f_{1^*\lbrace 1 1^* 1^* \rbrace}} z_{1}\bar{z}_1^2 .
\end{align}\end{subequations}
In fact, as mentioned in Section~\ref{sec:styles}, the only resonant sets are $\Is = \lbrace 1 1 1^* \rbrace$, which resonates with $\Rs = \lbrace 1 \rbrace$ and the set $\Is = \lbrace 1 1^* 1^* \rbrace$, which resonates with $\Rs = \lbrace 1^* \rbrace$. It follows that the only nonzero coefficients $f_{r\Is}$ of the reduced dynamics are those appearing in equations~\eqref{eq:cnf_reddyn_singlemode}. A very peculiar property of this method can be deduced by noticing that in the equation for $\dot{z}_1$ only odd order terms in the form $z_{1}(z_{1}\bar{z}_1)^m$ with $m$ integer will appear at the generic order; similarly, in the equation for $\dot{\bar{z}}_1$ only odd order terms in the form $\bar{z}_1(z_{1}\bar{z}_1)^m$ will appear. This property has an important consequence that holds in the case of single mode reduction and that can be observed by expressing the normal coordinates in polar form, as shown in~\cite{Haller2016,PONSIOEN2018}. 

Plugging Eqs.~\eqref{eq:cnf_reddyn_polar} into Eqs.~\eqref{eq:cnf_reddyn_singlemode} yields:
\begin{subequations}\begin{align}
&(\dot{\rho}+\iu \dot{\alpha}\rho)\e^{+\iu \alpha}  = 
(\lambda_1 \rho + \dfrac{1}{4}\p[3]{f_{1\lbrace 1 1 1^* \rbrace}}\rho^3)\e^{+\iu \alpha} ,
\\
&(\dot{\rho}-\iu \dot{\alpha}\rho)\e^{-\iu \alpha}  = 
(\bar{\lambda}_1 \rho + \dfrac{1}{4}\p[3]{f_{1^*\lbrace 1 1^* 1^* \rbrace}}\rho^3)\e^{-\iu \alpha} .
\end{align}\end{subequations}

Due to the form of the monomials in the reduced dynamics, $z_{1}(z_{1}\bar{z}_1)^m$ in the equation for $\dot{z}_1$ and $\bar{z}_1(z_{1}\bar{z}_1)^m$ in the equation for $\dot{\bar{z}}_1$, the exponential terms $\e^{+\iu \alpha}$ and $\e^{-\iu \alpha}$ can be collected and therefore eliminated from the equations. As anticipated, one could then directly obtain the values of $\dot{\rho}$ and $\dot{\alpha}$, which represent the decay ratio and frequency of the normal coordinates dynamics, as a function of the sole amplitude $\rho$. They read:
\begin{subequations}\label{eq:rho_alpha}\begin{align}
&\dot{\rho}= 
\dfrac{\lambda_1+\bar{\lambda}_1}{2} \rho +
\dfrac{\p[3]{f_{1\lbrace 1 1 1^* \rbrace}}+\p[3]{f_{1^*\lbrace 1 1^* 1^* \rbrace}}}{8}\rho^3 ,
\\
&\dot{\alpha} = 
\dfrac{\lambda_1-\bar{\lambda}_1}{2\iu} +
\dfrac{\p[3]{f_{1\lbrace 1 1 1^* \rbrace}}-\p[3]{f_{1^*\lbrace 1 1^* 1^* \rbrace}}}{8\iu}\rho^2 .
\end{align}\end{subequations}

This result implies that, in the case of a single mode reduction, no numerical solution of the reduced dynamics is needed for the complex normal form style because, thanks to the symmetry of the formulation, the solution is available in explicit form. This property holds for generic order but does not extend to the case of multiple modes reduction.

We now move to the case of real normal form style. The reduced dynamics in the complex coordinates for a single mode case up to order three reads:
\begin{subequations}\label{eq:rnf_reddyn_singlemode}\begin{align}
&\dot{z}_1 = \lambda_1 z_1 + 
\p[3]{f_{1\lbrace 1 1 1^* \rbrace}} z_{1}^2\bar{z}_1 +
\p[3]{f_{1\lbrace 1 1^* 1^* \rbrace}} z_{1}\bar{z}_1^2 \, ,
\\
&\dot{\bar{z}}_1 = \bar{\lambda}_1 \bar{z}_1+ 
\p[3]{f_{1^*\lbrace 1 1 1^* \rbrace}} z_{1}^2\bar{z}_1 +
\p[3]{f_{1^*\lbrace 1 1^* 1^* \rbrace}} z_{1}\bar{z}_1^2\, .
\end{align}\end{subequations}

Two third-order monomials appear in both equations due to the fact that the resonant sets $\Is = \lbrace 1 1 1^* \rbrace$ and $\Is = \lbrace 1 1^* 1^* \rbrace$ are both considered resonant with $\Rs = \lbrace 1 1^* \rbrace$. An explicit expression of the solution is not available in this case due to the extra terms present in the equations. However, it can be shown that it is possible to write such system in a single oscillator form when Cartesian coordinates are used. Let us recall that, in the case of normal form style, the following property holds (see Eqs.~\eqref{eq:ff_mu}):
\begin{subequations}\begin{align}
\p{f_{r\Is}}+\p{f_{r^*\Is}} = - \phiv_r^\text{T}\M \p{\FU_\Is}.
\end{align}\end{subequations}

At order three, the right-hand side of the equation is equal to zero because all the lower order dynamic coefficients, only quadratic in this case, are zero; so the tensor $\p{\FU_\Is}$ for $\Is =\lbrace 1 1 1^* \rbrace$ and $\Is = \lbrace 1 1^* 1^* \rbrace$ vanishes. It follows that:
\begin{subequations}\label{eq:symm_RNF}\begin{align}
&\p[3]{f_{1\lbrace 1 1  1^* \rbrace}}+\p[3]{f_{1^*\lbrace 1 1  1^* \rbrace}} = 0,
\\
&\p[3]{f_{1\lbrace 1 1^*1^* \rbrace}}+\p[3]{f_{1^*\lbrace 1 1^*1^* \rbrace}} = 0.
\end{align}\end{subequations}

By summing Eqs.~\eqref{eq:rnf_reddyn_singlemode} and using this property, all the nonlinear monomials in the equation vanish, leading to:
\begin{equation}
\dot{z}_1 + \dot{\bar{z}}_1 
=
\lambda_1 z_1 +\bar{\lambda}_1 \bar{z}_1
=
\Re{\lambda_1}( z_1 +\bar{z}_1) + \iu \Im{\lambda_1} (z_1 - \bar{z}_1).
\end{equation}

Using Eqs.~\eqref{eq:cartes_def}, one can then write a linear equation for $\dot{\ar}_1$:
\begin{equation}\label{eq:dota1}
\dot{\ar}_1 = \Re{\lambda_1}\ar_1 - \Im{\lambda_1}\ar_{1^*} \, .
\end{equation}

As per the equation for $\dot{\ar}_1^*$, it can be derived from the difference between Eqs.~\eqref{eq:rnf_reddyn_singlemode} and it reads:
\begin{equation}\label{eq:dota1*}
\dot{\ar}_{1^*}
=
\Re{\lambda_1}\ar_{1^*} + \Im{\lambda_1} \ar_1
+ \sum_{{\Is}} \p[3]{\tilde{f}_{s{\Is}}} \; \p[3]{\tilde{\pi}_{{\Is}}},
\end{equation}
where the third order monomials have been collected in their generic form due to their lengthy expressions.

Finally, by differentiating Eq.~\eqref{eq:dota1} with respect to time and using Eq.~\eqref{eq:dota1*}, the oscillator-like reduced dynamics can be obtained as
\begin{equation}\label{eq:rnf_osclikedyn}
\ddot{\ar}_1 +
2\xi_1\omega_1 \,\dot{\ar}_1 +
\omega_1^2 \, \ar_1 
+ (\omega_1 \sqrt{1-\xi_1^2})\,\left(\sum_{{\Is}} \p[3]{\tilde{f}_{s{\Is}}} \; \p[3]{\tilde{\pi}_{{\Is}}}\right)= 0,
\end{equation}
where the expressions for real and imaginary parts of the eigenvalues have been used. It is worth remarking that this oscillator-like form of the equations can be obtained without approximations only if there is a linear differential relationships like Eq.~\eqref{eq:dota1} between the Cartesian coordinates. Otherwise, during the last substitution operation, terms of order higher than three would be generated. In such case, in order not to lose important coefficients, it is advisable to solve the equations for the Cartesian coordinates in first-order ODE form. 

It has been shown that the real normal form style reduced dynamics equations can be easily transformed into second-order oscillator-like equations in the case of a single mode reduction up to order three. In the general case of multiple modes reduction, this is still possible with no approximations, provided that there are no second-order internal resonances between the modes. On the contrary, for the graph style, the linear relationship between the Cartesian coordinates is always verified. In fact, the following property holds in general:
\begin{subequations}\label{eq:graph_symm}\begin{align}
\p{f_{r\Is}}+\p{f_{r^*\Is}} = 0,\qquad\forall r\in[1,2n], \forall\Is \, .
\end{align}\end{subequations}

In the case of graph style, this does not come from the fact that the reduced dynamics coefficients inside $\M \p{\FU_\Is}$ are zero, but rather from the fact that the mapping tensors inside $\M \p{\FU_\Is}$ are always orthogonal to the master modes, as shown in Eq.~\eqref{eq:graph_orthog}, thus their projection on $\phiv_r$ is zero. Moreover, in the case of the graph style the modal displacement $\ur_i$ coincides with the Cartesian coordinate $\ar_i$ for each $i\in[1,n]$; so not only an oscillator-like equation can always be obtained, but it would even coincide with the dynamics of the master modes. In the case of single mode reduction up to third order, the reduced dynamics would then read:
\begin{equation}\label{eq:graph_osclikedyn}
\ddot{\ur}_1 +
2\xi_1\omega_1 \,\dot{\ur}_1 +
\omega_1^2 \, \ur_1 
+ (\omega_1 \sqrt{1-\xi_1^2})\,\left(\sum_{{\Is}} \p[2]{\tilde{f}_{s{\Is}}} \; \p[2]{\tilde{\pi}_{{\Is}}}
+ \sum_{{\Is}} \p[3]{\tilde{f}_{s{\Is}}} \; \p[3]{\tilde{\pi}_{{\Is}}}
\right)= 0,
\end{equation}
where now the monomials ${\tilde{\pi}_{{\Is}}}$ are both second and third order and they are expressed in terms of $\ur_1$ in place of $\ar_1$.

In this Section, the reduced-order dynamics produced by the presented styles have been derived for a simple case of single mode reduction, up to third order. It is important to highlight that, although the reduced dynamics in the normal coordinates are different, once the original coordinates are reconstructed through the nonlinear mappings, their expression will be equivalent up to a certain order of expansion, as shown for example in~\cite{ROMGEOMNL}. This point will also be further discussed with numerical results in the next Section.

\section{Numerical Results}\label{sec:results}

In this section we discuss a series of applications of both academic and industrial interest that can be addressed with the  proposed method. The first example, reported in Section \ref{sec:arch}, concerns curved arches where the increase in curvature yields a transition of the response from softening to hardening. Interestingly, a softening behaviour is obtained at small amplitude, then a hardening behaviour is recovered at larger amplitudes. This example has been selected to show that such softening-hardening transition can be captured with a single-mode reduction and a development at least to order five.
The second structure, analysed in Section \ref{sec:cantilever},  is a cantilever beam. This example is challenging since inertia nonlinearities play an important role. Besides, preliminary results reported in~\cite{YichangVib} show that the second-order direct normal form (DNF) as implemented in~\cite{artDNF2020} was not able to capture the hardening behaviour up to very large amplitudes. Elaborating on this example, we show that the first mode manifold shows a folding point at large amplitudes. As a consequence, graph style parametrisation is not able to recover the correct behaviour, whereas normal form styles can. Finally we investigate an application of industrial interest, namely a MEMS micromirror subjected to large rotations, which highlights the potential impact of the presented method. Remarks on the computational performance of the approach are presented in Sec.~\ref{sec:performance}. All the examples are benchmarked against full-order harmonic balance finite elements (HBFEM) solutions computed using a custom fortran code~\cite{opreni2021piezo}. On the other hand, the reduced models are integrated using the harmonic balance libraries NLVib~\cite{krack2019harmonic} and ManLab~\cite{cochelin2009high}. Validations are also performed with the continuation package BifurcationKit~\cite{veltz2020bk}.

\subsection{Transition from softening to hardening behaviour: the case of a shallow arch}\label{sec:arch}

The development of ROMs is known to be easier for flat symmetric structures such as beams and plates since nonlinear quadratic couplings occur between bending and in-plane modes only and the slow/fast assumption is well fulfilled~\cite{Vizza3d,YichangVib}. On the other hand, arches and shells are known to integrate more couplings due to the loss of symmetry of the centre line and represent a challenge for reduction methods since the slow/fast assumption generally fails~\cite{TOUZE:JSV:2006,TOUZE:CMAME:2008,YichangNODYCON}. The aim of the present Section is to consider an arch with increasing curvature such that the small amplitude response is softening, and then turns back to hardening at larger amplitudes. Even though this kind of behaviour is typical of a single-mode model, it cannot be captured with third-order approximations as those used for example in~\cite{artDNF2020,AndreaROM}, since the change from softening to hardening requires at least the fifth-order term. The higher-order expansion is thus expected to solve this problem and proposes a single-mode ROM capable of reproducing such feature.

\begin{figure}[ht]
    \centering
    \includegraphics[width = .99\linewidth]{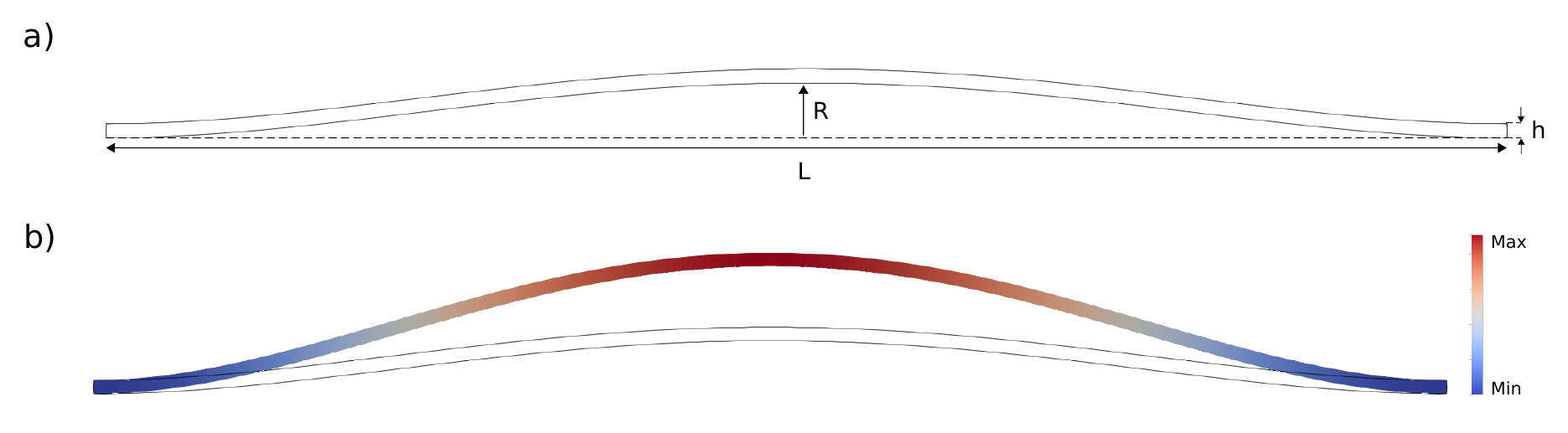}
    \caption{(a) Schematic representation of the shallow arch geometry, with thickness $h$ = 6.4 $\mu$m, and length $L$ = 640 $\mu$m. The rise of the arch $R$ for the three configurations is $R$\textsubscript{I} = 2.40 $\mu$m, $R$\textsubscript{II} = 2.88 $\mu$m, $R$\textsubscript{III} = 3.36 $\mu$m. The out-of-plane thickness of the arch is 32 $\mu$m for all geometries. (b) Magnified displacement field associated to the first eigenmode of the structure.}\label{fig:arch_geom_and_modes}
\end{figure}

A schematic representation of the layout is depicted in 
Fig.~\ref{fig:arch_geom_and_modes}(a). The systems are curved clamped-clamped beams of length 640 $\mu$m with an in-plane thickness of 6.4 $\mu$m and an out-of-plane thickness of 32 $\mu$m. The structures differ by their rise value $R$ and we denote as $\Ir$, $\Ir\Ir$ and $\Ir\Ir\Ir$
the arches with $R=2.40\,\mu$m, $R=2.88\,\mu$m and $R=3.36\,\mu$m, respectively. A solution obtained on a flat beam with a purely hardening response is presented as well for reference. The selected curvatures are small, thus corresponding to shallow arches, since the rise value $R$ is in the worst case close to half the thickness of the beam.

The reference flat beam structure (Ref) has the same geometry as the arches, but the rise is set to zero. All the structures are made of polycrystalline silicon, which is modelled as isotropic with a Young's modulus of 160 GPa, a Poisson's ratio of 0.22 and a density of 2320 kg/m$^3$.

Geometries are discretised using finite elements with 15-nodes (quadratic wedge elements). The total number of nodes is equal to 1161 for all geometries, corresponding to 3483 degrees of freedom.

\begin{table}[htb]
\centering
\caption{Shallow arches examples: rise of the different structures analysed and list 
of the first eigenfrequencies}
\begin{tabular}{|c|c|c|c|c|}
\hline
Geometry & Ref  & I & II & III \\ \hline
Rise [$\mu$m] & 0.00 & 2.4 & 2.88 & 3.36  \\ \hline
$\omega_1$ [rad/$\mu$s] & 0.8418 & 0.9223 & 0.9551 & 0.9923  \\ \hline
$\omega_2$ [rad/$\mu$s] & 2.3194 & 2.3191 & 2.3189 & 2.3188  \\ \hline
$\omega_3$ [rad/$\mu$s] & 4.1418 & 4.1359 & 4.1333 & 4.1303  \\ \hline
$\omega_4$ [rad/$\mu$s] & 4.5446 & 4.5607 & 4.5678 & 4.5763  \\ \hline
$\omega_5$ [rad/$\mu$s] & 7.5083 & 7.5067 & 7.5060 & 7.5052  \\ \hline
$\omega_6$ [rad/$\mu$s] & 10.252 & 10.267 & 10.274 & 10.282  \\ \hline
\end{tabular}
\label{tab:arc_eig}
\end{table}

The values of the first six eigenfrequencies of the structures are collected in Table~\ref{tab:arc_eig} to highlight the absence of low-order resonances between eigenmodes.
Single-mode reduction is used and the fundamental mode is selected as master. This corresponds, for the four cases, to the symmetric first bending mode.  A schematic representation of the displacement field associated to the mode is depicted in Fig.~\ref{fig:arch_geom_and_modes}(b) for illustration. 

The parametrisation order is spanned from 3 to 15 to analyse the convergence of the method. Both forced-damped and undamped solutions reported in the remainder of the section are validated using the HBFEM \cite{opreni2021piezo} applied to the full-order system and a Fourier expansion order equal to 7, which corresponds to a total of 52245 nodal unknowns.

\subsubsection{Undamped Response: backbone curves}

\begin{figure}[htb]
    \centering
    \includegraphics[width = .99\linewidth]{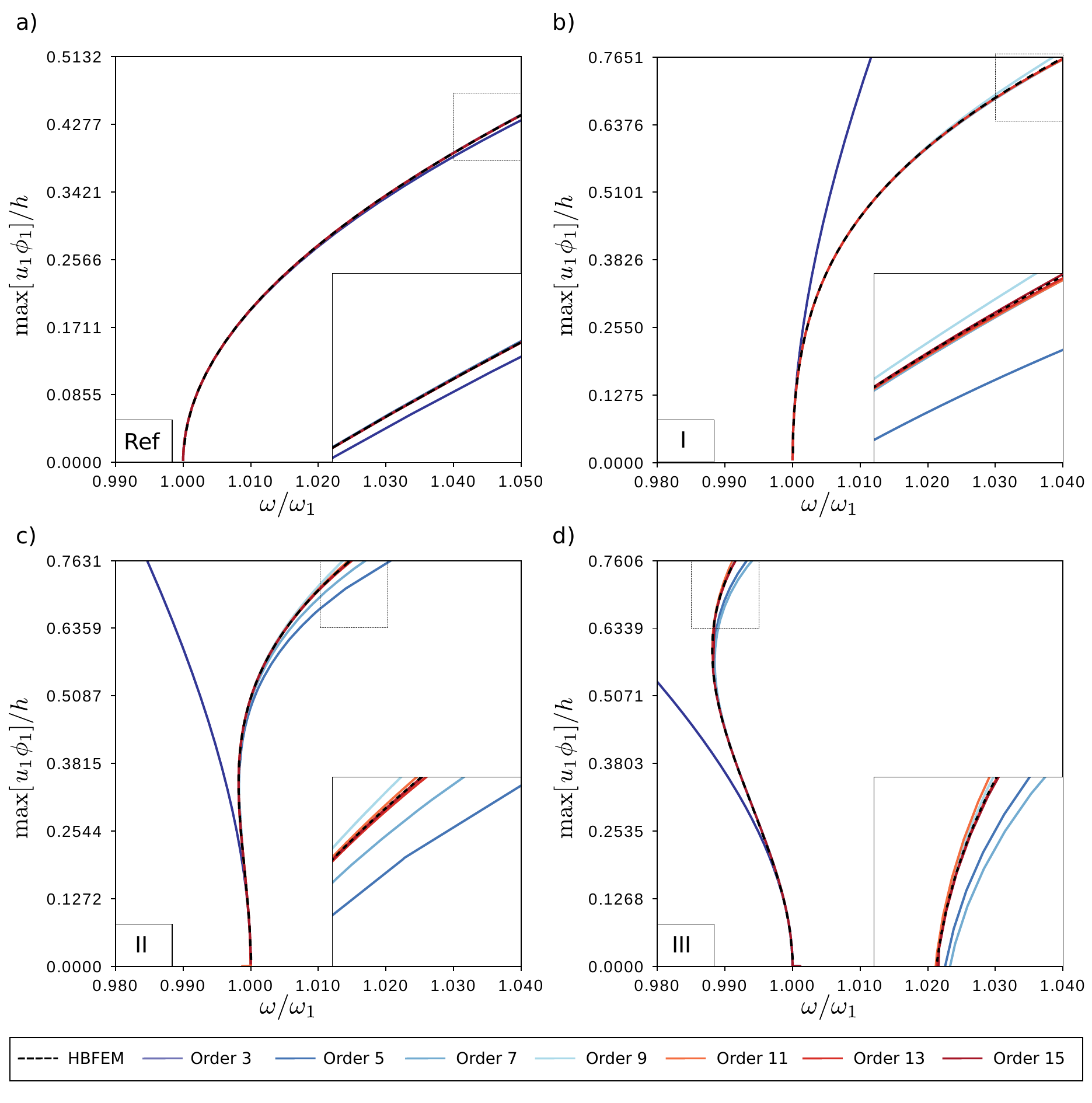}
    \caption{Backbone curves predicted for the arched structures, for the fundamental mode selected as master coordinate. Comparison between full-order HBFEM simulations and ROMs obtained from RNF style parametrisation. (a) Reference flat beam. (b) Arch I. (c) Arch II. (d) Arch III. The modal displacement $\ur_1$ is normalised by the maximum of the first eigenvector $\phiv_1$ and the thickness $h$; the nonlinear frequency $\omega$ is normalised by the linear frequency $\omega_1$ of the first mode.}\label{fig:arch_undamped_results}
\end{figure}

We first focus on the analysis of the backbone curves. The results obtained for the structures under consideration are reported in Fig.~\ref{fig:arch_undamped_results}, where a real normal form (RNF) is used as parametrisation style. 

The reference flat beam structure reported in Fig. \ref{fig:arch_undamped_results}(a) displays a hardening behaviour, typical of flat, symmetric structures. The third-order expansion already gives an excellent prediction in this case, and convergence up to half the beam's thickness is obtained with an order 5. On the other hand, third-order expansions fail for all the arch-structures under consideration. For instance, even for the purely hardening arch $\Ir$, whose backbone curve is reported in Fig. \ref{fig:arch_undamped_results}(b), the presence of large quadratic terms make low-order expansions deviate from the reference HBFEM solution. Indeed, third-order expansions start to become unacceptable for oscillation amplitudes equal to 0.2 of the normalised modal displacement, which corresponds to 20\% of the arches thickness, being the motion dominated by the master mode.
The effectiveness of low-order expansions is further reduced for arches with higher rise values. Indeed, arches $\Ir\Ir$ and $\Ir\Ir\Ir$ actually show an initially softening behaviour, followed by a hardening response at higher amplitudes, a behaviour that is totally missed by a third-order parametrisation, which diverges towards smaller frequency values. On the other hand, higher order parametrisations are able to correctly capture this transition and a very good agreement with the reference HBFEM solution is achieved with asymptotic expansions of order 13-15.

\subsubsection{Forced-Damped Response}

The forced-damped response of the system is analysed by selecting mass-proportional damping and modal forcing as proposed in \cite{AndreaROM}:
\begin{equation}
    \M \Vt + \C \V + \K \U +\G(\U,\U) + \H(\U,\U,\U)  = \kappa \M \phiv_{m} \cos{(\Omega t)},
\end{equation}
with $\C$ = $\omega_1/Q\M$. A complete treatment of the forcing in the invariant manifold procedure is challenging and generally leads to increasing the computational burden dramatically since the manifold reduction depends directly on the forcing frequency $\Omega$, as shown for example in~\cite{JIANG2005H,BreunungHaller18}. Instead, modal forcing as proposed in~\cite{touze03-NNM,TOUZE:JSV:2006} is a viable assumption that offers versatile ROMs for numerical continuation. The treatment of damping-forcing is also further commented in~\cite{artDNF2020}. \red{The treatment of forcing proposed in ~\cite{touze03-NNM,TOUZE:JSV:2006}, can be seen as a zero-order truncation of the perturbation procedure detailed in~\cite{BreunungHaller18} for the case of RNF and graph style parametrisation; it does not apply to the CNF style which would produce a time dependent manifold even at zero-order truncation due to the different treatment of resonant terms.} Here, only mass-proportional Rayleigh damping is considered. $Q=500$ is the quality factor of the system, $\kappa$ is a load multiplier and $\phiv_{m}$ is the master linear mode. The value of $\kappa$ for the reference flat beam is chosen equal to $0.06\omega_1^{2}$, while for the three curved structures 
$\kappa=0.09\omega_1^{2}$. 
The load multiplier for the reference flat beam is lower in order to avoid internal 
resonances that could be induced by the strongly hardening response of the structure. 

\begin{figure}[htb]
    \centering
    \includegraphics[width = .99\linewidth]{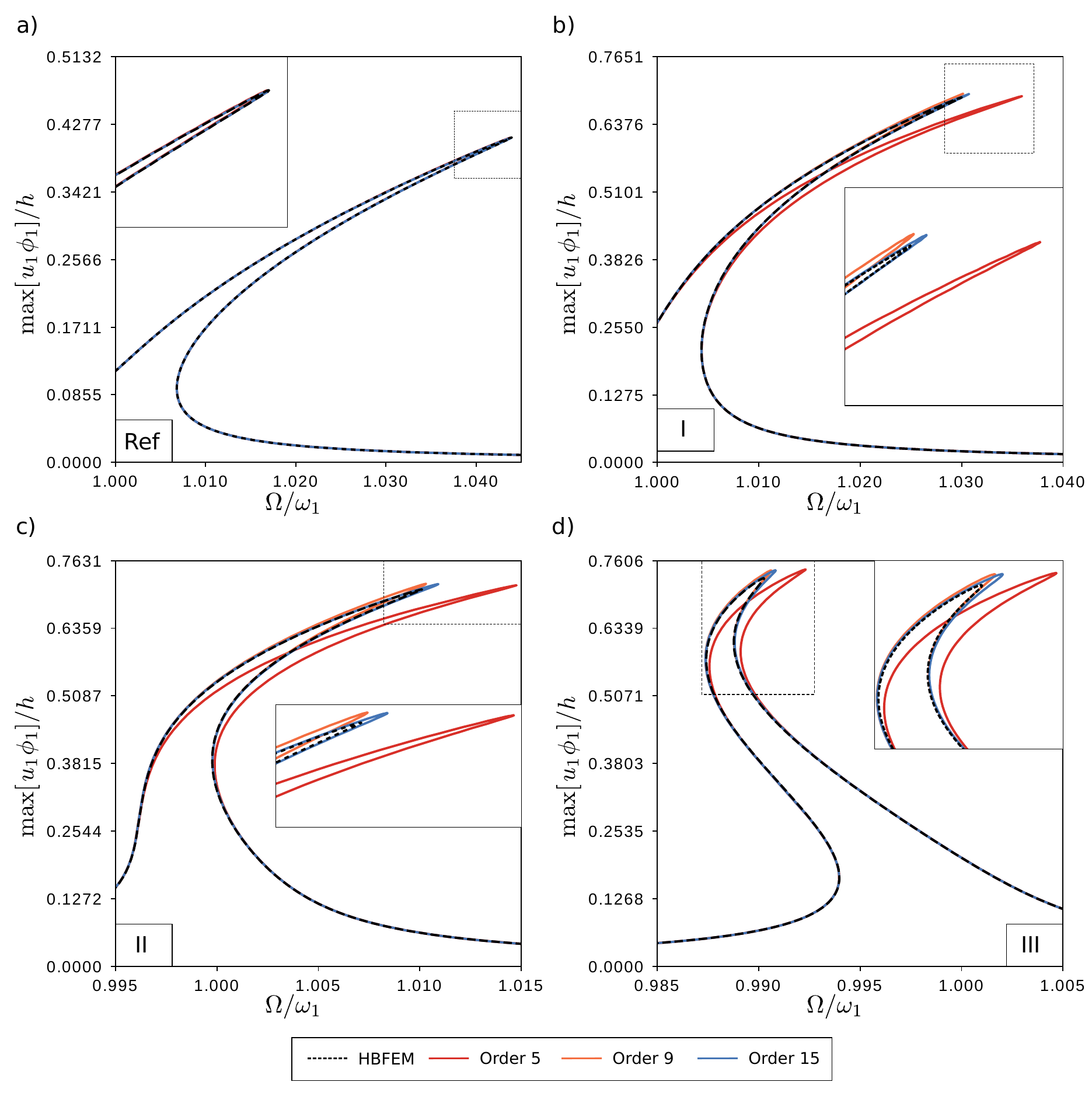}
    \caption{Frequency-response curves for the arched structures. Comparison between full-order HBFEM simulations and ROMs obtained from RNF style parametrisation. (a) Reference flat beam, (b)  Arch I. (c) Arch II. (d) Arch III.  The modal displacement $\ur_1$ is normalised by the maximum of the first eigenvector $\phiv_1$ and the thickness $h$; the forcing frequency $\Omega$ is normalised by the linear frequency $\omega_1$ of the first mode.}\label{fig:arch_damped_results}
\end{figure}

The reduced models are again obtained using a real normal form parametrisation on a single master mode. Fig.~\ref{fig:arch_damped_results} shows the result for the four cases under investigation. For the sake of simplicity, only three orders of asymptotic expansions are reported: orders 5, 9 and 15. As expected from the previous analysis on the backbone curves, the order-five solution is accurate enough only for the flat beam, but slightly departs from the reference solution as soon as a very small curvature is considered in case I. On the other hand, orders 9 and 15 gives excellent results up to the amplitudes  shown in Fig.~\ref{fig:arch_damped_results}, which corresponds to almost three fourth of the arch thickness.


The only minor discrepancies are observed at the peak of the frequency response curves and can be attributed to the treatment of the forcing term in the ROM. Indeed, if the parametrisation procedure was initially applied to the forced system, it would yield a reduced dynamics with coefficients that depend on the frequency and amplitude of the applied forcing, hence making standard continuation approaches not applicable and the final results less appealing in terms of versatility and post-processing for design purposes.
\red{It is also worth remarking that zero-order treatments of the forcing are reported also in \cite{JAIN2021How,MingwuLi2021_1} where the more general case of a non-modal forcing is included. However, as underlined for example in~\cite{JIANG2005H,artDNF2020}, even in presence of general forcing types, not accounting for the time dependence of the manifold does not seem to yield a sensible loss of accuracy, at least for moderate loads.}

As a short conclusion on this example, it has been demonstrated that the difficult case of an arch structure showing transition from softening to hardening behaviour can be finely predicted thanks to the arbitrary order expansions proposed in this article. As expected, the minimal order needed to retrieve the change in behaviour is 5, and to obtain convergence, higher-orders around 7 to 9 are needed. Shallow arch structures have been used in a number of contexts as a benchmark example to highlight the difficulties that ROMs can encounter for retrieving this behaviour, see {\em e.g.}~\cite{LacarboBuck98,nayfehcarboChin99,LACARBONARA2004nnm,VizzaMDNNM,MARCONI2021}. In general, the fact that the slow/fast assumption is not fulfilled is a strong obstacle for using single mode reduction. As an example, a very similar arch structure as the one studied here is detailed in~\cite{MARCONI2021} with a linear reduction method combining selected vibration modes plus their modal derivatives. Softening behaviour returning to hardening at higher amplitudes has also been successfully reported, with a reduction basis composed of 5 linear modes and 15 modal derivatives. Here we demonstrate that, with a nonlinear reduction method, a single-mode ROM is sufficient, but it needs an expansion order that is at least larger than five. For the sake of comparison, we want to highlight that, in the case of single mode reduction, at each order $p$, $p$+1 nonlinear mapping vectors have to be calculated. Therefore, the expansion up to order 5 of this example required the computation of 18 nonlinear mapping tensors plus the linear mode; however, the advantage of the method is that the mappings evaluation is an offline cost and the produced reduced-order model is a single degree of freedom oscillator.

\subsection{Manifold folding in presence of large transformations: the case of a cantilever beam}\label{sec:cantilever}

The freedom to choose the parametrisation style, as long as the homological equations are satisfied, introduces several options. As discussed in previous sections, if one adopts a graph style parametrisation, then modal displacement and velocity are in a one-to-one relation with the coordinates of the reduced dynamics. That is, a graph is built between the modal coordinates of the master and the slave. On the other hand, a normal form style of parametrisation implies that master coordinates are only identity tangent to the modal coordinates of the master mode. Theoretically speaking, the difference between these two parametrisation styles is assumed to be negligible until the manifold encounters a folding point. 
In the present section we report an example of a folding manifold, corresponding to the fundamental bending mode of a cantilever beam experiencing very large transformations.
We assume that the Saint-Venant Kirchhoff constitutive model holds for arbitrary large transformations, such that only geometric nonlinearities are at hand up to very large amplitudes.
Modelling cantilever beams has always been challenging for both model derivation and reduction methods, see {\em e.g.}~\cite{KimCantilever,CTENOC04,thomas16-ND,Meier2019,FAROKHI2020,YichangVib}. In particular, the motion of the structure becomes qualitatively different from that of the first bending mode and the displacement field changes in such a way that the modal coordinate of the first eigenmode is not a good measure of the oscillation amplitude. Indeed, when the vibration amplitude becomes comparable with the length of the structure, the modal amplitude of the bending mode tends to saturate, hence the modal amplitude does not increase anymore with the oscillation amplitude. It is worth stressing that this saturation is not associated to internal resonances between modes, thus leaving the invariant manifold two-dimensional. The aim of this Section is to clearly underline the failure of the graph style parametrisation in a case where the associated manifold presents a folding point.

\begin{figure}[ht]
    \centering
    \includegraphics[width = .99\linewidth]{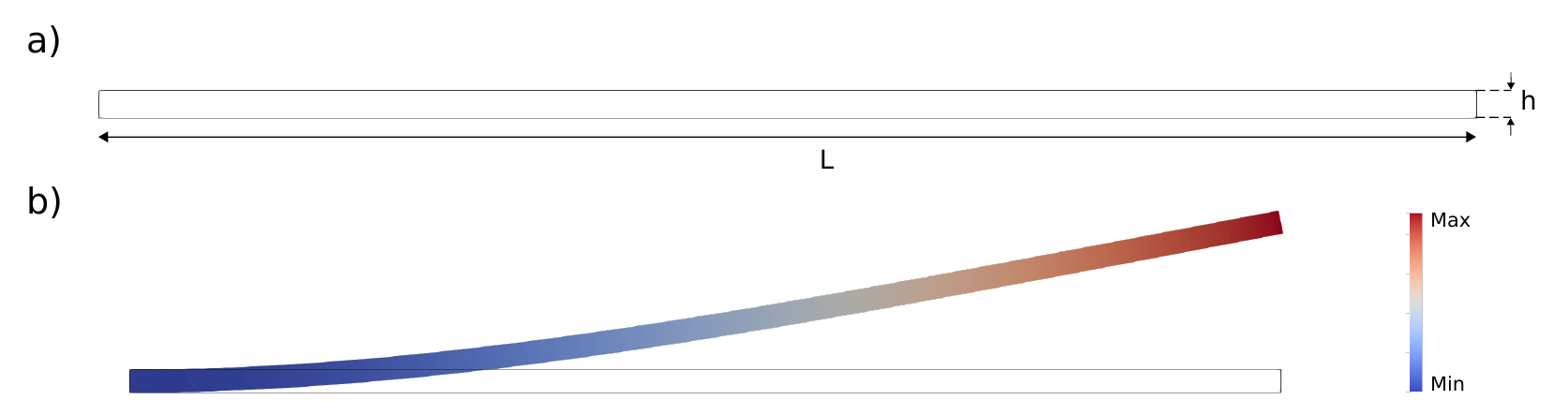}
    \caption{(a) Schematic representation of the cantilever geometry. $L$ = 1 m, $t$ = 0.02 m. The out-of-plane thickness is 0.05 m. (b) magnified displacement field of the eigenmode of interest.}\label{fig:cant_geom_and_modes}
\end{figure}

The geometry of the cantilever under consideration is reported in Fig.~\ref{fig:cant_geom_and_modes}(a). The structure is made of titanium with a Young modulus of 104 GPa, a Poisson ratio of 0.3, and a density of 4400 kg/m$^3$. The total length of the beam is 1~m and its thickness along the bending direction is 0.02~m. The out-of-plane thickness is 0.05~m. As for the arches, the geometry of the system is discretised using a FE procedure with 15 nodes wedge elements. The total number of nodes of the full order model is equal to 621, corresponding to 1863 degrees of freedom.

The first eigenmode of the structure corresponds to the typical bending mode represented in Fig.~\ref{fig:cant_geom_and_modes}(b). The eigenfrequency associated to the mode is equal to 99.18 rad/s.

\subsubsection{Folding of the fundamental mode's manifold}

The investigations are here focused on the conservative system in order to highlight the behaviour of the first mode up to very large amplitudes in terms of geometry of its associated manifold and corresponding backbone curves. The three parametrisation styles discussed in previous sections are used for reducing the system to its fundamental mode. The parametrisation order is spanned between 3 and 25 to verify the convergence of the methods.

\begin{figure}[ht]
    \centering
    \includegraphics[width = .99\linewidth]{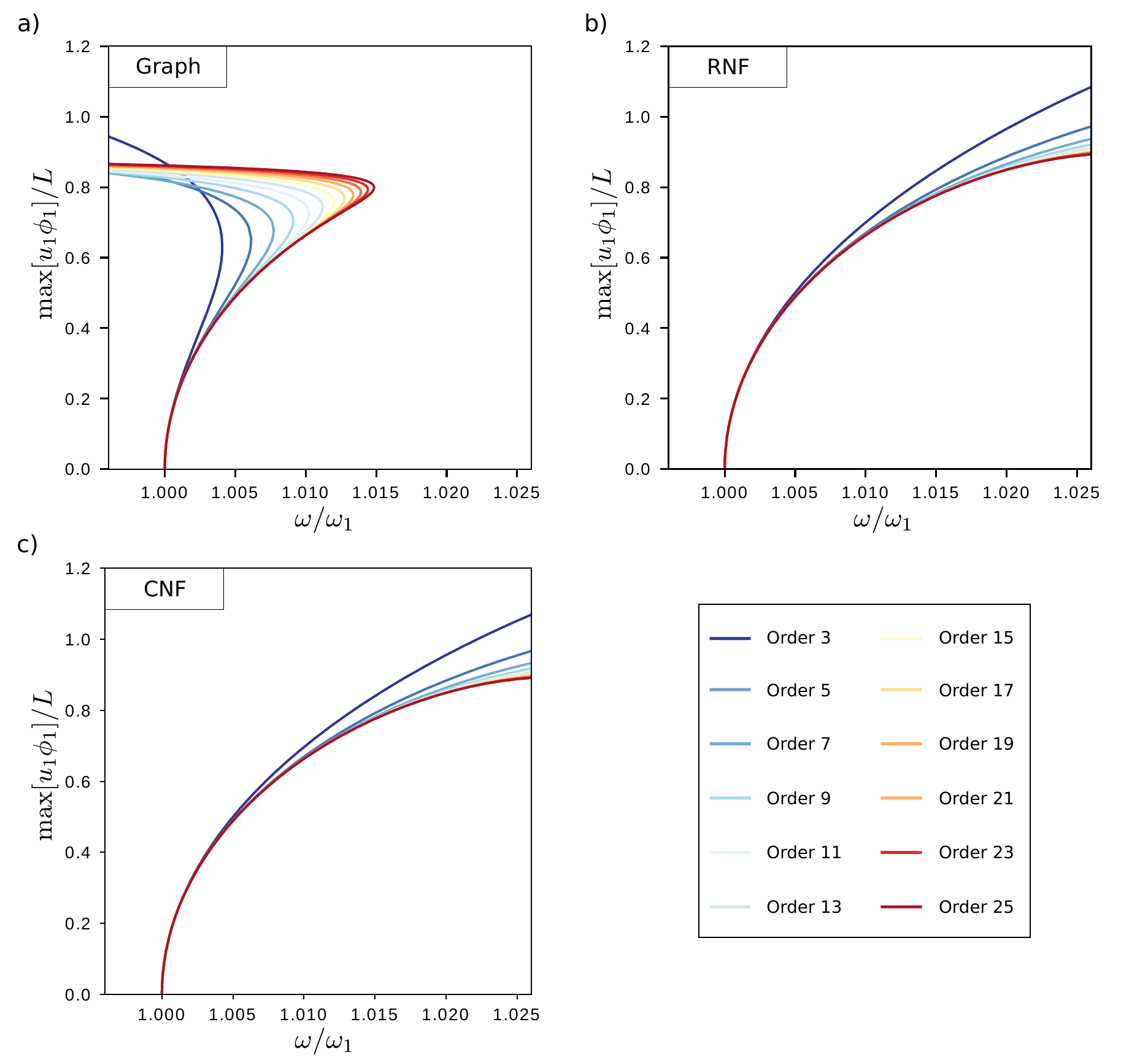}
    \caption{Backbone curves of the fundamental mode of the cantilever beam. Convergence of the asymptotic expansion upon a change in order for different parametrisations: (a) graph style, (b) RNF style, and (c) CNF style. Data are reported for expansions from order 3 to order 25. The modal displacement $\ur_1$ is normalised by the maximum of the first eigenvector $\phiv_1$ and the length $L$; the nonlinear frequency $\omega$ is normalised by the linear frequency $\omega_1$ of the first mode.}\label{fig:cant_undamped_1}
\end{figure}

The convergence results for the different parametrisation styles are reported in Fig.~\ref{fig:cant_undamped_1} from which it emerges that, while the results obtained for the two normal form styles are comparable within the considered frequency range, the graph style parametrisation shows an abrupt transition from hardening to softening behaviour. Focusing on the graph style's result shown in Fig.~\ref{fig:cant_undamped_1}(a), one clearly observes a convergence along the hardening behaviour with the expansion order, until a plateau is reached as a limit when increasing the orders up to 25. This abrupt behaviour appears to be a non-physical effect which needs further investigations and highlights the breaking of the graph style solution to a normalised vibration amplitude around 0.85, which corresponds to a physical displacement equal to 0.85 of the beam length $L$. On the other hand, normal form styles yield hardening curves within the frequency ranges investigated and convergence is achieved at low orders. 

It is worth mentioning that previous investigations reported in~\cite{YichangVib}, using a second-order DNF with a real normal form style not adopted in this study, showed the same behaviour as the graph style solution of order-three reported in Fig.~\ref{fig:cant_undamped_1}(a), with a correct hardening behaviour at the beginning followed by a strong departure to softening. This is explained by the fact that reduced-order dynamics of second-order DNF is equivalent to that of the graph style at third-order. Interestingly, the two normal form styles adopted in this study do not encounter this problem and show a correct behaviour from the very first orders. Further analytical investigations on the different styles and the failure of the first orders are reported in Section~\ref{subsec:multiplescales}. Using multiple scales method up to the second-order (in time scale, corresponding to reduced dynamics to order five) are analysed in order to highlight how the different solutions depart one from another, from the analytical expressions of the coefficients.

\begin{figure}[ht]
    \centering
    \includegraphics[width = .99\linewidth]{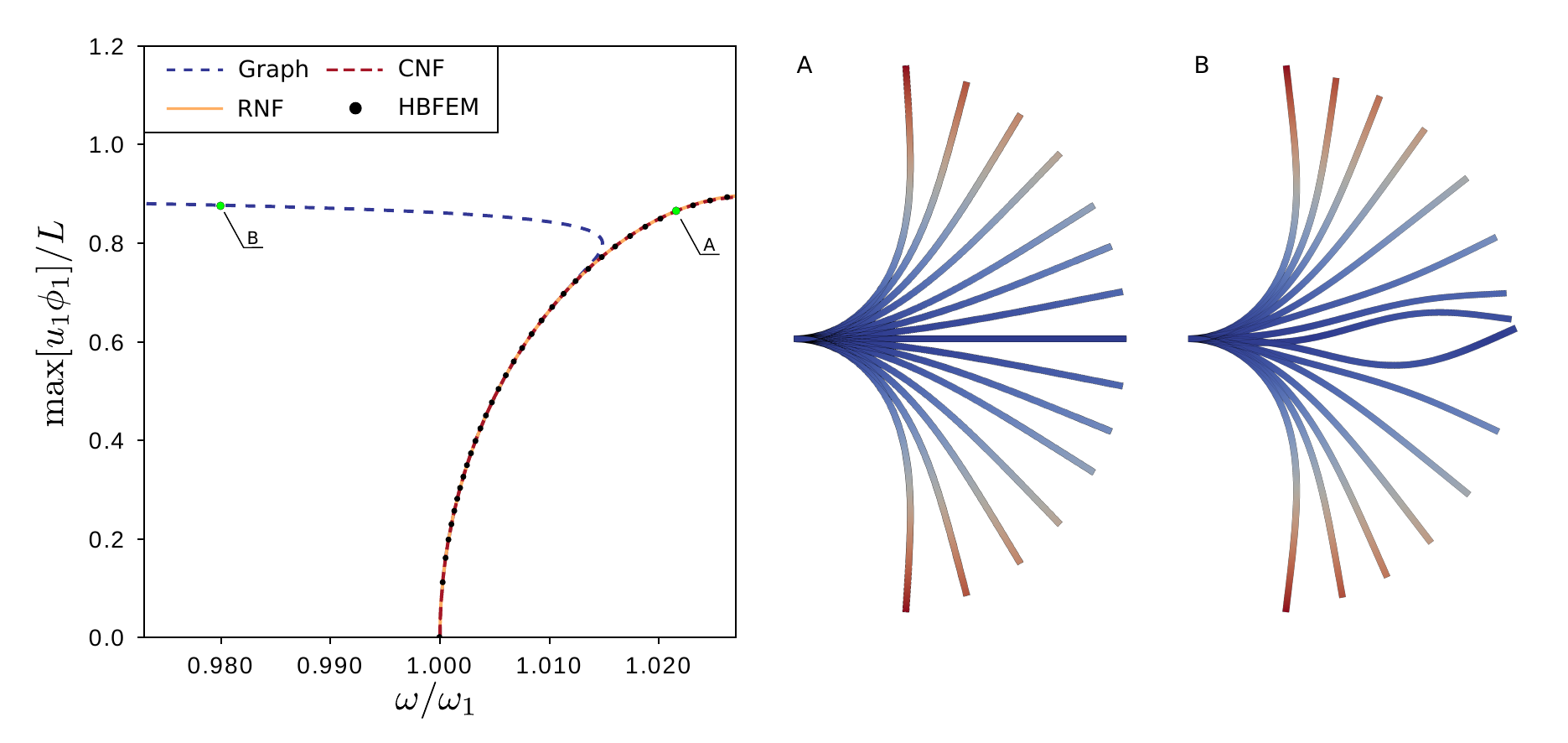}
    \caption{Backbone curve of fundamental mode of the cantilever beam. Comparison between full-order model solution (reference) computed with an HBFEM method (black points) and ROMs up to order 25 (converged solution) for graph (dotted blue), RNF and CNF (orange and dashed red) styles. The modal displacement $\ur_1$ is normalised by the maximum of the first eigenvector $\phiv_1$ and the length $L$; the nonlinear frequency $\omega$ is normalised by the linear frequency $\omega_1$ of the first mode. Panels A and B show snapshots of the deformation field of the structure at points A and B shown in the left (A: normal form style, B: graph style).}\label{fig:cant_undamped_2}
\end{figure}

\begin{figure}[p]
    \centering
    \includegraphics[width = .99\linewidth]{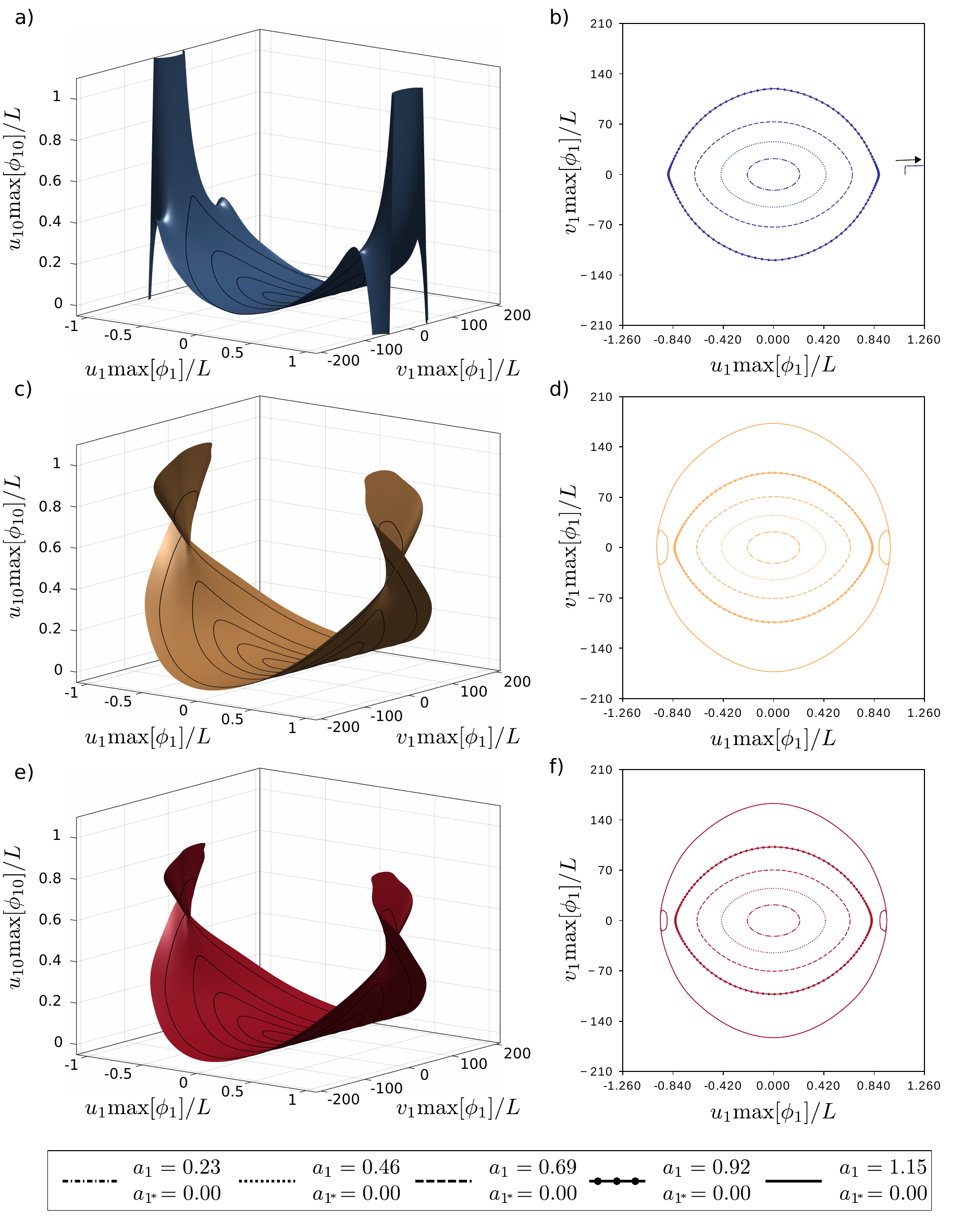}
    \caption{Cantilever beam. Left column: 3d representation of the invariant manifold in space $(\ur_1,\vr_1,\ur_{10})$, where $(\ur_1,\vr_1)$ are the modal displacement and velocity of the fundamental master mode while $\ur_{10}$ is the modal displacement of mode 10. The modal displacements $\ur_k$ are normalised by the maximum of their corresponding eigenvector $\phiv_k$ and the length $L$; the modal velocity $\vr_1$ is normalised by the maximum of the first eigenvector $\phiv_1$ and the length $L$. Right column: 2d projections of representative trajectories on the master mode plane. First row: graph style, second row: real normal form style, third row: complex normal form. The 5 trajectories depicted in pictures (b), (d), and (e) are also drawn on the manifolds as black lines, and are obtained for increasing normal amplitudes $\ar$. The arrow in figure (b) highlights that the trajectory initiated using as starting condition $\ar_1=1.15$ and $\ar_{1^*}=0$ diverges.}\label{fig:cant_manifold_comparison}
\end{figure}

The results obtained with the reduced-order models are next compared with full-order HBFEM simulations. The Fourier expansion order is set to 7 in order to provide sufficient accuracy while filtering potential internal resonances with high order modes which are typically observed in undamped systems~\cite{artDNF2020}. The comparison between order 25 expansions provided by the different parametrisation styles and the reference HBFEM solution is presented in Fig.~\ref{fig:cant_undamped_2}. The response predicted by the reference HBFEM solution is always hardening and it is perfectly reproduced by the reduced model parametrised with the normal form styles. Furthermore, the reliability of the ROM can also be inspected from the displacement field reconstructed from the reduced model. As shown from the physical reconstruction of the structure response at point A using a complex normal form style, the displacement field remains qualitatively similar to that of the first eigenmode. On the other hand, by reconstructing the displacement of the graph-style parametrisation at point B in Fig.~\ref{fig:cant_undamped_2}, it is possible to observe that nonphysical effects are present, especially close to the points of maximum normal velocity of the system. 

In order to understand the origin of the failure of the graph style parametrisation, Fig.~\ref{fig:cant_manifold_comparison} reports the shape of the invariant manifolds obtained with the different solutions. Three-dimensional representations are shown in the left column by using three modal coordinates, the first two being those of the master mode $(\ur_1,\vr_1)$, $\ur_1$ being the modal displacement and $\vr_1$ the associated modal velocity defined in Eqs.~\eqref{eq:modal_displ} and~\eqref{eq:modal_vel}. As third coordinate, the 10$^\mathrm{th}$ mode has been selected, $\ur_{10}$ being the associated modal displacement. Mode 10 corresponds to the first axial mode of this cantilever. It is strongly coupled to the first bending mode due to the classical axial-bending coupling, and one can see that the range of amplitudes reached along this coordinate $\ur_{10}$ is very important. Finally, axial-bending coupling being quadratic allows recovering a symmetric shape for the projection of the manifolds into this 3D representation. In the right column of Fig.~\ref{fig:cant_manifold_comparison}, a set of five trajectories computed for different initial conditions of the reduced dynamics coordinates are reported. They are computed for increasing values of the real-valued normal coordinate $\ar_1$ of the reduced dynamics. The initial condition is set with zero velocity corresponding to $\ar_{1^{*}}=0$. More specifically, these trajectories are computed with the ROM (single mode reduction), and then back-projected to the original coordinates thanks to the nonlinear mappings.  All manifolds are computed on a cantilever beam having the same geometry but with a more refined mesh of 2515 nodes in order to avoid any mesh effect on the shape of the computed manifold. 

The manifolds represented in Fig.~\ref{fig:cant_manifold_comparison} undoubtedly show that the failure of the graph style is due to a folding of the manifold in phase space. By its definition, the graph style parametrisation imposes an injective relation between the modal coordinates of the master mode and the slaves, and is thus unable to follow a folding in phase space.  This in turn causes a divergence of the manifold at folding points, as highlighted in Fig.~\ref{fig:cant_manifold_comparison}(a). The trajectories sampled with different initial conditions of the normal coordinates in Fig.~\ref{fig:cant_manifold_comparison}(b) show a similar result. In particular, trajectories initiated over the folding point diverge, which is in contrast with the nature of the stable fixed point of the system. On the other hand, the two normal form style parametrisations do not diverge in presence of a folding of the manifold as a benefit of the non-injective relation between modal and normal coordinates which allows capturing such processes. This is evidenced also by the trajectories reported in  Figs.~\ref{fig:cant_manifold_comparison}(d-e) showing intersections along the $(\ur_1,\vr_1)$ plane, which would not be possible using a graph-style parametrisation since planar systems develop trajectories that do not intersect. It is worth stressing that over the folding point the two normal form styles show a small difference in the computed solutions, that can be seen from the loops having a slight amplitude difference. This is in line with the fact that the two methods treats differently the resonant terms, however the deviations between the two methods are very small and can be appreciated only in very large amplitude, after the folding point, which is a critical case.


\subsubsection{Multiple Scales Expansion}\label{subsec:multiplescales}

\red{This Section is devoted to analyse the failure of the backbone obtained with graph style, to highlight the intimate relation between the geometry of the manifold and the dynamics on it; in fact, it will be shown that the lack of terms in the nonlinear mappings that caused the failure of the manifold in graph style comes with an excess of terms in the reduced dynamics that causes a failure of the backbone.} To this purpose, the backbone of the different ROMs are analysed with a multiple scales solution. The aim is to show how the departure of the different methods can be analysed in terms of the first developments of the backbone curve (amplitude-frequency relationship), by showing how the coefficients of the different styles intervene in the solution. Since the backbone curves of Fig.~\ref{fig:cant_undamped_1} show qualitatively similar results with the order of expansion for all the styles, we will focus our attention here on the reduced dynamics equations provided by each style up to the third-order. 

Recalling Eqs.~\eqref{eq:cnf_reddyn_singlemode}, the reduced dynamics in complex form for the CNF style reads:
\begin{subequations}\begin{align}
&\dot{z}_1 = \iu \omega_1z_1 + f z_1^2z_{1^*},
\\
&\dot{z}_{1^*} = -\iu \omega_1z_{1^*} - f z_1z_{1^*}^2,
\end{align}\end{subequations}
where we denoted as $f = \p[3]{f_{1\lbrace 111^*\rbrace}}$ the only nonlinear coefficient appearing in the dynamics, since, for a conservative system, one has $\p[3]{f_{1^*\lbrace 11^*1^*\rbrace}} = -\p[3]{f_{1\lbrace 111^*\rbrace}}$. Note that $f$ is purely imaginary in the considered case of a conservative system.  Recalling Eqs.~\eqref{eq:rho_alpha}, which rewrites the reduced dynamics in polar form,  the decay ratio is thus vanishing: $\dot{\rho} =0$, and the dynamical equation on the phase makes appear the nonlinear oscillation frequency $\omega$ such that $\dot{\alpha} = \omega$. Finally, the  amplitude-frequency relationship simply reads:
\begin{equation}\label{eq:bbCNFexpl}
\omega = 
\omega_1 (1-\dfrac{\iu f}{4 \omega_1}\rho^2) \, ,
\end{equation}
which is real since  $f$ is purely imaginary.

In the case of the RNF style, the reduced dynamics has two additional terms, due to the treatment of resonances in real formulation, as shown in Eqs.~\eqref{eq:rnf_reddyn_singlemode}. It reads:
\begin{subequations}\begin{align}
&\dot{z}_1 = \iu \omega_1z_1 + f z_1^2z_{1^*} + f z_1z_{1^*}^2 \, ,
\\
&\dot{z}_{1^*} = -\iu \omega_1z_{1^*} - f z_1z_{1^*}^2 - f z_1^2z_{1^*} \, .
\end{align}\end{subequations}

Indeed, in light of Eqs.~\eqref{eq:symm_RNF}, the additional reduced dynamics coefficients $\p[3]{f_{1\lbrace 11^*1^*\rbrace}}$ and $\p[3]{f_{1^*\lbrace 111^*\rbrace}}$ are equal to $-\p[3]{f_{1^*\lbrace 11^*1^*\rbrace}}$ and $-\p[3]{f_{1\lbrace 111^*\rbrace}}$, respectively. As shown in Eq.~\eqref{eq:rnf_osclikedyn}, the second-order oscillator-like dynamics in the Cartesian normal coordinates simply reads:
\begin{equation}
\ddot{a}_1 +
\omega_1^2 \, a_1 +
\omega_1^2 a_1
\left(-\dfrac{\iu f}{2 \omega_1}a_1^2 -\dfrac{\iu f}{2 \omega_1} \left(\frac{\dot{a}_1}{\omega_1}\right)^2\right) = 0.
\end{equation}

An analytical expression of the amplitude-frequency  relationship of such a system can be derived by means of the multiple scales method, which will allow one to compare the coefficients with those of the CNF. Before showing the multiple scales solution of this system, we give the expressions of the graph style dynamics. The reduced dynamics in complex coordinates reads: 
\begin{subequations}\begin{align}
&\dot{z}_1 = \iu \omega_1z_1 + f z_1^2z_{1^*} + f z_1z_{1^*}^2 + \hat{f} z_1^3 + \hat{f}  z_{1^*}^3 \, ,
\\
&\dot{z}_{1^*} = -\iu \omega_1z_{1^*} - f z_1z_{1^*}^2 - f z_1^2z_{1^*} - \hat{f}  z_1^3 - \hat{f}  z_{1^*}^3\,.
\end{align}\end{subequations}
One can remark that, in this case, no quadratic terms are present, simply because the original system representing the dynamics of the flexural mode of a flat cantilever does not display such self-quadratic terms for symmetry reason, see {\em e.g.}~\cite{givois2019,Vizza3d}. Concerning the cubic terms,
 four additional monomials appear. A new coefficient $\hat{f}$ has been added to replace the values of $\p[3]{f_{1\lbrace 1111\rbrace}}$ and $\p[3]{f_{1\lbrace 11^*1^*1^*\rbrace}}$ which can be shown to be equal. The symmetry between the first and second equation stems instead from Eq.~\eqref{eq:graph_symm}. 

The second-order oscillator-like equation in the graph style case reads:
\begin{equation}
\ddot{a}_1 +
\omega_1^2 \, a_1 +
\omega_1^2 a_1
\left(
-\dfrac{\iu (f+\hat{f})}{2 \omega_1}
a_1^2 
-\dfrac{\iu (f-3\hat{f})}{2 \omega_1} 
\left(\frac{\dot{a}_1}{\omega_1}\right)^2
\right) = 0.
\end{equation}

In order to give an explicit expression for the backbone curve for both RNF and graph styles, a multiple scales expansion can be performed on third-order oscillator-like equations. Let us consider a generic form for the equation fitting to both styles: 
\begin{equation}
\ddot{a}_1 +
\omega_1^2 \, a_1 +
\omega_1^2 a_1
\left(
c_{30}
a_1^2 
+c_{12}
\left(\frac{\dot{a}_1}{\omega_1}\right)^2
\right) =0,
\end{equation}
where the coefficient $c_{30}$ represent the coefficient multiplying $a_1^3$ and the coefficient $c_{12}$ that multiplying $a_1 \dot{a}_1^2$. The amplitude-frequency relationship in this case would read:
\begin{equation}
\omega = \omega_1 ( 1 + \p[2]{\Gamma} \rho^2 + \p[4]{\Gamma} \rho^4),
\end{equation}
with:
\begin{subequations}\begin{align}
&\p[2]{\Gamma} = \dfrac{3c_{30}+c_{12}}{8},
\\
&\p[4]{\Gamma} = \dfrac{-15 c_{30}^2 + 14 c_{30} c_{12} + c_{12}^2}{256}.
\end{align}\end{subequations}

Replacing the coefficients $c_{30}$ and $c_{12}$ of the RNF and graph styles by their values, both methods yield a $\p[2]{\Gamma}$ value that matches that of the CNF and reads:
\begin{equation}
\p[2]{\Gamma} = -\dfrac{\iu f}{4 \omega_1}
\end{equation}

However, the value of $\p[4]{\Gamma}$ is different between the two. In the RNF case it is simply zero because $c_{30}=c_{12}$, thus matching the backbone predicted by the CNF which, by definition of developments up to third-order, has no fourth-order frequency dependence.
In the case of the graph style, it reads:
\begin{equation}
\p[4]{\Gamma_\text{Graph}} = \dfrac{(4 f + 3 \hat{f})\hat{f}}{64 \omega_1^2}.
\end{equation} 

The numerical value for $f$ and $\hat{f}$ for the present cantilever case are $6.81662 \iu$ and $-66.9832 \iu$, respectively. At low amplitudes, the frequency amplitude curve is dominated by the $\p[2]{\Gamma}$ contribution of $\rho^2$, so the three methods are in perfect agreement. The numerical value of $\p[2]{\Gamma}$ for the three styles is positive, as shown by the hardening effect at low amplitudes. Conversely, the numerical value of $\p[4]{\Gamma_\text{Graph}}$ is negative, therefore when the amplitude reaches a certain level, the backbone curve obtained with the graph style will make a transition from hardening to softening when the $\rho^4$ term becomes more important than the $\rho^2$ one. This effect is however nonphysical, as shown in the comparison with the full model results of Fig.~\ref{fig:cant_undamped_2}. 

In this perspective, not only the  simpler form of the mappings obtained with graph style is unable to replicate the occurrence of a folding of the manifold, but also the additional reduced dynamics coefficients of the graph style cause an nonphysical softening effect in the backbone curve.

\begin{figure}[t]
    \centering
    \includegraphics[width = .55\linewidth]{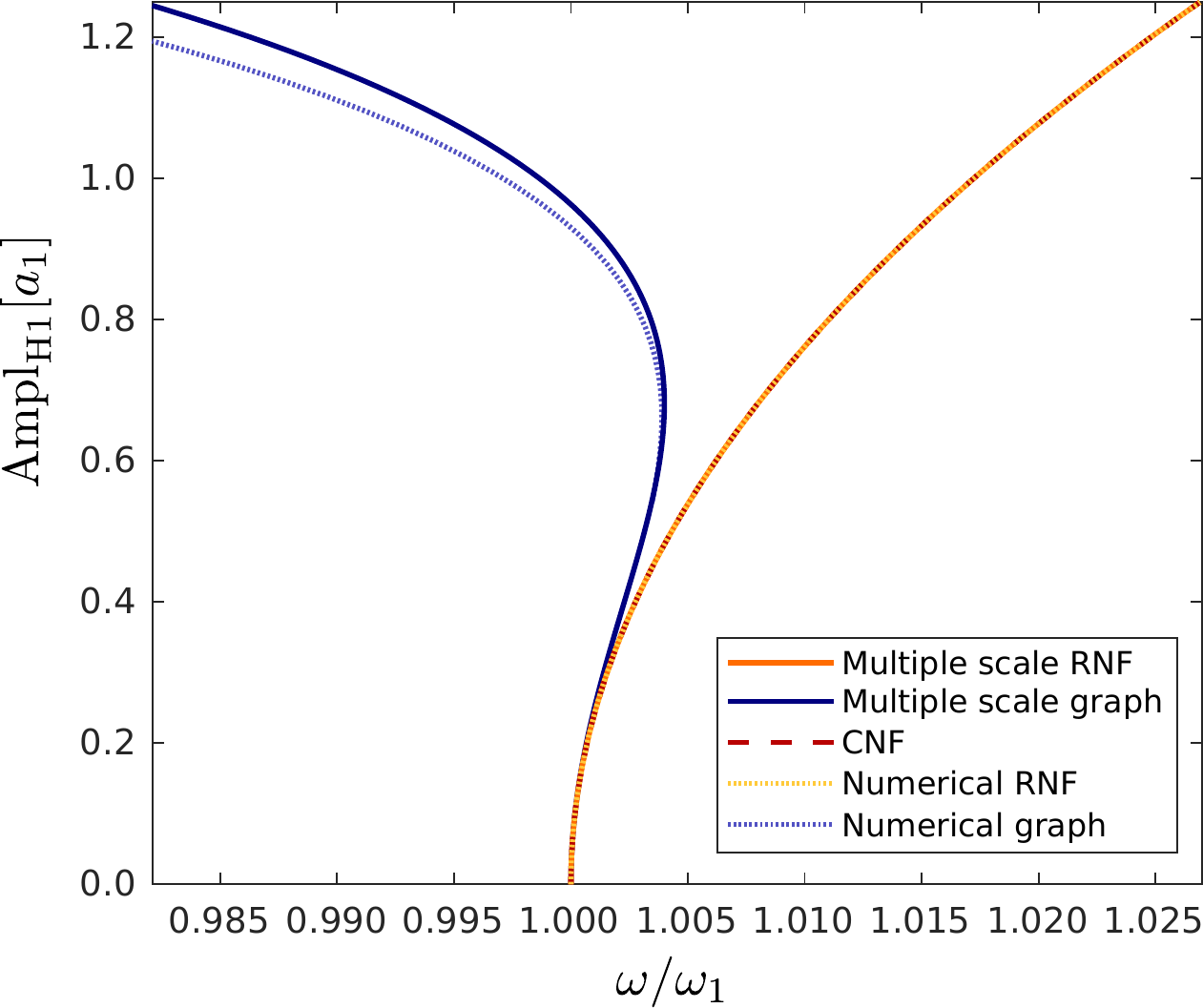}
    \caption{Comparison of numerical and analytical backbone curves for the reduced dynamics up to order three for the three styles under study (CNF, RNF and graph). Amplitude of the first harmonic of $a_1$, which coincides with $\rho$ for the analytical solutions, plotted against the nondimensional nonlinear frequency. For graph style and normal form styles numerical solutions are obtained with HBM including 30 harmonics and multiple scales development are obtained up to $\p[4]{\Gamma}$. For complex normal form style numerical and analytical backbones coincide.}\label{fig:MS_cant}
\end{figure}

Fig.~\ref{fig:MS_cant} illustrates the analytical developments by comparing the backbone provided by the multiple scales expansions to those numerically obtained with the HBM on the reduced dynamics up to order three. It can be observed that a very good agreement is found, showing that multiple scales solutions indeed closely follows the numerical one up to large amplitudes, stressing that the discussion on the different values of $\p[4]{\Gamma}$ given by each method is meaningful. It must be underlined that the expansion order of the perturbation method (up to power four in terms of amplitude: $\p[4]{\Gamma} \rho^4$) is larger than the expansion order of the reduced dynamics (limited at order three). This has been deliberately selected in order to account for using high-order numerical solutions (as HBM with numerous harmonics) on reduced dynamics having lower order, thus replicating the effect of accurate numerical methods such as HBM or collocation. The present development shows that, the sudden change of behaviour of the graph style backbone when approaching the folding point can be explained by the incorrect sign of the the next-order term of the expansion.

\subsection{MEMS Micromirror}\label{sec:mirror}

\begin{figure}[ht]
    \centering
    \includegraphics[width = .85\linewidth]{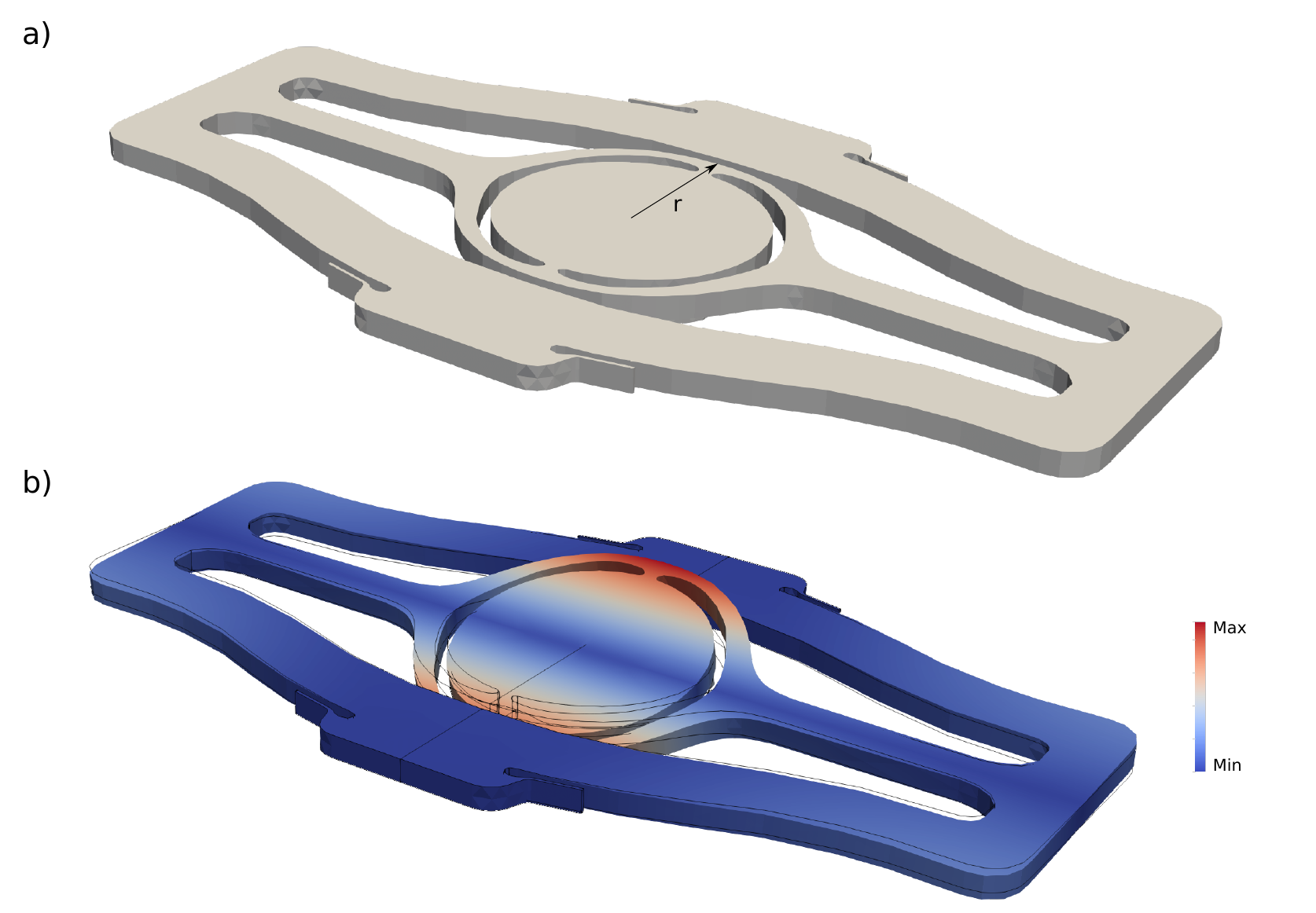}
    \caption{(a) MEMS micromirror geometry. (b) Magnified displacement field associated to the third eigenmode of the structure.}\label{fig:perseus_geom_and_modes}
\end{figure}

In the present section we apply the reduction method to a case of remarkable industrial interest. MEMS micromirrors are core components in many high-end industrial applications and their performance requirements are steadily increasing. Therefore, accurate estimation of the structure nonlinear response is of paramount importance during the design stage of the device.\\
An example of MEMS micromirror is reported in Fig.~\ref{fig:perseus_geom_and_modes}(a). The device is developed by STMicroelectronics\textregistered. It is made of monocrystalline silicon, whose mechanical properties are detailed in~\cite{hopcroft2010young}. The real structure is actuated using eight lead-zirconate-titanate (PZT) patches and an extensive study of their effect on the dynamic response of the system is detailed in~\cite{opreni2021piezo}. Since this section aims at showing the model-order reduction technique only, the piezoelectric force is replaced by a modal  loading following the procedure used for the arch-structures in Section~\ref{sec:arch}. To reproduce the operating conditions of the real device, a quality factor $Q$ of 1000 is selected, together with a load multiplier $\kappa$ equal to 1, 1.5, 2.00, 2.50, and 3.00 $\mu$m/$\mu$s$^2$. The highest $\kappa$ value allows reaching a rotation of approximately 12$^{\circ}$, which corresponds to the maximum operating range of the device.\\

The geometry of the structure, illustrated in Fig. \ref{fig:perseus_geom_and_modes}(a), is discretised as before using quadratic elements with 15 nodes. The number of nodes of the final geometry is 3244, which corresponds to 9732 degrees of freedom. The structure is parametrised along its third eigenmode, corresponding to a torsion of the central reflective surface. The associated resonance frequency is 0.1839 rad/$\mu$s$^2$, which is larger than the one reported in~\cite{opreni2021piezo} since a coarser mesh is used within the present work. The displacement field corresponding to the eigenmode is illustrated in Fig.~\ref{fig:perseus_geom_and_modes}(b).

\begin{figure}[ht]
    \centering
    \includegraphics[width = .99\linewidth]{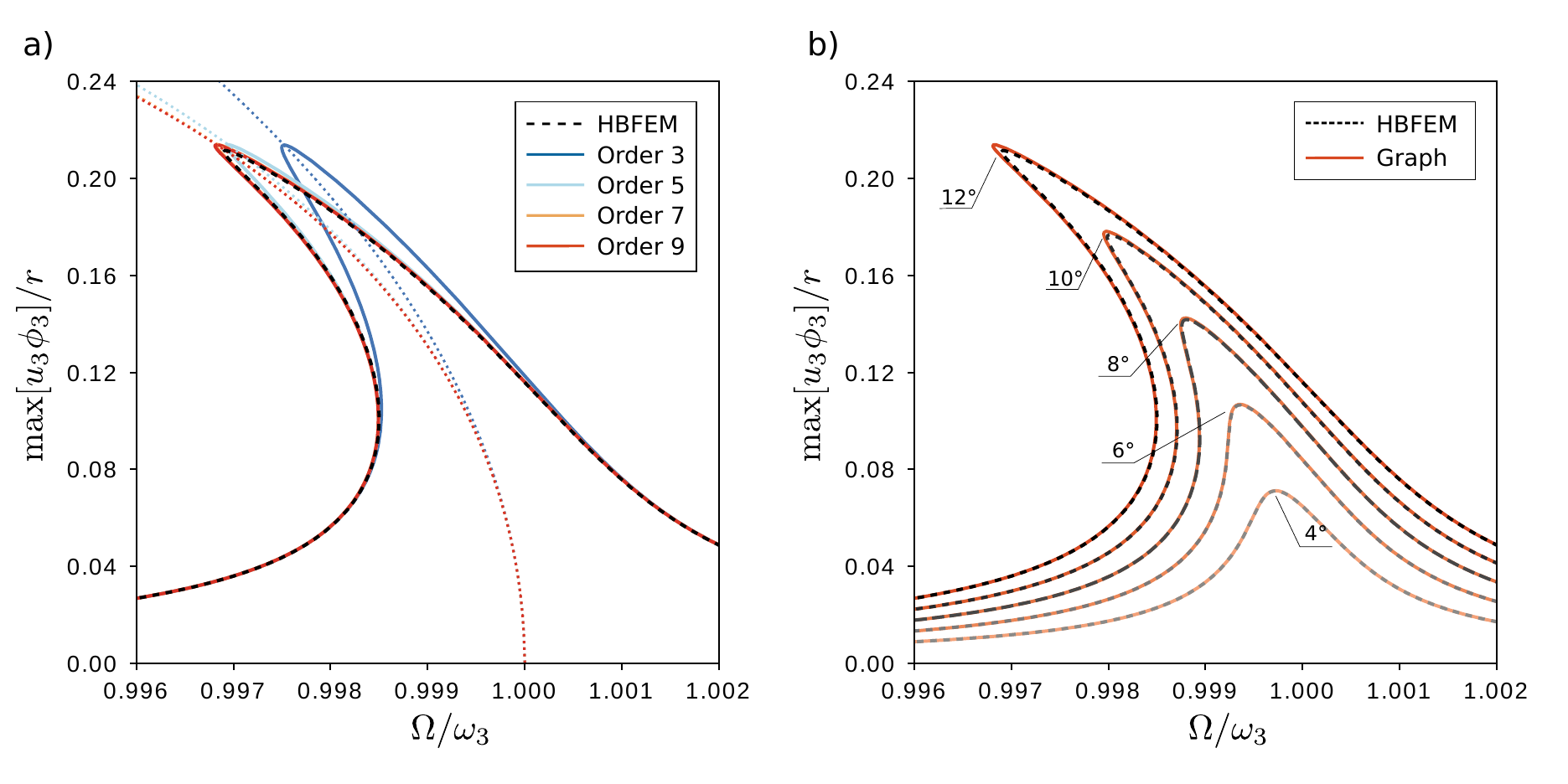}
    \caption{Frequency-response functions for the MEMs micromirror. The modal displacement $\ur_3$ is normalised by the maximum of the third eigenvector $\phiv_3$ and the radius $r$; the nonlinear frequency $\omega$ is normalised by the linear frequency $\omega_3$ of the third mode. (a) Convergence of a graph style parametrisation by varying the expansion order for $\kappa=3\mu$m/$\mu$s$^2$ and comparison with reference HBFEM. (b) Comparison between HBFEM solution and order 9 expansion for load multiplier $\kappa$ values equal to 1.0, 1.5, 2.0, 2.5, and 3.0 $\mu$m/$\mu$s$^2$.  Tags report the physical rotation angle reached by the device at the peak of the FRF.}\label{fig:perseus_res}
\end{figure}

Reduction is computed with the graph style and the convergence of the result is measured by spanning the parametrisation order from 3 to 9. The resulting frequency-response functions (FRFs) for $\kappa$ equal to 3.00 $\mu$m/$\mu$s$^2$ are plotted in Fig.~\ref{fig:perseus_res}(a) and compared to the full-order reference solution obtained by direct HBFEM method. At this level of amplitudes, the third-order expansion clearly deviates from the reference solution, and  a satisfactory convergence is reached from order 7, hence highlighting how high order expansions are essential to cover a wide range of industrial applications.

The data obtained from the ROM with an order 9 graph style parametrisation and the HBFEM solution for all $\kappa$ values are collected in Fig.~\ref{fig:perseus_res}(b). The comparison highlights the excellent performance of the proposed reduction procedure for any forcing level, which proves its ability to model structures of arbitrary geometrical complexity. For this case, the different styles of parametrisations gave similar results, justifying the choice of showing only the graph solution. Interestingly, a very slight departure from the full-order solution can be observed at the maximal amplitude of the FRF, for the largest values of the forcing. Since the only remaining approximation in the method is on the treatment of the forcing, this slight departure is attributed to the assumption of a time-varying manifold which is not dependent on the normal coordinate. Nevertheless, the results of the computations reported clearly shows that this effect is small and might be neglected with good reason.

The computational gain obtained in this example is remarkable. The time required to derive the reduced model is equal to 2 minutes and 39 seconds, and less than one minute was necessary to compute all the frequency response functions. On the other hand, over three days were required to obtain the full order HBFEM solutions. The next Section investigates further this aspect by detailing the computational burden involved by scaling to large models composed of millions of dofs.

\subsection{Remarks on computational performance}\label{sec:performance}

The cost of model-order reduction comprises the offline phase of the reduction procedure and the computing time required to solve the reduced model itself.  
In the present framework, the latter is negligible, since the reduced model contains only a single master mode, regardless of the size of the original system. Therefore, in this section we consider the parametrisation of the MEMS micromirror as a benchmark to investigate the offline computational performance of the method. For reference, all the analyses were performed on a desktop workstation with an AMD\textregistered Ryzen 5950X processor and 128 GB RAM.

First, we report the computing time required to obtain the reduced model for a fixed mesh size and for different parametrisation orders, to highlight the time required to achieve high-order expansions on structures with a moderate number of nodes. Secondly, we show how the computing performance changes for a given expansion order, by varying the mesh refinement. All the analyses are performed on the same model detailed in Section~\ref{sec:mirror}.

\begin{figure}[ht]
    \centering
    \includegraphics[width = .99\linewidth]{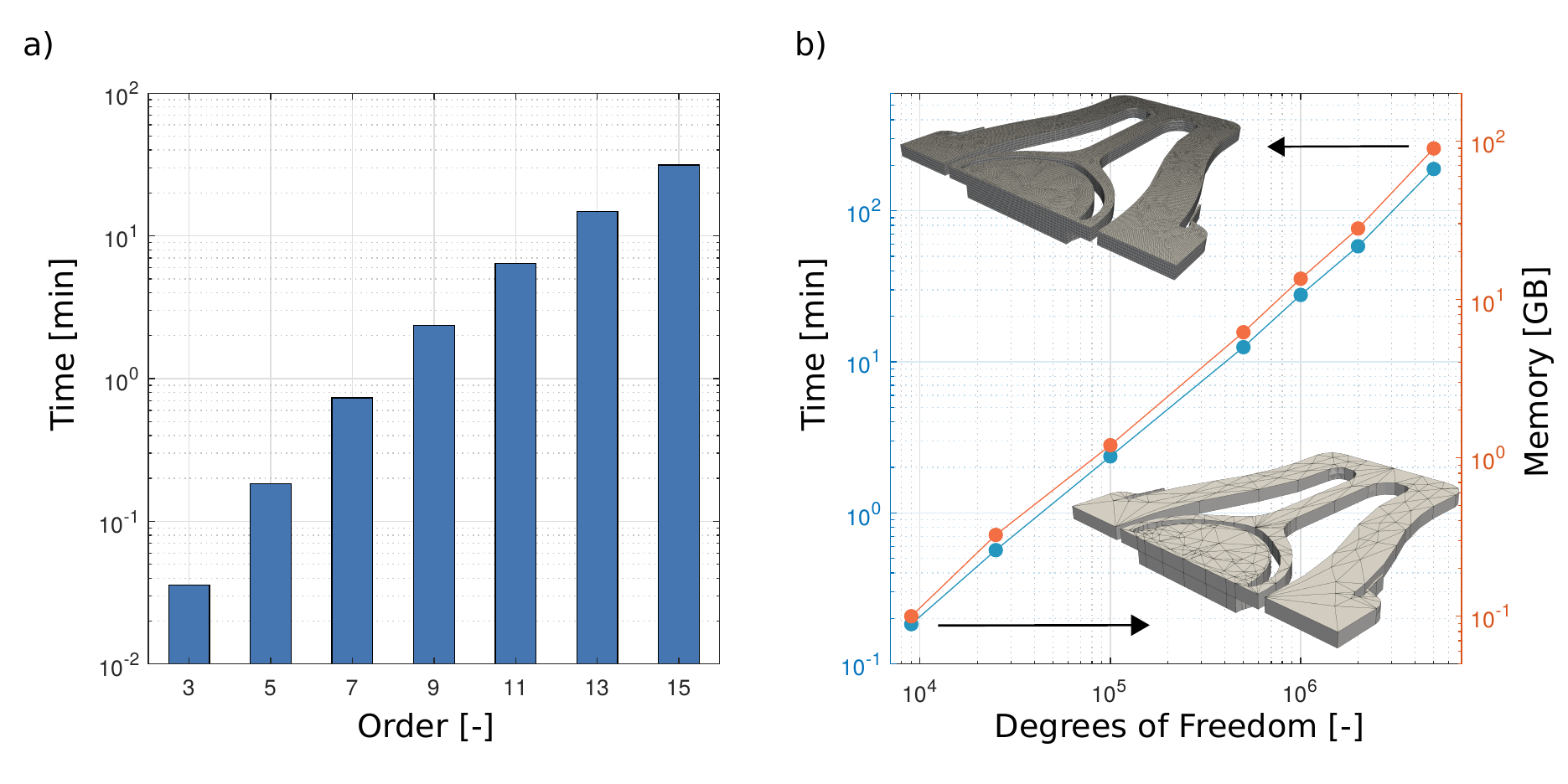}
    \caption{(a) Computing times required to reach a given parametrisation order, for a mesh of 9732 degrees of freedom. (b) Computing times required to perform an order 5 parametrisation with RNF style, for meshes of increasing refinement. All the analyses were performed on a desktop workstation with an AMD\textregistered Ryzen 5950X processor and 128GB RAM.}\label{fig:performance}
\end{figure}

The computing times obtained by varying the expansion order for a mesh of 9732 degrees of freedom are collected in Fig.~\ref{fig:performance}(a). The number of linear systems to be solved at each order scales as $p^{2n-1}$, with $p$ the expansion order, and $n$ the number of master modes, which in the present contribution is equal to one. 
For a single master mode reduction, the number of linear systems to be solved at each order $p$ scales linearly with the order, so the total number of linear systems that need to be solved for an expansion up to order $p$, scales quadratically with $p$.
However, the number of operations required to compute the right hand side of each homological equation does not scale quadratically, hence the trend is slightly less than quadratic. The peak memory consumption required by the parametrisation procedure for this example was equal to 600 MB.

The computing times and memory requirements to obtain an order 5 parametrisation for different mesh refinements are reported in Fig. \ref{fig:performance}(b). As highlighted from the charts on a log-log scale, the memory use is almost linear and the same is observed for the analysis time. We stress that, in order to obtain an order-5 ROM for a system having 5 millions degrees of freedom, the proposed method requires approximately 3 hours.

\section{Conclusion}

In this contribution, an arbitrary order expansion allowing to compute accurate ROMs relying on invariant manifold theory has been proposed. In essence, the method provides a nonlinear mapping that relates directly the physical space (nodes of the FE model) to an invariant-based span of the phase space, and the resulting reduced dynamics on the manifold. Both expansions of the mapping and the reduced dynamics can be computed up to a generic order of expansion, hence reaching convergence in numerous test cases. In the solving of the parametrisation method, three different styles of solutions have been developed, and numerous computational details have been highlighted in order to decrease the computational burden.

\red{A special emphasis has been put on initial problems provided by a structural FE discretisation, with a damping matrix diagonalised by the linear modes of the conservative problem. Since damping is taken into account in the developments, the invariant manifold computations
presented in this contribution are thus truncated expansions of the unique SSM, such submanifold being reached only if the order of the expansion is larger than the spectral quotient, a number that is very large in practical applications.}

The method has been applied to three cases with different complexities. The first case is an arch with increasing curvature, and it has been demonstrated that backbone curves, that are initially softening and then turning back to hardening behaviour at larger amplitudes, can be reproduced accurately with a single-mode reduction and an expansion order of at least 5. A detailed study of the fundamental bending mode of a cantilever beam has been investigated, and it has been found that the associated invariant manifold encounters a folding point at large amplitudes. As a consequence, the graph style parametrisation is not able to reproduce the behaviour. On the other hand, the two normal form styles are able to pass this folding point and accurately follow the full-order solution up to extreme amplitudes. Finally, the case of a MEMS micromirror of industrial relevance has been selected. In the operating range of the device, third-order expansions were clearly not sufficient, such that the previously proposed methods based on DNF in~\cite{artDNF2020,AndreaROM} are not accurate enough. Reaching an order seven in the expansion shows a perfect convergence.

As a matter of fact, only one assumption is still remaining in the derivation of the proposed ROM, now that higher-orders have been included: the treatment of the forcing. \red{Taking the forcing directly into account from the beginning of the derivation, as shown in \cite{BreunungHaller18,JAIN2021How} for instance, would lead to a time-varying manifold. The calculations would then require an additional parametrisation to account for said variations, having as a main consequence that the ROM  needs to be recomputed for each value of the external excitation frequency, as shown for example in~\cite{JIANG2005H,MingwuLi2021_1}.} This is a severe drawback in that the ROM loses a number of appealing features as compared to those used in the present study where continuation solutions are easily accessible to compute FRFs. The remedy is to make the classical assumptions of small forcing as the one already used in~\cite{touze03-NNM,TOUZE:JSV:2006,Haller2016,VERASZTO,PONSIOEN2018}. The computations on the test cases considered in this paper shows that the assumption is responsible for a very slight departure of the ROM to the reference solution, and only at very large amplitude values.

As already emphasised for example in~\cite{ROMGEOMNL}, the parametrisation method and the invariant manifold theory offers a sound theoretical background to derive accurate ROMs for nonlinear structures. Invariance property is key and enforcing its fulfilment from the beginning is a guarantee to produce ROMs with powerful predictive capacities. Strong results from dynamical systems theory ensures that the long-time behaviour of the solutions lies in the vicinity of these manifolds, consequently approximating them is the best solution to produce effective ROMs. As underlined for example in~\cite{TOUZE:JFS:2007}, linear methods like POD (Proper Orthogonal Decomposition) are necessarily approximating the invariant manifolds by adding new basis vectors to capture its curvature, such that spatial features retrieved by proper orthogonal modes are here directly computed from the model and can be traced back in the nonlinear mappings.

\section*{Funding}
The work received no additional funding.

\section*{Conflict of interest}
The authors declare that they have no conflict of interest.
\section*{Data availability statement}
The data that support the findings of this study are available from the corresponding author, upon reasonable request.
\bibliographystyle{unsrt}
\bibliography{biblioROM}

\appendix
\section{Complex eigenproblem properties}
\label{app:eigs}
The aim of this section is to detail the derivation of the complex eigenproblem properties given in the text.

We can write th complex right eigenvector in a more compact form as:
\begin{equation}
\v{Y}_s = \begin{bmatrix}
\phiv_s\lambda_s\\ \phiv_s
\end{bmatrix}
\end{equation}
with $s\in [1,2N]$. As for the complex left eigenvector:
\begin{equation}
\v{X}_r = \frac{1}{\lambda_r - \bar{\lambda}_r}
\begin{bmatrix}
\phiv_r \\
- \phiv_r\bar{\lambda}_r
\end{bmatrix}
\end{equation}
with $r\in [1,2N]$.

The derivation of the first orthogonality property of Eq.~\eqref{eq:massnorm_complex} reads:
\begin{equation}
\v{X}_r^\text{T}
\begin{bmatrix}
\M & \0
\\
\0 & \M
\end{bmatrix}
\v{Y}_s = 
\frac{\phiv_r^\text{T} \M \phiv_s\lambda_s - \bar{\lambda}_r\phiv_r^\text{T} \M \phiv_s}{\lambda_r - \bar{\lambda}_r}
=
\frac{\lambda_r- \bar{\lambda}_r}{\lambda_r - \bar{\lambda}_r}\delta_{sr}
=\delta_{sr}
\end{equation}
where the mass orthogonality property of real eigenvectors $\phiv_r^\text{T} \M \phiv_s = \delta_{sr}$ has been used.

The derivation of the second orthogonality property of Eq.~\eqref{eq:stiffnorm_complex} reads:
\begin{equation}
\v{X}_r^\text{T}
\begin{bmatrix}
\C & \K
\\
-\M & \0
\end{bmatrix}
\v{Y}_s = 
\frac{ \phiv_r^\text{T} 
( \C \phiv_s \lambda_s+ \K \phiv_s )
- \bar{\lambda}_r\phiv_r^\text{T} 
(-\M \phiv_s \lambda_s)}{\lambda_r - \bar{\lambda}_r}
=
\frac{2\xi_r\omega_r\lambda_r + \omega_r^2 + \bar{\lambda}_r\lambda_r}{\lambda_r - \bar{\lambda}_r}
\delta_{sr}\label{eq:stiff_orth_demo}
\end{equation}
where the stiffness orthogonality property of real eigenvectors $\phiv_r^\text{T} \C \phiv_s = 2\xi_r\omega_r \delta_{sr}$  and the damping orthogonality property of real eigenvectors $\phiv_r^\text{T} \K \phiv_s = \omega_r^2 \delta_{sr}$ have been used. Recalling Eqs.~\eqref{eq:def_complex_eigvals}, one has that: 
\begin{equation}
\bar{\lambda}_r\lambda_r = 
\left(-\xi_r \omega_r + \iu \omega_r \sqrt{1-\xi_r^2}\right)
\left(-\xi_r \omega_r - \iu \omega_r \sqrt{1-\xi_r^2}\right)
=
\xi_r^2\omega_r^2  + \omega_r^2 (1-\xi_r^2)
=
\omega_r^2
\end{equation}
and that:
\begin{equation}
-(\bar{\lambda}_r+\lambda_r) = 
2\xi_r \omega_r\;.
\end{equation}
By plugging the last two into Eq.~\eqref{eq:stiff_orth_demo}, writes:
\begin{equation}
\v{X}_r^\text{T}
\begin{bmatrix}
\C & \K
\\
-\M & \0
\end{bmatrix}
\v{Y}_s = 
\frac{-(\bar{\lambda}_r+\lambda_r)\lambda_r +2\bar{\lambda}_r\lambda_r}{\lambda_r - \bar{\lambda}_r}
\delta_{sr}=
\frac{-(\lambda_r-\bar{\lambda}_r)\lambda_r}{\lambda_r - \bar{\lambda}_r}
\delta_{sr}
=
-\lambda_r\delta_{sr}
\end{equation}
which coincides with Eq.~\eqref{eq:stiffnorm_complex}.

The derivation of the right eigenproblem reads:
\begin{equation}
\left(
\lambda_s
\begin{bmatrix}
\M & \0
\\
\0 & \M
\end{bmatrix}
+
\begin{bmatrix}
\C & \K
\\
-\M & \0
\end{bmatrix}
\right)
\v{Y}_s
=
\begin{bmatrix}
(\M\lambda_s^2 + \C\lambda_s + \K)\phiv_s
\\ 
\M\phiv_s\lambda_s - \M\phiv_s\lambda_s
\end{bmatrix}
=\0
\end{equation}
where the top rows correspond to the second order eigenproblem of Eq.~\eqref{eq:lin_eig_secondorder} and the bottom ones are identical.

The derivation of the left eigenproblem reads:
\begin{equation}
\v{X}_r^\text{T}
\left(
\lambda_r
\begin{bmatrix}
\M & \0
\\
\0 & \M
\end{bmatrix}
+
\begin{bmatrix}
\C & \K
\\
-\M & \0
\end{bmatrix}
\right)
=
\frac{1}{\lambda_r - \bar{\lambda}_r}
\begin{bmatrix}
\phiv_r^\text{T} ( \M (\lambda_r+\bar{\lambda}_r) + \C)
&
\phiv_r^\text{T} (\K - \bar{\lambda}_r \lambda_r\M )
\end{bmatrix}
\end{equation}
that is verified by noticing that the matrices are symmetric and that:
\begin{subequations}\begin{align}
( \M (\lambda_r+\bar{\lambda}_r) + \C)\phiv_r = 
( - 2\xi_r \omega_r \M + \C ) \phiv_r = \0
\\
(\K - \bar{\lambda}_r \lambda_r\M )\phiv_r = 
(\K - \omega_r^2 \M) \phiv_r = \0
\end{align}\end{subequations}

Finally, we demonstrate Eq.~\eqref{eq:symm_of_syst} in the present case of damping matrix $\C$ diagonalised by the eigenvectors of the conservative system. One has to demonstrate that:
\begin{equation}
((\lambda_r) \M + \C)\phiv_r = (- \bar{\lambda}_r)\M \phiv_r
\end{equation}
If both sides are premultiplied by any $\phiv_{s\neq r}^\text{T}$ they are both zero so they are equal. If both sides are premultiplied by $\phiv_{r}^\text{T}$ they read:
\begin{equation}
 \lambda_r + 2\zeta_r \omega_r = - \bar{\lambda}_r
\end{equation}
which leads to:
\begin{equation}
\lambda_r + \bar{\lambda}_r = - 2\zeta_r \omega_r 
\end{equation}
which is true by definition of the complex eigenvalues. Being the matrix of all eigenvectors an isomorphism, this is a general result.

\section{Real normal form styles}
\label{app:oldrealNF}

This section is devoted to the presentation of the real normal form style derived in~\cite{touze03-NNM,TOUZE:JSV:2006,artDNF2020}, which is different from the real normal form style introduced in this contribution. In order to explain clearly the difference between these two real styles, let us call RNF the real normal style used in the present contribution, and FRNF for full real normal form, the one that has been used in~\cite{touze03-NNM,TOUZE:JSV:2006,artDNF2020}. The FRNF mainly differs from RNF by the fact that more monomials are considered as resonant, which can be simply understood by interpreting the resonance condition differently. Indeed, in FRNF, the resonance condition used can be rewritten as:
\begin{equation}\label{eq:old_nf_cond}
\left(\sum_{i_k\in\Is}\pm\lambda_{i_k}\right)^2 \approx \left(\pm\lambda_r\right)^2
\end{equation}
for any value of $\pm$. The fact that any value of the sign is included in the condition, makes this style different from the RNF style presented in this article; for instance, the set $\Is = \lbrace 1,1,1 \rbrace$ is not considered resonant with $\lambda_1$ nor $\lambda_{1^*}$ in the RNF, but it is in the condition of Eq.~\eqref{eq:old_nf_cond}. \red{This choice have two main effects}: all the calculations can be realised in a complete real formalism, without the need of any complex at any stage of the development (see {\em e.g.}~\cite{touze03-NNM,TOUZE:JSV:2006}); on the other hand, many symmetries are lost \red{and more terms are present in the reduced dynamics.}

More \red{precisely}, the RNF requires a complex parametrisation in $\v{z}$, whereas the real coordinates used in FRNF can be seen as the Cartesian normal coordinates $\v{a}$. In fact, each monomial in $\v{a}$ is composed of several monomials in $\v{z}$; for instance, the monomial $a_1^3$ is a linear combination of $z_1^3$, $z_1^2 z_{1^*}$, $z_1 z_{1^*}^2$, and $z_{1^*}^3$. If we denote by ${\tilde{\WU}}$ the coefficients of the nonlinear mappings written with Cartesian coordinates to differ them from those of the mappings written with complex coordinates ${\WU}$, one has that  the mapping $\p[3]{\tilde{\WU}_{111}}$ that multiplies the monomial $a_1^3$ in the reconstruction, is a linear combination of $\p[3]{{\WU}_{111}}$, $\p[3]{{\WU}_{111^*}}$, $\p[3]{{\WU}_{11^*1^*}}$, and $\p[3]{{\WU}_{1^*1^*1^*}}$. The same reasoning holds for the reduced dynamics. A fully real development is then only possible if the treatment of the sets $\lbrace 111^* \rbrace$, $\lbrace 11^*1^* \rbrace$, $\lbrace 111 \rbrace$, and $\lbrace 1^*1^*1^* \rbrace$, is identical, which is what the condition given by Eq.~\eqref{eq:old_nf_cond} imposes. It is possible to show the same would apply in the presence of multiple modes and internal resonances.

To better understand the differences between the three normal form styles presented (CNF, RNF and FRNF), let us present the results of the different parametrisation styles on a simple Duffing oscillator, with an asymptotic expansion up to third order. The starting point is the equation of motion written as:
\begin{equation}
\ddot{u} + \omega_0^2 u + \gamma u^3 = 0.
\end{equation}
To make the developments and explanations as simple as possible, no quadratic terms have been considered such that no second-order terms will be present neither in the mapping, nor in the reduced dynamics.

The mappings for each of the styles, written with complex normal coordinates, generically writes, up to the third-order, as:
\begin{equation}
    u = z+\bar{z} + 
    \sum_{i_1 = 1,1^*}\sum_{i_2 = 1,1^*}\sum_{i_3 = 1,1^*}
    \Psi^{(3)}_{i_1i_2i_3} \;z_{i_1}z_{i_2}z_{i_3},
\end{equation}
following the general formula of the expansions introduced in the main text. In this simple case, $u$ is now a scalar such that $\Psi^{(3)}_{i_1i_2i_3}$ is also scalar. Using the fact that in this simple case one has  $z_1 = z$ and $z_{1^*} = \bar{z}$, and thanks to  the property of Eq.~\eqref{eq:symm_complex} together with the fact that all coefficients $\Psi$ are purely real in the case of conservative systems, the nonlinear mapping can be simply rewritten for the case of the Duffing equation as
\begin{equation}\label{eq:mapdufcomplexco}
    u = z+\bar{z} + 
    \Psi_0 \; (z^3+\bar{z}^3) + \Psi_2\; (z^2\bar{z} + z\bar{z}^2 ),
\end{equation}
\red{where only two coefficients, simply rewritten as $\Psi_0$ and $\Psi_2$, are needed to incorporate all the $\Psi^{(3)}_{i_1i_2i_3}$}.

In order to fully compare all the different possible representations, let us also introduce the same mapping but with Cartesian coordinates, which reads:
\begin{equation}
    u = a + 
    \sum_{i_1 = 1,1^*}\sum_{i_2 = 1,1^*}\sum_{i_3 = 1,1^*}
    \tilde{\Psi}^{(3)}_{i_1i_2i_3} \;a_{i_1}a_{i_2}a_{i_3},
\end{equation}
where the coefficients are\red{, as well,} noted as scalars and with a tilde in order to distinguish them from the complex formulation. If one sets $a_1 = a$ and $a_{1^*} = b$, the mapping in this case can be simply rewritten as:
\begin{equation}
    u = a + 
    {\tilde{\Psi}_0}\; a^3 + {\tilde{\Psi}_2}\; a b^2,
\end{equation}
where the remaining coefficients have been noted ${\tilde{\Psi}_0}$ and ${\tilde{\Psi}_2}$, also using the fact that dissipative monomials $a^2b$ and $b^3$ are vanishing in the case of a conservative system~\cite{touze03-NNM}.

For the sake of completeness and in order to show an important property of the real normal form styles (both RNF and FRNF), let us also introduce the mapping in polar coordinates. Using Eq.~\eqref{eq:cnf_reddyn_polar}, one can easily pass from complex to polar coordinates, such that the mapping in polar form reads:
\begin{equation}
    u = \rho \cos(\alpha) + 
    \dfrac{1}{4}\,\Psi_2\; \rho^3 \cos(\alpha)
    +
    \dfrac{1}{4}\,\Psi_0\; \rho^3 \cos(3\alpha),
    \label{killing}
\end{equation}
where the coefficients $\Psi_0$ and  $\Psi_0$ are the same as those introduced in Eq.~\eqref{eq:mapdufcomplexco}

The general formulation for the reduced dynamics up to the third-order and with  complex normal coordinates writes:
\begin{subequations}
\begin{align}
    &\dot{z} =  +\iu\omega_0\, z + 
    {{f_0}}\; z^3 + {{f_1}}\; z^2\bar{z} + 
    {{f_2}}\; z\bar{z}^2 + {{f_3}}\; \bar{z}^3,
    \\
    &\dot{z} =  -\iu\omega_0\, \bar{z} - 
    {{f_3}}\; z^3 - {{f_2}}\; z^2\bar{z} - 
    {{f_1}}\; z\bar{z}^2 - {{f_0}}\; \bar{z}^3,
\end{align}\end{subequations}
where the property of Eq.~\eqref{eq:symm_complex_f} has been used together with the fact that all $f_p$ coefficients are purely imaginary in the case of conservative systems.

If one uses now Cartesian coordinates, the general form of the reduced dynamics writes:
\begin{subequations}\label{eq:reduceddyancart}
\begin{align}
    &\dot{a} =  -\omega_0\,b +
    {{\tilde{f}_1}}\; a^2b + {{\tilde{f}_3}}\; b^3,\label{eq:reduceddyancarta}
    \\
    &\dot{b} =  +\omega_0\,a +
    {{\tilde{f}_0}}\; a^3 + {{\tilde{f}_2}}\; ab^2.\label{eq:reduceddyancartb}
\end{align}\end{subequations}

Now that the general form of the equations have been established for the simple case of a Duffing oscillator, let us specify the values of all the coefficients (mapping and reduced dynamics) depending on the different normal form style used: complex normal form CNF, as well as the two real normal styles: RNF and FRNF. The analytical values of the coefficients are reported in Tables~\ref{tab:compcoef} and~\ref{tab:cartesiancoefs}.

\begin{table}
\begin{center}
\begin{tabular}{|c|cc|cccc|}
\hline
\rule[-7pt]{0pt}{20pt}& $\Psi_0$ & $\Psi_2$ &
$f_0$ & $f_1$ & $f_2$ & $f_3$
\\\hline
\rule[-13pt]{0pt}{28pt} CNF 
    & $+\dfrac{\gamma}{8\omega_0^2}$ 
    & $-\dfrac{3\gamma}{4\omega_0^2}$ 
    & $0$
    & $\iu \dfrac{3\gamma}{2\omega_0}$
    & $0$
    & $0$
\\\hline
\rule[-13pt]{0pt}{28pt} RNF 
    & $+\dfrac{\gamma}{8\omega_0^2}$ 
    & $0$ 
    & $0$
    & $\iu \dfrac{3\gamma}{2\omega_0}$
    & $\iu \dfrac{3\gamma}{2\omega_0}$
    & $0$
\\\hline
\rule[-13pt]{0pt}{28pt} FRNF 
    & $0$ 
    & $0$ 
    & $\iu \dfrac{\gamma}{2\omega_0}$
    & $\iu \dfrac{3\gamma}{2\omega_0}$
    & $\iu \dfrac{3\gamma}{2\omega_0}$
    & $\iu \dfrac{\gamma}{2\omega_0}$
\\\hline
\end{tabular}
\end{center}
\caption{Analytical expression of the coefficients of the different normal styles in complex coordinates.}\label{tab:compcoef}
\end{table}

\begin{table}
\begin{center}
\begin{tabular}{|c|cc|cccc|}
\hline
\rule[-7pt]{0pt}{20pt}& $\tilde{\Psi}_0$ & $\tilde{\Psi}_2$ & $\tilde{f}_0$ & $\tilde{f}_1$ & $\tilde{f}_2$ & $\tilde{f}_3$
\\\hline
\rule[-13pt]{0pt}{28pt} CNF 
    & $-\dfrac{5\gamma}{32\omega_0^2}$ 
    & $-\dfrac{9\gamma}{32\omega_0^2}$ 
    & $+\dfrac{3\gamma}{8\omega_0}$ 
    & $-\dfrac{3\gamma}{8\omega_0}$ 
    & $+\dfrac{3\gamma}{8\omega_0}$ 
    & $-\dfrac{3\gamma}{8\omega_0}$ 
\\\hline
\rule[-13pt]{0pt}{28pt} RNF 
    & $+\dfrac{\gamma}{32\omega_0^2}$ 
    & $-\dfrac{3\gamma}{32\omega_0^2}$ 
    & $+\dfrac{3\gamma}{4\omega_0}$ 
    & $0$
    & $+\dfrac{3\gamma}{4\omega_0}$ 
    & $0$
\\\hline
\rule[-13pt]{0pt}{28pt} FRNF 
    & $0$ 
    & $0$ 
    & $+\dfrac{\gamma}{\omega_0}$ 
    & $0$
    & $0$
    & $0$
\\\hline
\end{tabular}
\end{center}
\caption{Analytical expression of the coefficients of the different normal styles in Cartesian coordinates.}\label{tab:cartesiancoefs}
\end{table}

From the values of the coefficients, one can draw the following interesting conclusions when comparing the different normal form styles:
\begin{itemize}
    \item Since $\tilde{f}_1=\tilde{f}_3=0$ for both real normal form styles, Eq.~\eqref{eq:reduceddyancarta} reduces to $\dot{a} = - \omega_0 \;b$. Thanks to this important simplification, one is then able to express the reduced dynamics as a single oscillator in second-order form without approximations, thus explaining the name of {\em real} normal form styles.
    \item Both real normal forms RNF and FRNF have $\Psi_2=0$. This has an important consequence which can be easily interpreted thanks to Eq.~\eqref{killing}: the amplitude of the fundamental harmonics does not depends on the nonlinearity and does not change with the asymptotic expansion. This important property has already been remarked for example in~\cite{NeildNF00,NeildNF01} for the RNF and in~\cite{VizzaMDNNM} for the FRNF style. It has also been named as {\em killing the fundamental} in~\cite{WaggReview} and it is a general property of the real normal form styles that contrasts with CNF.
    \item Since only $\tilde{f}_0$ is not vanishing in FRNF style with Cartesian coordinates (see Tab.~\ref{tab:cartesiancoefs}), this simply means that the FRNF of a Duffing oscillator is left unchanged, and all the mappings coefficients are simply vanishing. This is also a direct consequence of the interpretation of the resonance condition as Eq.~\eqref{eq:old_nf_cond}. As stated in the general comment, the formulation is fully real in this case. Importantly, at this order and without quadratic nonlinearity, FRNF is thus equivalent to the graph style in terms of reduced-order dynamics. This remark also explains why the backbone of the cantilever computed in~\cite{YichangVib} with DNF (direct normal form, which uses FRNF style) shows the same folding point as the one found here with the graph style. Importantly, using either CNF or RNF corrects this behaviour from the third-order and allows retrieving a better solution for the backbone of the cantilever.
\end{itemize}

\end{document}